\documentclass[a4paper,11pt]{amsart}
\usepackage{amssymb,amsmath,latexsym, upgreek,comment}
\usepackage{epsfig,caption}
\usepackage{mathrsfs}
\usepackage{hyperref}
\usepackage{color,graphics,subfigure}

\usepackage{epic,eepic,wrapfig,ifthen,enumerate}
\usepackage{url}
\newtheorem{thm}{Theorem}[section]
\newtheorem{corollary}[thm]{Corollary}
\newtheorem{lemma}[thm]{Lemma}
\newtheorem{prop}[thm]{Proposition}
\newtheorem{defn}[thm]{Definition}

\newtheorem{remark}[thm]{Remark}

\numberwithin{equation}{section}

\newcommand{\formula}[2][nolabel]
{\ifthenelse{\equal{#1}{nolabel}}
 {\begin{align*} #2 \end{align*}}
 {\ifthenelse{\equal{#1}{}}
  {\begin{align} #2 \end{align}}
  {\begin{align} \label{#1} #2 \end{align}}
 }
}

\def\pf{{\medskip\noindent {\bf Proof. }}}
\def\qed{{\hfill $\Box$ \bigskip}}

\renewcommand{\bar}{\overline}
\newcommand{\cal}[1]{\mathcal{#1}}
 \def\sB {{\cal B}} 
 \def\sE {{\cal E}} \def\sF {{\cal F}}
  \def\sI {{\cal I}}
  
\def\sM {{\cal M}} \def\sN {{\cal N}}

 \def\bH {{\mathbb H}}

\def\R {{\mathbb R}}

\def\Z {{\mathbb Z}}
\def\P{{\mathbb P}}
\def\E{{\mathbb E}}

\def\eps{\varepsilon}

\def\R{{\mathbb R}}

\def\E{{\mathbb E}}
\def\F{{\mathbf F}}
\def\P{{\mathbb P}}

\def\eps{\varepsilon}
\def\wh{\widehat}
\def\wt{\widetilde}
\def\pf{\noindent{\bf Proof.} }

\def\ub{{\overline{\beta}}}
\def\lb{{\underline{\beta}}}

\def\diam{{\text{\rm diam}}}
\def\1{{\bf 1}}

\def\nn{\nonumber}

\def\bn{{\bf n}}

\allowdisplaybreaks

\makeatletter
\@addtoreset{equation}{section}

\makeatother

\begin{document}
\bibliographystyle{plain}

\title[Heat kernel estimates for degenerate Markov processes in bounded sets]
{ \bf   
Heat kernel estimates for Markov processes in bounded sets with jump kernels decaying at the boundary}

\author{Soobin Cho, \quad Panki Kim, \quad Renming Song \quad and \quad Zoran Vondra\v{c}ek}

\address[Cho]{Department of Mathematics, University of Illinois Urbana-Champaign, Urbana, IL 61801, USA}
\curraddr{}
\email{soobinc@illinois.edu}

\address[Kim]{Department of Mathematical Sciences and Research Institute of Mathematics,
	Seoul National University,	Seoul 08826, Republic of Korea}
\thanks{Panki Kim is supported by the National Research Foundation of Korea(NRF) grant funded by the Korea government(MSIP) (No. RS-2023-00270314).}
\curraddr{}
\email{pkim@snu.ac.kr}

\address[Song]{
	Department of Mathematics, University of Illinois Urbana-Champaign, Urbana, IL 61801,
	USA}
\curraddr{}
\thanks{Research of Renming Song is supported in part by a grant from
	the Simons Foundation \#960480.}
\email{rsong@illinois.edu}

\address[Vondra\v{c}ek]
{
	Dr.~Franjo Tu{\dj}man Defense and Security 
	University
	 and Department of Mathematics, Faculty of Science, University of Zagreb, Zagreb, Croatia
}
\curraddr{}
\thanks{Zoran Vondra\v{c}ek is supported in part by the Croatian Science Foundation under the project IP-2025-02-8793.} 
\email{vondra@math.hr}

\date{}

\begin{abstract}
In this paper, we study two types of purely discontinuous symmetric Markov processes $X$ in bounded smooth subsets of $\R^d$: 
conservative processes and processes killed either  upon approaching the boundary of  the set or by a killing potential $\kappa$. 
The jump kernel of $X$ is of the form $J(x,y)=\sB(x,y)|x-y|^{-d-\alpha}$, $\alpha\in (0,2)$, where the function $\sB(x,y)$ decays to 0 at the boundary 
and is described in terms of two $O$-regularly varying functions and one slowly varying function. 
Under the conditions, introduced in \cite{CKSV24},  on $\sB(x,y)$ and on the killing potential $\kappa$, 
we establish sharp two-sided estimates on the heat kernel of $X$: in Lipschitz sets when $X$ is conservative, and in $C^{1,1}$ open sets for the killed process. 
\end{abstract}

\maketitle

\bigskip

\maketitle



\bigskip
\tableofcontents


\section{Introduction}\label{s:intro}

\subsection{Introduction and background}\label{ss:intro-back}
Let $D$ be an open subset of the Euclidean space $\R^d$.
The general theory of non-local singular operators of the type
\begin{equation}\label{e:operator-1}
	L^{\mathcal{B}}_{\alpha}f(x)=\lim_{\epsilon\to 0} \int_{D,\, |y-x|>\epsilon}\big(f(y)-f(x)\big) 
	\mathcal{B}(x,y)|x-y|^{-d-\alpha}\,dy,
\end{equation}
and
\begin{equation}\label{e:operator-2}
L f(x)=L^{\mathcal{B}}_{\alpha}f(x) - \kappa(x) f(x),	
\end{equation}
where the function $\mathcal{B}(x,y)$ vanishes at the boundary of $D$ and the killing potential $\kappa$  may be subcritical or critical, was very recently developed in \cite{CKSV24}.
The kernel of the operator is given by $J(x,y):=\sB(x,y)|x-y|^{-d-\alpha}$, and it decays to zero at the boundary of $D$.
In the case when $D$ is a $C^{1,1}$ open set, the main results obtained in \cite{CKSV24} 
include the boundary Harnack principle (and its possible failure), as well as sharp two-sided Green function estimates. 

The goal of the present paper is to essentially complete the study in \cite{CKSV24} by establishing sharp two-sided heat kernel estimates, or, in probabilistic terms, sharp two-sided estimates of the transition density of
the associated jump Markov process, in case $D$ is sufficiently smooth and bounded. 

Qualitative knowledge of heat kernels has numerous important consequences in partial differential equations, probability, functional analysis, and geometry. Among the many papers devoted to this topic, we mention \cite{BKKL19, BGR10, BGR14, CKK11, CKS10, CKS14, CK03, CK08,  CKW20a, CKW21, GHH18, GHH}. In all of these works,  the function $\sB(x,y)$ is assumed to be bounded between two positive constants, which can be viewed as a uniform ellipticity condition
for the operator $L^{\mathcal{B}}_{\alpha}$. As a consequence, 
the main contribution to the heat kernel at time $t$ comes from a single (large) jump from $x$ to $y$.

An example of a Markov process with a jump kernel $J(x,y)$ decaying at the boundary is a $\delta$-stable, 
$\delta\in (0, 2]$, 
process killed upon exiting an open set $D$, 
subsequently subordinated by an independent 
$\gamma$-stable subordinator, $\gamma\in (0, 1)$. 
The generator of the resulting process is $-((-\Delta)^{\delta/2}|_D)^\gamma$, which reduces to the spectral fractional Laplacian $-(-\Delta|_D)^\gamma$ when $\delta=2$. Spectral fractional Laplacians have been studied intensively in recent years in the PDE literature, see, for instance, \cite{BFV, BSV, BV, Gr} and the references therein.
Potential theory for this type of   
 Markov processes
was studied in \cite{KSV18-a, KSV18-b}. 
In the case $\delta \in (0, 2)$, there were two surprising results in \cite{KSV18-b}: 
(a) the jump kernel of the process exhibits a  phase transition at $\gamma=1/2$ and (b) the boundary Harnack principle holds for $\gamma\in(1/2,1)$ and fails for $\gamma\in (0, 1/2]$. 
The sharp two-sided heat kernel estimates 
established in \cite{CKSV22} revealed 
a somewhat anomalous qualitative behavior, quite 
different from that of the uniformly elliptic case. 

The surprising results 
in \cite{ CKSV22, KSV} motivated us to consider 
general operators of the form \eqref{e:operator-1} and \eqref{e:operator-2}. 
We first studied the case when $D=\R^d_+$, the upper half-space, and the corresponding Markov processes are self-similar, 
see \cite{KSV, KSV20, KSV21, CKSV23}.
The function $\sB(x,y)$ was assumed to depend on the distance between the points $x$ and $y$, as well as their distances to the boundary,
and on four parameters $\beta_1,\beta_2,\beta_3,\beta_4$ that roughly govern the strength of the boundary decay of $\sB$. 
The boundary Harnack principle (and its possible failure) together with sharp two-sided Green function estimates were established in \cite{KSV20}, while the more delicate heat kernel estimates were 
obtained recently in \cite{CKSV23}.
In the case when the killing potential $\kappa(\cdot)\equiv 0$ and the corresponding process is conservative,
it was shown that when both points $x$ and $y$ are close to the boundary, the usual “one big jump’’ form of the heat kernel estimates
 no longer holds. Moreover, depending on the relationship between 
the parameters $\beta_1$, $\beta_2$ and $\alpha$, 
it may happen that the dominant contribution to the heat kernel arises not from a single jump, but from two jumps. In the case 
when the corresponding process is not conservative,
the usual approximate factorization property was established, with a decay rate depending (in a rather intricate way) on the killing potential $\kappa$.

The present paper generalizes the results of \cite{CKSV23} in two directions. First, the upper half-space 
$\R^d_+$ is replaced by a bounded smooth set $D\subset \R^d$ (either Lipschitz or $C^{1,1}$). The main difficulty in this setting is the lack of scaling. The general framework addressing this problem was developed in \cite{CKSV24} and will be briefly reviewed below. Second, instead of power and logarithmic functions (governed by the parameters $\beta_1, \beta_2,\beta_3,\beta_4$), we consider 
functions satisfying appropriate weak scaling conditions, 
thereby allowing for some oscillations. The corresponding theoretical foundation was again established in \cite{CKSV24}.

Below we describe the setting and state the main results of the paper, first for the conservative case, and then for the killed case.

\subsection{Conservative case}\label{ss:conservative-case}

Let $D\subset \R^d$, $d\ge 2$, be a Lipschitz open set with localization radius $\wh{R}$ and Lipschitz constant $\Lambda_0$
(see Definition \ref{df:lipschitz}(i)
for a precise formulation).  Denote by  $\overline{D}$ the closure of $D$, and by $\delta_D(x)$ the Euclidean distance between  
$x$
and the boundary $\partial D$. 
We will assume  that  the jump measure of the process $\overline{Y}$, which we will construct, is absolutely continuous with respect to the Lebesgue measure on  $\overline{D}$. Since $D$ is Lipschitz, the Lebesgue measure of $\partial D$
is zero and the value of the jump kernel on $\partial D$ does not matter. 
Let $\alpha \in (0,2)$. We consider the symmetric bilinear form
\begin{equation*}
	\sE^0(u,v):= \frac{1}{2}\iint_{D \times D}  (u(x)-u(y))(v(x)-v(y)) \frac{\sB(x,y)}{|x-y|^{d+\alpha}}dxdy,
\end{equation*}
where $\sB:D\times D \to (0,\infty)$ is a symmetric   bounded Borel function. 
	
For a Borel set $A \subset \R^d$ and  $p \in [1,\infty]$, we denote by $L^p(A)$ the 
space $L^p(A, dx)$ and
by  $\mathrm{Lip}_c(A)$ the family of all Lipschitz functions 
of compact support in $A$. 
Note that, by the boundedness of $\sB$, $\sE^0(u,u)<\infty$ for all $u \in \mathrm{Lip}_c(\overline D)$. Let  $\overline \sF$ be the closure of 
$\mathrm{Lip}_c(\overline D)$ in 
$L^2(\overline D)=L^2(D)$ under  the norm $(\sE^0_1)^{1/2}$, where  $\sE^0_1:=\sE^0 + \lVert \cdot \rVert_{L^2(D)}^2$.  Then  $(\sE^0, \overline \sF)$ is a regular Dirichlet form  
on $L^2(\overline D)$ which is irreducible and conservative,  
see \cite[Chapter 1]{FOT} and \cite{CKSV24}. 
Associated with the regular Dirichlet form $(\sE^0, \overline \sF)$, there is a conservative Hunt process 
$\overline Y=(\overline Y_t,t \ge 0; \P_x,x \in \overline D\setminus  \sN)$. 
Here $ \sN\ $ is an  exceptional set for $\overline Y$. 

\medskip
 
We now recall the notion of  scaling conditions and Matuszewska indices.
\begin{defn}\label{def:Matuszewska}
\rm	Let $\Phi$ be a positive function on $(0,1]$.
	
	\noindent (i) We say that $\Phi$ satisfies the \textit{upper scaling condition with index $\overline \beta \in \R$ at zero}   if there exists a constant $C\ge 1$ such that
	\begin{align*}
		\Phi(r)/\Phi(s) \le C(r/s)^{\overline\beta} \quad \text{for all $0<s\le r\le 1$}.
	\end{align*}
	The \textit{upper Matuszewska index $\beta^*$ at  zero}  of $\Phi$ is the infimum of all $\overline\beta$ such that $\Phi$ satisfies the upper scaling condition with index $\overline\beta$.

		\noindent (ii) We say that $\Phi$ satisfies the \textit{lower scaling condition with index $\underline \beta \in \R$ at zero}   if there exists a constant $C\ge 1$ such that
	\begin{align*}
		\Phi(r)/\Phi(s) \ge C^{-1}(r/s)^{\underline\beta} \quad \text{for all $0<s\le r\le 1$}.
	\end{align*}
	The \textit{lower Matuszewska index $\beta$ at  zero}  of $\Phi$ is the supremum of all $\underline\beta$ such that $\Phi$ satisfies the lower scaling condition with index $\underline\beta$.
\end{defn}

Functions of the type described in Definition \ref{def:Matuszewska} are often referred to as  $O$-regularly varying functions; see, for example, \cite{BGT87}. 

A positive function $\Phi$ is said to be \textit{almost increasing} (resp. \textit{almost decreasing}) on the interval $(0,1]$, 
if there exists $C\ge 1$ such that $\Phi(r) \ge C^{-1} \sup_{s\in (0,r]}\Phi(s)$ (resp.  $\Phi(r) \le C \inf_{s\in (0,r]}\Phi(s)$) for all $r\in (0,1]$.
Note that $\Phi$ is almost increasing on $(0,1]$ if and only if $\Phi$ satisfies the lower scaling condition with index $\underline \beta\ge 0$, and  $\Phi$ is almost decreasing on $(0,1]$ if and only if $\Phi$ satisfies the upper scaling condition with 
index $\overline \beta\le 0$.  

Let  $\Phi_1$ and $\Phi_2$ be almost increasing positive Borel  functions on $(0,1]$ with the  lower and upper Matuszewska indices 
$(\beta_1, \beta_1^*)$ and $(\beta_2, \beta_2^*)$ respectively. Note that it holds that $ \beta_1^* \ge \beta_1\ge 0$ and   $ \beta_2^* \ge \beta_2 \ge0 $.
Let $\ell$ be a positive Borel  function on $(0,1]$ with the following property:  for any $\eps>0$, there exists a constant  $C=C(\eps)\ge1$ such that 
\begin{align}\label{e:ell-scaling-pre}
C^{-1} \bigg( \frac{r}{s}\bigg)^{ -\eps \wedge    \beta_1}
	\le \frac{\ell(r)}{\ell(s)} 
	\le C \bigg( \frac{r}{s}\bigg)^{  \eps  \wedge \beta_2} \quad \text{for all} \;\, 0<s\le r\le 1.			
\end{align}
Note that $\ell$  is almost increasing if $\beta_1=0$, and $\ell$ is almost decreasing if $\beta_2=0$.  We extend the functions $\Phi_1,\Phi_2,\ell$ to $(0,\infty)$ by setting
$\Phi_1(r)= \Phi_2(r)=\ell(r)=1$ for all $r \ge 1$.

We consider the following condition which coincides with the assumption \textbf{(B4-c)} in \cite{CKSV24}.
\medskip 

\setlength{\leftskip}{0.17in}

\noindent\hypertarget{A1}{{\bf (A1)}} There exist  comparison constants such that for all $x,y \in D$,
\begin{align*}
	\sB(x,y) \asymp  \Phi_1\bigg(\frac{\delta_D(x) \wedge \delta_D(y)}{|x-y|}\bigg)\,
	\Phi_2\bigg(\frac{\delta_D(x) \vee \delta_D(y)}{|x-y|}\bigg)\,
	\ell\bigg(\frac{\delta_D(x)\wedge \delta_D(y)}{(\delta_D(x)\vee\delta_D(y))\wedge |x-y|}\bigg).
\end{align*}

\setlength{\leftskip}{0in}

\medskip

Throughout this paper, we use the notation $a\wedge b:=\min\{a, b\}$, $a\vee b:=\max\{a, b\}$.  The notation $f\asymp g$ means that the quotient $f/g$ is bounded between two positive constants in a specified region.

We remark that \hyperlink{A1}{{\bf (A1)}}  implies conditions \textbf{(B2-a)}, \textbf{(B2-b)}, \textbf{(B4-a)} and \textbf{(B4-b)} from \cite{CKSV24}, cf.~\cite[Lemma 9.2]{CKSV24}. 
Throughout this paper, except in Section \ref{s:aux}, we  assume that $\sB$ satisfies \hyperlink{A1}{{\bf (A1)}}, where $\Phi_1$ and $\Phi_2$ are almost increasing functions with  lower and upper Matuszewska indices $(\beta_1, \beta_1^*)$ and $(\beta_2, \beta_2^*)$, respectively, and $\ell$ satisfies \eqref{e:ell-scaling-pre}. 
Then, by \cite[Proposition 4.2, Remark 4.9 and Lemma 9.2]{CKSV24},  the process $\overline Y$ 
can be refined to be a strongly Feller process which can start from every point in $\overline{D}$, and the exceptional set $\sN$ can be taken to be the empty set.  Moreover, $\overline Y$ 
has a jointly continuous transition density  $\overline p(t,x,y)$ defined on 	$(0,\infty) \times \overline D  \times \overline D$.

\bigskip

For positive functions $f,g,h$ on $(0, \infty)$, we define $A_{f, g, h}:[0,\infty) \times D\times D \to (0,\infty)$ by\begin{align}\label{e:def-A(f,g,h)}	\begin{split}	A_{f, g, h}(t,x,y)	&:=f\bigg(\frac{(\delta_D(x)\wedge \delta_D(y))\vee t^{1/\alpha}}{|x-y|}\bigg) \, g\bigg(\frac{\delta_D(x)\vee \delta_D(y)\vee t^{1/\alpha}}{|x-y|}\bigg)\\	&\quad\; \times h\bigg(\frac{(\delta_D(x)\wedge\delta_D(y))\vee t^{1/\alpha}}{(\delta_D(x)\vee \delta_D(y)\vee t^{1/\alpha})\wedge |x-y|}\bigg).\end{split}\end{align}Note that, by  \hyperlink{A1}{{\bf (A1)}},  $A_{\Phi_1, \Phi_2, \ell}(0,x,y)\asymp \sB(x,y)$ for $x,y \in D$.

For $x\in D$, we fix a point $Q_x\in \partial D$  
satisfying $|x-Q_x|=\delta_D(x)$. 
Define 
$\mathbf{n}_x:=(x-Q_x)/|x-Q_x|$. Since $D$ is Lipschitz, there exist  constants $\eta_1\in (0, 1]$ and $\eta_2\in (0,1)$  such that
\begin{equation}\label{e:lifting-property}
	\eta_2(\delta_D(x) + r) \le 	\delta_D(x + r \bn_x)  \le \delta_D(x) + r \quad \text{for all $x\in D$ and $r\in (0,\eta_1]$}.
\end{equation} 
In particular, we have $x+r\bn_x\in D$ for all $x\in D$ and $r\in (0,\eta_1]$. Define 
\begin{equation}\label{e:def-eps1}
	\eps_1:=\frac14 \wedge \frac{\eta_1}{\diam(D)}.
\end{equation}

The first main result of this paper in the conservative case is stated below.

\begin{thm}\label{t:main}
Suppose $D\subset \R^d$, $d\ge 2$, is a  bounded Lipschitz open set and  $\sB$ satisfies \hyperlink{A1}{{\bf (A1)}}. 
(i) If  $\beta_1^*<\alpha +  \beta_1$, then for any $T>0$,
there exist comparison constants depending on $T$ such that the following hold  for all $t\in (0,T]$ and $x,y\in \overline D$.
\begin{enumerate}[\indent (a)]
	\item (On-diagonal estimates)  If $t^{1/\alpha} \ge \eps_1|x-y|/2$, then
	\begin{align}\label{e:HKE-on}
		\overline p(t,x,y) \asymp t^{-d/\alpha}.
	\end{align}
	\item (Off-diagonal estimates) If $t^{1/\alpha} < \eps_1|x-y|/2$, then
	\begin{align}\label{e:HKE-off}
		\begin{split} 
			\overline{p}(t,x,y)\asymp   \frac{t^2}{|x-y|^{d+\alpha}}&\bigg[ \int_{ t^{1/\alpha}}^{\eps_1|x-y|}A_{\Phi_1,\Phi_2,\ell}(t,x,x+u\mathbf{n}_x)\,A_{\Phi_1,\Phi_2,\ell}(t,x+u\mathbf{n}_x,y)\frac{du}{u^{\alpha+1}}\\
			& \quad + \int_{ t^{1/\alpha}}^{\eps_1|x-y|}A_{\Phi_1,\Phi_2,\ell}(t,y,y+u\mathbf{n}_y)\,A_{\Phi_1,\Phi_2,\ell}(t,y+u\mathbf{n}_y,x)\frac{du}{u^{\alpha+1}} \bigg] .
		\end{split} 
	\end{align}
\end{enumerate}

\noindent 
(ii)  For all $t>T$ and $x,y\in D$,
$$
\overline{p}(t,x,y)\asymp 1.
$$
\end{thm}

\medskip

The assumption $\beta_1^*<\alpha +  \beta_1$ in Theorem \ref{t:main} means that the gap between the upper Matuszewska index $\beta_1^*$ and the lower Matuszewska index $\beta_1$ is strictly less than the parameter $\alpha$.
In Theorem \ref{t:main}, this assumption $\beta_1^*<\alpha +  \beta_1$ is used only in establishing the upper bound in \eqref{e:HKE-off}. The lower bounds in Theorem \ref{t:main} hold without the assumption $\beta_1^*<\alpha +  \beta_1$ --- this assumption is not used in the derivation of the lower bounds in Section \ref{s:lb}.
We believe this assumption is essential, rather than merely a limitation of our method for obtaining this type of estimate. 
{\it All subsequent main results will be proved under this assumption.}

Observe that the off-diagonal estimates \eqref{e:HKE-off} contain two terms with $x$ and $y$ interchanged. Without any additional assumptions, it is not possible to determine which term dominates.

Define \begin{equation}\label{e:Phi_0}\Phi_0(r):=\Phi_1(r)\ell(r), \quad r>0. \end{equation}

 \begin{remark}\label{r:1.3}{\rm 
	 (i)  Using the simple properties \eqref{e:upper-bound-A} and \eqref{e:interior-lower-bound-A}
	 of the function $A_{\Phi_1,\Phi_2,\ell}$,
	  one can combine \eqref{e:HKE-on} and \eqref{e:HKE-off} as follows: 
	 for all $t\in (0,T]$ and $x,y \in \overline D$,
	\begin{align}\label{e:HKE-off-2}
		\begin{split} 
			 \overline{p}(t,x,y)&
			\asymp    \left( t^{-d/\alpha} \wedge \frac{t}{|x-y|^{d+\alpha}}\right)   \\
			& \times   \bigg[ \big(t\wedge  |x-y|^{\alpha} \big)  \int_{t^{1/\alpha} \wedge (\eps_1|x-y|/2)}^{\eps_1|x-y|}A_{\Phi_1,\Phi_2,\ell}(t,x,x+u\mathbf{n}_x)\,A_{\Phi_1,\Phi_2,\ell}(t,x+u\mathbf{n}_x,y)\frac{du}{u^{\alpha+1}}\\
			& \quad+   \big(t\wedge  |x-y|^{\alpha} \big)  \int_{t^{1/\alpha} \wedge (\eps_1|x-y|/2)}^{\eps_1|x-y|}A_{\Phi_1,\Phi_2,\ell}(t,y,y+u\mathbf{n}_y)\,A_{\Phi_1,\Phi_2,\ell}(t,y+u\mathbf{n}_y,x)\frac{du}{u^{\alpha+1}} \bigg] .
		\end{split} 
	\end{align}	
	
	\noindent  (ii) 	
	We have the following  alternative 
	form 
	of \eqref{e:HKE-off} 
	(see Lemma \ref{l:equivalences}):
 For all $t\in (0,T]$ and $x,y\in \overline D$ 
	with $r:=|x-y|$ satisfying $t^{1/\alpha}<\eps_1r/2$,
	\begin{align*}
&	\overline{p}(t,x,y)	\asymp   
 \frac{t}{r^{d+\alpha}} 
\bigg[  A_{\Phi_1,\Phi_2,\ell}(t,x,y)\\
		& +  t\Phi_1\bigg(\frac{\delta_D(x)\vee t^{1/\alpha}}{r}\bigg)  
		\int_{(\delta_D(x)\vee \delta_D(y)\vee t^{1/\alpha})\wedge r}^{r}\Phi_0\bigg(\frac{\delta_D(y)\vee t^{1/\alpha}}{u} \bigg)   \Phi_2\bigg(\frac{u}{r}\bigg) \ell \bigg(\frac{\delta_D(x)\vee t^{1/\alpha}}{u}\bigg) \frac{du}{u^{\alpha+1}}\\
		& +  t \Phi_1\bigg(\frac{\delta_D(y)\vee t^{1/\alpha}}{r}\bigg) 
		\int_{(\delta_D(x)\vee \delta_D(y) \vee  t^{1/\alpha})\wedge r}^{r}\Phi_0\bigg(\frac{\delta_D(x)\vee t^{1/\alpha}}{u} \bigg)   \Phi_2\bigg(\frac{u}{r}\bigg) \ell \bigg(\frac{\delta_D(y)\vee t^{1/\alpha}}{u}\bigg) \frac{du}{u^{\alpha+1}}\bigg].
	\end{align*}
	\noindent (iii) Using \hyperlink{A1}{{\bf (A1)}} and the scaling properties of $\Phi_1,\Phi_2$ and $\ell$,  
	we can obtain
	another equivalent form 
	of  \eqref{e:HKE-off}: Let $J(x,y):=\sB(x,y)|x-y|^{-d-\alpha }$ 
	be the jump kernel of $\sE^0$. 	For all $t\in (0,T]$ and $x,y\in \overline D$ satisfying $t^{1/\alpha}<\eps_1|x-y|/2$,
	\begin{align*}
		\begin{split} 
			&\overline p(t,x,y)\asymp  t^2 \bigg(  \int_{ t^{1/\alpha} }^{\eps_1|x-y|} J(x+ t^{1/\alpha} \bn_x, x+ u\bn_x)\,J(x+ u\bn_x,y + t^{1/\alpha} \bn_y) \, u^{d-1}du \\
			&\qquad\qquad \qquad \;+  \int_{ t^{1/\alpha}}^{\eps_1|x-y|} J(y + t^{1/\alpha}\bn_y,y+u \bn_y)\,J(y+ u\bn_y,x + t^{1/\alpha}\bn_x) \, u^{d-1}du
			\bigg) .
		\end{split}
	\end{align*}
}
\end{remark}

\medskip

The next theorem is a refinement of the previous one, in which we distinguish two cases: (i) $ \beta^*_2 <\alpha+\beta_1$, 
and (ii) $ \beta_2 >\alpha+\beta_1^*$.  Clearly, these two cases do not 
exhaust all possible relationships among the parameters 
$\beta_1$, $\beta_1^*$, $\beta_2$, $\beta_2^*$ and $\alpha$.  
In each of these two cases, we can give simpler forms of the estimates for $\overline{p}(t,x,y)$.
In case (i), the main contribution to the heat kernel comes from a single (big) jump, while in the case (ii) the main contribution may come from two jumps. Note that the first case corresponds to \cite[Theorem 1.2(i)]{CKSV23}, while the second corresponds to \cite[Theorem 1.2(ii)]{CKSV23}.
Part (iii) of the theorem below deals with a borderline case.

\begin{thm}\label{t:main11}
 Suppose $D\subset \R^d$, $d\ge 2$, is a  bounded Lipschitz open set and  $\sB$ satisfies \hyperlink{A1}{{\bf (A1)}} with $\beta_1^*<\alpha +  \beta_1$. 
Then for any $T>0$, there exist comparison constants depending on $T$ such that  the following estimates hold  for all $t\in (0,T]$ and $x,y\in \overline D$:  

 	\noindent (i)  When  $ \beta^*_2 <\alpha+\beta_1$, 
 	\begin{align}\label{e:main11-case1}
 		&\overline p(t,x,y)\asymp   \left( t^{-d/\alpha} \wedge \frac{t}{|x-y|^{d+\alpha}}\right)  A_{\Phi_1, \Phi_2,  \ell} (t,x,y).
 	\end{align}
 	
 	\noindent (ii)		When $\beta_2>\alpha+\beta_1^*$, 
	 \begin{align}\label{e:main11-case2}
 		\overline p(t,x,y) \asymp 	 \left(t^{-d/\alpha}\wedge \frac{t}{|x-y|^{d+\alpha}}\right) 
 		\left[ A_{\Phi_1, \Phi_2, \ell}(t,x,y)+ \left(1\wedge \frac{t}{|x-y|^{\alpha}}\right)A_{\Phi_0, \Phi_0,1}(t,x,y)\right] .
	\end{align}
	\noindent (iii)	 
	When $\Phi_1(r)=(1\wedge r)^{\beta_1}$ and $\Phi_2(r)=(1\wedge r)^{\beta_2}\phi(r)$, with $\beta_1\ge 0$, $\beta_2=\alpha+\beta_1$ and $\phi$ being a positive function on $(0, 1]$ with  lower and
upper Matuszewska indices  both equal to  $0$, 	
	\begin{align}\label{e:main11-case3}
 		&\overline p(t,x,y)\asymp
 		\left(t^{-d/\alpha}\wedge \frac{t}{|x-y|^{d+\alpha}}\right)\\
 		& \times  \bigg[ A_{\Phi_1, \Phi_2,\ell}(t,x,y)
 		+ \left(1\wedge \frac{t}{|x-y|^{\alpha}}\right)	 A_{\Phi_0,\Phi_0,1}(t,x,y)\int_{(\delta_D(x)\vee \delta_D(y)\vee t^{1/\alpha})\wedge |x-y|}^{|x-y|}\frac{du}{u}\nn\\
 			&\hspace{2in} \times  \bigg( \frac{\ell((\delta_D(x) \vee t^{1/\alpha})/u) \,\ell((\delta_D(y) \vee t^{1/\alpha})/u)\, \phi(u/|x-y|)}{\ell((\delta_D(x) \vee t^{1/\alpha})/|x-y|)\, \ell((\delta_D(y) \vee t^{1/\alpha})/|x-y|) }  \bigg) 	\bigg]. \nn 	
 	\end{align}
 \end{thm}

 \medskip

In general, neither of the  two terms in the brackets on the right-hand side of \eqref{e:main11-case2} dominates the other. Also, neither of the  two terms in the brackets on the right hand side of \eqref{e:main11-case3} dominates the other in general; see Lemma \ref{l:two-jumps-non-dominant}.
 	
 \begin{remark}\label{r:1.5}
 Note that the right-hand sides of both \eqref{e:main11-case2} and \eqref{e:main11-case3} are comparable to $t^{-d/\alpha}$ if $t^{1/\alpha}\ge \eps_1|x-y|/2$. Using \hyperlink{A1}{{\bf (A1)}} and the scaling properties of $\Phi_1,\Phi_2$ and $\ell$,  
we can obtain the following alternative formulations of the off-diagonal estimates in \eqref{e:main11-case1} and \eqref{e:main11-case2}
 (cf. Remark \ref{r:1.3}(iii)):   	 	When $\beta_2^*<\alpha+\beta_1$, for all $t\in (0,T]$ and $x,y\in \overline D$ satisfying $t^{1/\alpha}<\eps_1|x-y|/2$,
 \begin{align*}
 	\overline p(t,x,y) \asymp tJ(x+ t^{1/\alpha}\bn_x, y + t^{1/\alpha} \bn_y).
 \end{align*}
 When $\beta_2>\alpha+\beta_1^*$, for all $t\in (0,T]$ and $x,y\in \overline D$ satisfying $t^{1/\alpha}<\eps_1|x-y|/2$,
 	 	\begin{align*}
 	 		\overline p(t,x,y)& \asymp  \Big [ tJ(x+ t^{1/\alpha}\bn_x, y + t^{1/\alpha} \bn_y)\\
 	 		&\qquad  + t^2 |x-y|^d J(x+ t^{1/\alpha}\bn_x, x+ \eps_1|x-y| \bn_x)\, J(x+ \eps_1 |x-y|\bn_x, y + t^{1/\alpha} \bn_y)  \Big] \\
 	 		& \asymp  \Big [ tJ(x+ t^{1/\alpha}\bn_x, y + t^{1/\alpha} \bn_y)\\
 	 		&\qquad  + t^2\int_{B(x+ \frac{3\eps_1}{4}|x-y| \bn_x ,\frac{\eps_1}{4} |x-y| ) } J(x+ t^{1/\alpha}\bn_x, z)\, J(z, y + t^{1/\alpha} \bn_y) \,dz \Big].
 	 	\end{align*}
 \end{remark}

 \begin{remark}
 	Under the setting of Theorem \ref{t:main11}(iii), suppose in addition that 
 	\begin{align*}
 		\ell(r) = \frac{1}{\log^{\beta_3}(e+1)} \log^{\beta_3} \bigg(e+ \frac{1}{r\wedge 1}\bigg) \quad \text{and} \quad \phi(r) = \frac{1}{\log^{\beta_4}(e+1)}  \log^{\beta_4} \bigg( e  + \frac{1}{r\wedge 1}\bigg)
 	\end{align*}
 	with $\beta_3,\beta_4\ge 0$, and that $\beta_1>0$ if $\beta_3>0$. Then by \cite[Lemma A.9]{CKSV23}, the right-hand side of \eqref{e:main11-case3} is comparable to
 	\begin{align*}
 			&
 		\left(t^{-d/\alpha}\wedge \frac{t}{|x-y|^{d+\alpha}}\right) \bigg[ A_{\Phi_1, \Phi_2,\ell}(t,x,y)\\
 		& \quad 
 		+ \left(1\wedge \frac{t}{|x-y|^{\alpha}}\right)	 A_{\Phi_0,\Phi_0,1}(t,x,y) \log^{\beta_4+1}\bigg(e+\frac{|x-y|}{(\delta_D(x)\vee \delta_D(y)\vee t^{1/\alpha})\wedge |x-y|}\bigg) \bigg].
 	\end{align*}
 	This recovers the form of  the half-space heat kernel estimates established  in \cite[Theorem 1.2(iii)]{CKSV23}.
 \end{remark}  

\subsection{Killed case}\label{ss:killed-case}

We now discuss the case in which the process $\overline{Y}$ is killed either upon approaching the boundary $\partial D$ or by an appropriate killing potential $\kappa$.

Let $\sF^0$ be the closure of  $\mathrm{Lip}_c(D)$ in $L^2(D)$ under $\sE^0_1$. Then $(\sE^0,\sF^0)$ is a regular Dirichlet form. Let $Y^0=(Y^0_t,t \ge 0; \P_x,x \in D \setminus \sN_0)$ be the Hunt process associated with $(\sE^0,\sF^0)$, where  $\sN_0$ is an  
exceptional set for $Y^0$. Since $Y^0$ is a subprocess of $\overline{Y}$, we can take $\sN_0$ to be the empty set. By \cite[Proposition 4.4]{CKSV24}, $\sF^0=\overline \sF$ if and only if $\alpha\le 1$. 
In case $\alpha\in (1,2)$, the process $Y^0$ is equal to the process $\overline{Y}$ killed upon approaching the boundary $\partial D$, while for $\alpha\in (0,1]$, $Y^0=\overline{Y}$ and never hits the boundary.

Another process we are going to study
is obtained by killing $Y^0$ via  a killing potential $\kappa$. Let $\kappa$ be a non-negative locally integrable Borel function on $D$. The symmetric form $(\sE^\kappa, \sF^\kappa)$ is defined by
	\begin{align*}
		\sE^\kappa(u,v)&=\sE^0(u,v) + 
		\int_{D} u(x)v(x)\kappa(x)dx, \\
		\sF^\kappa&= 	\wt \sF^0 \cap L^2(D, 
			\kappa(x) dx),\nn
	\end{align*}
	where $\wt \sF^0$ is the family of all $\sE^0_1$-quasi-continuous functions in $\sF^0$. Then $(\sE^\kappa,\sF^\kappa)$ is a regular Dirichlet form on $L^2(D)$ with Lip$_c(D)$ as a special standard core,  see \cite[Theorems 6.1.1 and 6.1.2]{FOT}.
Let  $Y^\kappa=(Y^\kappa_t,t \ge 0; \P_x,x \in D \setminus  \sN_\kappa)$ be the Hunt process associated with 
$(\sE^\kappa, \sF^\kappa)$ where $\sN_\kappa$ is an exceptional set for $Y^\kappa$. We denote by $\zeta^\kappa$ the lifetime of $Y^\kappa$, and define $Y^\kappa_t=\partial$ for $t \ge \zeta^\kappa$, where $\partial$ is a cemetery point added to the state space $D$. 
Note that $Y^\kappa$ includes $Y^0$, when $\alpha\in (1, 2)$, as a special case. The process $Y^\kappa$ can be regarded as the part process of $\overline Y$ killed at $\zeta^\kappa$. Hence, $Y^\kappa$ can be refined to be a Hunt process 
which can start from every point in $D$.

We now specify the assumptions used in this subsection. Let $D\subset \R^d$ be a $C^{1,1}$ open set with characteristics $(\wh{R}, \Lambda)$; see Definition 
\ref{df:lipschitz}(ii) for a precise formulation. Without loss of generality, we assume that   $\wh R\le 1 \wedge (1/(2\Lambda))$. 
The first condition concerns  the killing potential $\kappa$ and coincides with assumption \textbf{(K3)} in \cite{CKSV24}.

\setlength{\leftskip}{0.17in}

\medskip
\noindent  \hypertarget{K}{{\bf (K)}} There exist constants $\kappa_0\ge 0$,   $\eps_0 >0$ and $C_1>0$	such that for all $x \in D$,
	\begin{align*}	\begin{cases}			|\kappa(x) - \kappa_0 \sB(x,x) \delta_D(x)^{-\alpha}| \le C_1\delta_D(x)^{-\alpha+\eps_0} &\text{ if }  \delta_D(x) < 1,\\[4pt]			\kappa(x) \le C_1 &\text{ if } \delta_D(x) \ge 1.		\end{cases}	\end{align*} 		When $\alpha \le 1$, we further assume that $\kappa_0>0$. 
	
\medskip
\setlength{\leftskip}{0in}	
Note that condition \hyperlink{K}{{\bf (K)}} states that in the case when $\kappa_0>0$, the killing potential is \emph{critical} in the sense that near the boundary the main contribution comes from the term $\kappa_0\sB(x,x)\delta_D(x)^{-\alpha}$. 
Assuming $\kappa_0>0$ when $\alpha \le 1$ ensures that, in this case, the killing potential $\kappa$  
is non-trivial.
 
We will explain in Subsection \ref{ss:setup} that, 
under an additional assumption -- \hyperlink{A2}{{\bf (A2)}}  -- 
there exists a strictly increasing and continuous function 
$C(\cdot\, ; \alpha, \sB):[(\alpha-1)_+, \alpha+\beta_1)\to [0,\infty)$ such that 
$\kappa_0=C(q; \alpha, \sB)$ for a unique $q$. 
That is, to every $\kappa_0\ge 0$ we can associate a 
parameter $q\in [(\alpha-1)_+, \alpha+\beta_1)$. The parameter $q$ is a fundamental quantity, as it governs the decay rate of certain harmonic functions near the boundary (see Theorem \ref{t:Dynkin-improve}). 
	
Assumptions \hyperlink{A2}{{\bf (A2)}} and \hyperlink{A3}{{\bf (A3)}} are introduced precisely in Subsection 
\ref{ss:setup}  and  coincide, respectively,  with \textbf{(B5)} and \textbf{(B3)} in \cite{CKSV24}. Condition \hyperlink{A2}{{\bf (A2)}} is quite technical and may be viewed as a substitute for the flattening of the boundary 
method (which does not work in the current setting). 
The interested reader is referred to \cite[Section 1]{CKSV24} for a detailed explanation of this assumption. We note that, in this paper, we do not use these two assumptions explicitly, but only their consequences -- most notably, Theorem \ref{t:Dynkin-improve} and Proposition \ref{p:bound-for-integral-new}.

For examples of processes that satisfy all the assumptions \hyperlink{A1}{{\bf (A1)}}, \hyperlink{A2}{{\bf (A2)}}, \hyperlink{A3}{{\bf (A3)}} and \hyperlink{K}{{\bf (K)}}, see \cite[Section 11]{CKSV24}.

Suppose that assumptions  \hyperlink{A1}{{\bf (A1)}}, \hyperlink{A2}{{\bf (A2)}}, \hyperlink{A3}{{\bf (A3)}} and 
\hyperlink{K}{{\bf (K)}} are satisfied. As a consequence, all the assumptions from \cite{CKSV24} are valid.
Then, by \cite[Proposition 4.15, Remark 4.19  and Lemma 9.2]{CKSV24},  the process $Y^{\kappa}$ has a jointly continuous transition density  $p^{\kappa}(t,x,y)$ defined on $(0,\infty) \times  D  \times D$.

We can now state the third main result of the paper. 

\begin{thm}\label{t:main2}
 Suppose $D\subset \R^d$, $d\ge 2$, is a bounded $C^{1,1}$ open set,   $\kappa$ satisfies \hyperlink{K}{{\bf (K)}}, and $\sB$ satisfies \hyperlink{A1}{{\bf (A1)}}, \hyperlink{A2}{{\bf (A2)}} and  \hyperlink{A3}{{\bf (A3)}} with $\beta_1^*<\alpha +  \beta_1$.  Then for any $T>0$, there exist comparison constants depending on $T$ such that the following estimates hold.

  \noindent (i) For all $t\in (0,T]$ and $x,y\in D$,
 \begin{align*}
 	p^\kappa(t,x,y)
 	& \asymp  \left( 1 \wedge \frac{\delta_D(x)}{t^{1/\alpha}}\right)^q \left( 1 \wedge \frac{\delta_D(y)}{t^{1/\alpha}}\right)^q \overline p(t,x,y) \asymp \P_x(\zeta^\kappa>t)  \P_y(\zeta^\kappa>t) \overline p(t,x,y),
 \end{align*}
  where $q \in [(\alpha-1)_+,\alpha+\beta_1)$ is the  positive constant satisfying \eqref{e:C(alpha,p,F)}.
  
 \noindent (ii)  
 For all $t\ge T$ and $x,y\in D$,
  \begin{align*}
 	p^\kappa(t,x,y)
 	& \asymp   e^{-\lambda_1 t}\, \delta_D (x)^{q}\, \delta_D (y)^{q} \asymp  e^{-\lambda_1 t}\, \P_x(\zeta^\kappa>1)  \P_y(\zeta^\kappa>1),
\end{align*}
where $-\lambda_1<0$ is the largest eigenvalue of the infinitesimal generator of $Y^\kappa$.
\end{thm}

\subsection{Organization of the paper and roadmap of the proofs}\label{ss:org}

Section \ref{s:preliminary} contains several preliminary results and is divided into two subsections. 
In Subsection \ref{ss:setup} we give precise formulations of the assumptions \hyperlink{A2}{{\bf (A2)}} and \hyperlink{A3}{{\bf (A3)}}. For the motivation behind these assumptions and further explanation, we refer the reader to \cite{CKSV24}. Subsection \ref{ss:CKSV24} summarizes a number of useful results established in \cite{CKSV24} under assumptions \hyperlink{A2}{{\bf (A2)}} and \hyperlink{A3}{{\bf (A3)}}. In addition to some preliminary estimates for the heat kernel, we single out Theorem \ref{t:Dynkin-improve} and Proposition \ref{p:bound-for-integral-new}, which will play a central role in determining the exact decay rate of the heat kernel.

Section \ref{s:AA} provides several inequalities for integrals of functions $A_{\Phi_1,\Phi_2,\ell}$ of the type appearing in \eqref{e:HKE-off} and \eqref{e:HKE-off-2}. The main results are Lemmas \ref{l:two-jumps-lower-bound}–\ref{l:two-jumps-comparability}. Lemmas \ref{l:two-jumps-lower-bound} and \ref{l:two-jumps-comparability} compare expressions involving these integrals—first with $A_{\Phi_1,\Phi_2,\ell}$, and then with $A_{\Phi_0,\Phi_0,1}$. The former provides two lower bounds without any assumptions on the Matuszewska indices, while the latter derives two upper bounds depending on the relations between the indices. Lemma \ref{l:equivalences} gives an off-diagonal comparison of three types of expressions involving the function $A_{\Phi_1,\Phi_2,\ell}$ and its integrals. The results in this section depend only on the assumptions on the functions $\Phi_1$, $\Phi_2$, and $\ell$, together with the geometric properties of the set $D$; they do not rely on \hyperlink{A1}{{\bf (A1)}}, \hyperlink{A2}{{\bf (A2)}}, \hyperlink{A3}{{\bf (A3)}}, or \hyperlink{K}{{\bf (K)}}.

Sections \ref{s:pue}--\ref{s:proofs} are devoted to proving the sharp two-sided heat kernel estimates. 
The large-time estimates are fairly straightforward, and a short proof is given in Section \ref{s:proofs}. 
The main goal of this paper is to establish sharp two-sided small-time estimates. This is a  formidable task due to the singularity of  the function $\mathcal{B}(x,y)$ at the boundary, and requires very delicate arguments.
Note that, in the conservative case, the small-time heat kernel estimates  consist of two parts, cf.~\eqref{e:HKE-off-2}: 
the usual $t^{-d/\alpha}\wedge (t|x-y|^{-d-\alpha})$ estimate for the $\alpha$-stable process, and the contribution involving 
the function $A_{\Phi_1,\Phi_2, \ell}$, which arises from the function $\sB(x,y)$. 
In the killed case there is also a third component providing the boundary decay; cf.~Theorem \ref{t:main2}(i). 
In this paper,  these two cases
are treated simultaneously. The conservative case does not need 
\hyperlink{K}{{\bf (K)}}, \hyperlink{A2}{{\bf (A2)}} and  \hyperlink{A3}{{\bf (A3)}},  and the state space is a bounded Lipschitz.
The techniques employed to establish the lower bound differ from -- and are somewhat simpler than -- those used for the upper bound. Furthermore, the derivation of the lower bounds does not rely on the assumption that $\beta_1^*<\alpha+\beta_1$.

In Section \ref{s:pue} we prove a preliminary small-time upper bound on the heat kernel. The main result is Proposition \ref{p:UHK-rough} 
which sharpens the estimates of Lemma \ref{l:(5.6)implies(5.1)} by incorporating, in the killed case,  
the boundary decay terms. 
 We first show that this upper bound follows from a weaker inequality involving only one of the two boundary decay terms, see Lemma \ref{l:(5.6)implies(5.1)}. 
The proof of this weaker inequality proceeds in two steps:  
the first and simpler step
establishes the estimate for $q\in (0,\alpha)$, while the second employs an induction argument to extend it to the full range of the parameter $q$. A consequence of Proposition \ref{p:upper-heatkernel} is an upper bound of the lifetime probability, see Corollary \ref{c:life}.

Section \ref{s:lb} is devoted to the proof of the sharp lower bound. The main technical work is carried out in Subsection \ref{ss:lifted-points}, where various estimates involving the lifted points (introduced in \eqref{e:def-x(t)}) are derived. 
Subsection \ref{ss:slb} refines these estimates and establishes 
the full small-time heat kernel lower estimates. 
Corollary \ref{c:life-2} complements Corollary \ref{c:life} and provides two-sided survival probability estimates. 

Sharp upper bounds, which constitute the most difficult part of the proof, are addressed in Section \ref{s:ub}. 
The first step is to introduce the function $A_{f, 1,1}$ in the upper bound. In Lemma \ref{l:UHK-case1-induction-pre} this is achieved for $f(r)=(r\wedge 1)^{\lb}$ where $\lb\in ((q-\alpha)_+, \beta_1)$ if $\beta_1>0$ and $\lb=0$ if $\beta_1=0$. 
The proof proceeds by induction, which gradually  
increases
the exponent in the function $r\mapsto r^{\beta}$ until it reaches $\lb$.
The induction basis is provided by Proposition \ref{p:UHK-rough}. Lemma \ref{l:UHK-case1-induction} provides a similar estimate with the function $A_{\Phi_0,1,1}$. 
Again, the proof proceeds by induction with Lemma \ref{l:UHK-case1-induction-pre} serving as the induction basis. In the induction step we crucially use Lemma \ref{l:analog-of-10.10}, which explains the necessity of the assumption  $\beta_1^*<\alpha+\beta_1$. 
The next step in the proof uses Lemma \ref{l:general-upper-2} and separately estimates the two terms appearing in that lemma. The  estimate of the first term is simpler and is given in Lemma \ref{l:firstpart}. 
Estimating the second term, which involves the triple integral, is technically the most demanding part of the proof. It is carried out in Lemmas \ref{l:UHK-case1-main1} -- \ref{l:newlemma1-case1}, 
that carefully estimate the involved integrals. 
The sharp upper bound is given in Theorem \ref{t:UHK}.

The short Section \ref{s:proofs} completes the proofs of Theorems \ref{t:main}, \ref{t:main11}, and \ref{t:main2}.

The paper ends with a section on auxiliary results concerning the estimates of various integrals involving 
functions satisfying weak scaling conditions.  
This section is independent of the assumptions used throughout the rest of the paper.

We end this introduction with a few words on the additional notation. For $r\in \R$, we use the notation $\lfloor r \rfloor:=\max\{m\in \Z: m\le r\}$. For $x_0\in \overline{D}$ and $r>0$, 
we denote $B_{\overline D}(x_0, r)=B(x_0,r)\cap \overline{D}$ and $B_{D}(x_0, r)=B(x_0,r)\cap {D}$.


\section{Conditions {\bf (A2)}-{\bf (A3)}  and preliminary results}\label{s:preliminary}
In this section, we first give  the precise setup of our paper, including the formulations
of the assumptions \hyperlink{A2}{{\bf (A2)}} and \hyperlink{A3}{{\bf (A3)}}, and then recall several important results, 
based on  \hyperlink{A2}{{\bf (A2)}} and \hyperlink{A3}{{\bf (A3)}}, obtained in \cite{CKSV24}.

\subsection{Conditions {\bf (A2)}-{\bf (A3)}}\label{ss:setup}

For   $\gamma^*\ge \gamma$, we use $\sM(\gamma,\gamma^*)$  to denote the family of all positive functions on $(0,1]$ with  lower and upper Matuszewska  indices 
$(\gamma, \gamma^*)$,  and
\begin{align*}
	\sM^\uparrow(\gamma,\gamma^*):= \left\{ f\in \sM(\gamma,\gamma^*): \text{$f$ is almost increasing on $(0,1]$}\right\}.
\end{align*}
Note that $\sM^\uparrow(\gamma,\gamma^*)=\sM(\gamma,\gamma^*)$ if $\gamma>0$, $\sM^\uparrow(\gamma,\gamma^*)=\emptyset$ if $\gamma<0$, and $\sM^\uparrow(0,\gamma^*)\subsetneq \sM(0,\gamma^*)$ for all $\gamma^*\ge 0$.  We extend every $f\in \sM(\gamma,\gamma^*)$ to $(0,\infty)$ by setting $f(r)=1$ for all $r\ge 1$.

We  recall the notions of Lipschitz open sets and  $C^{1,1}$ open sets.
\begin{defn}\label{df:lipschitz}{\rm
			Let   $D\subset \R^d$ be an open set.

		\noindent (i)	We say that $D$ is a Lipschitz open set  with \textit{localization radius}  $\wh R>0$
		and \textit{Lipschitz constant} $\Lambda_0>0$,
		if for any $Q \in \partial D$, there exist a Lipschitz function $\Psi=\Psi^Q:\R^{d-1}\to \R$ with 
		\begin{align*}
			\Psi(\wt 0)=0 \quad \text{and} \quad |\Psi(\wt y)-\Psi(\wt z)| \le \Lambda_0|\wt y - \wt z| \; \text{ for all } \,  \wt y, \wt z \in \R^{d-1},
		\end{align*} 
		and an orthonormal coordinate system CS$_{Q}$ with origin at $Q$ such that
		\begin{align}\label{e:local-coordinate}
			B_D(Q, \wh R)=
			\left\{ y= (\wt y, y_d) \in B(0, \wh R) \text{ in CS$_{Q}$} : y_d>\Psi(\wt y) \right\}.
		\end{align}
		
		\noindent (ii) We say that  $D$	is a $C^{1,1}$ open set  		with characteristics  $(\wh R, \Lambda)$, if 		for any $Q \in \partial D$, there exist a $C^{1,1}$ function $\Psi=\Psi^Q:\R^{d-1} \to \R$ with 
		\begin{align*}		\Psi(\wt 0)= |\nabla \Psi(\wt 0)|=0 \quad \text{and} \quad |\nabla \Psi(\wt y)-\nabla\Psi(\wt z)| \le \Lambda |\wt y-\wt z| \; \text{ for all } \,  \wt y, \wt z \in \R^{d-1},		\end{align*}
			and an orthonormal coordinate system CS$_{Q}$ with origin at $Q$ such that	\eqref{e:local-coordinate} holds. 
	}
\end{defn}

\medskip
Define $\bH_{-1}:=\{(\wt y, y_d) \in \R^d: y_d>-1\}$. For $Q\in \partial D$, $\nu\in (0,1]$ and $r\in (0, \wh R/4]$, we let
\begin{align}\label{e:EQ}
	E^Q_\nu(r):= \left\{ y=(\wt y, y_d) \text{ in CS$_{Q}$}:  |\wt y|<r/4, \,    4r^{-\nu}|\wt y|^{1+\nu} <y_d < r/2 \right\}.
\end{align}
 
Depending on whether $\kappa_0>0$ or $\kappa_0=0$, 
we will assume different additional conditions on $\sB$.

\medskip

\textbf{Case $\kappa_0>0$:} In this case we consider the following condition which coincides with assumption \textbf{(B5-I)} in \cite{CKSV24}. 

\medskip

\setlength{\leftskip}{0.17in}

\noindent \hypertarget{A2-I}{{\bf (A2-I)}} There  exist  constants $\nu \in (0,1]$, $\theta_1,\theta_2,C_2>0$, and a non-negative Borel function $\F_0$ on $\bH_{-1}$   such that  for any $Q \in \partial D$ and    $x,y \in E^Q_\nu(\wh R/8)$ with $x=(\wt x,x_d)$ in CS$_Q$,
\begin{align*}
	& \big|\sB(x,y)- \sB(x,x)\F_0((y-x)/x_d) \big| + \big|\sB(x,y)- \sB(y,y)\F_0((y-x)/x_d) \big|\nn\\	&\le C_2\bigg(\frac{ \delta_D(x)\vee\delta_D(y) \vee |x-y|}{ \delta_D(x) \wedge \delta_D(y) \wedge |x-y|} \bigg)^{\theta_1}
	\big( \delta_D(x)\vee\delta_D(y) \vee|x-y|\big)^{\theta_2}.
\end{align*}

\setlength{\leftskip}{0in}
\medskip

Under condition \hyperlink{A2-I}{{\bf (A2-I)}}, we  define a  function $\F$ on $\bH_{-1}$ by 
\begin{align}\label{e:def-F0-transform}
	\F(y) =\frac{\F_0(y) + \F_0(-y/(1+y_d))}{2}, \quad \;\;y  = (\wt y,y_d) \in \bH_{-1}.
\end{align}
By \cite[Lemma 6.2]{CKSV24}, $\F$ is a bounded function. Moreover, we 	have
\begin{align}\label{e:F-symmetrization}
	\F(y)=\F(-y/(1+y_d)) \quad \text{for all} \;\, y \in \bH_{-1}.
\end{align}
With the function $\F$ in \eqref{e:F-symmetrization} and 
$p \in [(\alpha-1)_+, \alpha+\beta_1)$, we associate a  constant $C(p;\alpha,\F)$ defined by
\begin{align}\label{e:def-killing-constant}
	&	C(p;\alpha,f)=\int_{\R^{d-1}}\frac{1}{(|\wt u|^2+1)^{(d+\alpha)/2}}\int_0^1  	\frac{(s^p -1)(1-s^{\alpha-1-p})}{(1-s)^{1+\alpha}} f\big(((s-1)\wt u, s-1)\big) ds \, d \wt u. 
\end{align}
Assume additionally that 	
\begin{align}\label{e:K-upper-bound}	
	\kappa_0< \lim_{p \to \alpha+\beta_1} C(p;\alpha,\F).
\end{align}

It is shown in \cite[Lemma 6.3]{CKSV24} that	
$p \mapsto C(p;\alpha,\F)$ 
is a well-defined strictly increasing continuous function 
on $[(\alpha-1)_+, \alpha+\beta_1)$ 
and $C((\alpha-1)_+; \alpha,\F)=0$.  Therefore, under \eqref{e:K-upper-bound},  there exists a unique constant 
$q\in ((\alpha-1)_+, \alpha + \beta_1)$ such that 
\begin{align}\label{e:C(alpha,p,F)}
	\kappa_0=C(q;\alpha,\F).
\end{align}

\medskip

\textbf{Case $\kappa_0=0$:} In this case, instead of \hyperlink{A2-I}{{\bf (A2-I)}}, we introduce the following weaker condition  which coincides with assumption \textbf{(B5-II)} in  
\cite{CKSV24}:

\medskip

\setlength{\leftskip}{0.17in}

\noindent \hypertarget{A2-II}{{\bf (A2-II)}} There  exist  constants 
$\nu \in (0,1]$, $\theta_1,\theta_2,C_{2}>0$, $C_3>1$, $i_0 \ge 1$, and  non-negative  Borel  functions 
$\F_0^i:\bH_{-1} \to [0,\infty)$ and $\mu^i:D \to (0,\infty)$, $1\le i \le i_0$, such that 
\begin{align*}
	C_3^{-1} \le	\mu^i(x) \le C_3 \quad \text{for all} \;\, x \in D,
\end{align*}
and for any $Q \in \partial D$ and    $x,y \in E^Q_\nu(\wh R/8)$ with $x=(\wt x,x_d)$ in CS$_Q$,
\begin{align*}
	\begin{split}
		&\bigg|\sB(x,y)- \sum_{i=1}^{i_0}\mu^i(x)\F_0^i((y-x)/x_d) \bigg| 
		+  \bigg|\sB(x,y)- \sum_{i=1}^{i_0}\mu^i(y)\F_0^i((y-x)/x_d) \bigg| \\
		&\le C_2 \bigg(\frac{ \delta_D(x)\vee\delta_D(y) \vee |x-y|}{ \delta_D(x) \wedge \delta_D(y) \wedge |x-y|} \bigg)^{\theta_1}
		\big( \delta_D(x)\vee\delta_D(y) \vee|x-y|\big)^{\theta_2}.
	\end{split}
\end{align*}

\setlength{\leftskip}{0in}

Note that if \hyperlink{A2-I}{{\bf (A2-I)}} holds, then  \hyperlink{A2-II}{{\bf (A2-II)}} holds with $i_0=1$, $\F_0^1=\F_0$ and 
$\mu^1(x)=\sB(x,x)$.

\medskip

The conditions \hyperlink{A2-I}{{\bf (A2-I)}} and \hyperlink{A2-II}{{\bf (A2-II)}} always come in connection with \hyperlink{K}{{\bf (K)}}. We now combine these two conditions into 
condition \hyperlink{A2}{{\bf (A2)}}  which coincides with assumption \textbf{(B5)} in \cite{CKSV24}.

\medskip

\setlength{\leftskip}{0.17in}

\noindent \hypertarget{A2}{{\bf (A2)}}  If $\kappa_0>0$, then	\hyperlink{A2-I}{{\bf (A2-I)}} and \eqref{e:K-upper-bound}	 hold, and  if $\kappa_0=0$, then \hyperlink{A2-II}{{\bf (A2-II)}} holds.

\setlength{\leftskip}{0in}

\bigskip

If $\kappa_0=0$, then by 	\hyperlink{K}{{\bf (K)}}, we have $\alpha>1$, and 
\eqref{e:C(alpha,p,F)} holds with $q=\alpha-1$. Throughout the paper, if \hyperlink{A2}{{\bf (A2)}} is in force, we let $q\in [(\alpha-1)_+,\alpha+\beta_1)$ be the positive constant satisfying \eqref{e:C(alpha,p,F)}.

\medskip

Finally, we state the assumption \hyperlink{A3}{{\bf (A3)}} which coincides with \textbf{(B3)} in \cite{CKSV24}.

\setlength{\leftskip}{0.17in}

\medskip

\noindent	\hypertarget{A3}{{\bf (A3)}} If $\alpha \ge 1$, then there exist constants $\theta_0>\alpha-1$ and $C_4>0$ such that
\begin{align*}
	|\sB(x,x)-\sB(x,y)| \le C_4 \bigg(\frac{|x-y|}{\delta_D(x) \wedge \delta_D(y) \wedge  1
	} \bigg)^{\theta_0} 
	\quad \text{for all }  x,y \in D.
\end{align*}

\setlength{\leftskip}{0in}

\medskip

\emph{From now on,  we  	   assume  that $D\subset \R^d$, $d\ge 2$,  is a  bounded  Lipschitz open set, and except in Section \ref{s:AA}, 	 we 	   assume  that  $\sB$ satisfies \hyperlink{A1}{{\bf (A1)}} with $\Phi_1,\Phi_2$ and $\ell$,  where $\Phi_1 \in \sM^\uparrow(\beta_1,\beta_1^*)$ and $\Phi_2\in \sM^\uparrow(\beta_2,\beta_2^*)$,  and $\ell \in \sM(0,0)$ satisfies \eqref{e:ell-scaling-pre}. Throughout the paper,  whenever we work with $Y^\kappa$, we 	additionally assume that $D$ is a $C^{1,1}$ open set,   $\kappa$ satisfies \hyperlink{K}{{\bf (K)}} and $\sB$ satisfies \hyperlink{A2}{{\bf (A2)}} and \hyperlink{A3}{{\bf (A3)}}.
}

\medskip

\subsection{Preliminary results}\label{ss:CKSV24}

We define a function $\Phi_0$ on $(0, \infty)$ by
\begin{align}\label{e:def-Phi0}
	\Phi_0(r):=\Phi_1(r)\ell(r), \quad r>0.
\end{align}
From the definitions of   Matuszewska indices, using \eqref{e:ell-scaling-pre}, and the facts that $\Phi_1$ and $\Phi_2$ are almost increasing and $\Phi_1(r)=\Phi_2(r) = \ell(r)=1$ for $r\ge 1$, we see that for any $R\ge 1$, $\eta \in (0,1]$ and $\eps>0$,   there exists $C= C(R,\eta, \eps)\ge 1$ such that
\begin{align}
	C^{-1}\bigg( \frac{r}{s}\bigg)^{(\beta_1-\eps)_+}
	&\le \frac{\Phi_0(r)}{\Phi_0(s)}\le C
	\bigg( \frac{r}{s}\bigg)^{\beta_1^*+\eps} \quad \text{for all $0<\eta s\le r\le R$,} \label{e:Phi0-scaling}\\
	C^{-1}\bigg( \frac{r}{s}\bigg)^{(\beta_1-\eps)_+}
	&\le \frac{\Phi_1(r)}{\Phi_1(s)}\le C
	\bigg( \frac{r}{s}\bigg)^{\beta_1^*+\eps}\quad \text{for all $0<\eta s\le r\le R$,}\label{e:Phi1-scaling}\\
	C^{-1}\bigg( \frac{r}{s}\bigg)^{(\beta_2-\eps)_+}
	&\le \frac{\Phi_2(r)}{\Phi_2(s)}\le C
	\bigg( \frac{r}{s}\bigg)^{\beta_2^*+\eps}\quad \text{for all $0<\eta s\le r\le R$,}\label{e:Phi2-scaling}\\
	C^{-1} \bigg( \frac{r}{s}\bigg)^{ -\eps \wedge    \beta_1}
	&	\le \frac{\ell(r)}{\ell(s)} 
	\le C \bigg( \frac{r}{s}\bigg)^{  \eps  \wedge \beta_2}\quad\;\;\; \text{for all $0<\eta s\le r$.}\label{e:ell-scaling}
\end{align}
In particular, $\Phi_0 \in \sM^\uparrow(\beta_1,\beta_1^*)$. Moreover, 
since $\sB$ satisfies \hyperlink{A1}{{\bf (A1)}}, 
by \cite[Lemma 9.2]{CKSV24}, there exists $C>0$ such that
\begin{align}\label{e:B4-a}
	\sB(x,y) \le C \Phi_0 \left( \frac{\delta_D(x)\wedge \delta_D(y)}{|x-y|}\right) \quad \text{for all $x,y \in D$}.
\end{align} 

 For a relatively open subset $U$ of $\overline D$, define  $\overline \tau_U:=\inf\{t>0: \overline Y_t \notin U\}$  and  $ \tau^\kappa_U:=\inf\{t>0:  Y^\kappa_t \notin U\cap D\}$. Since $Y^\kappa$ is a subprocess of $\overline Y$, we have $\tau^\kappa_U\le \overline \tau_U$ for any relatively open subset $U$ of $\overline D$.  

The following sharp boundary estimates are among the most important intermediate results of \cite{CKSV24}, 
and they crucially depend on  assumptions \hyperlink{A2}{{\bf (A2)}} and \hyperlink{A3}{{\bf (A3)}}. 

Recall that, for $x_0\in \overline{D}$ and $r>0$, $B_{\overline D}(x_0, r)=B(x_0,r)\cap \overline{D}$ 
and $B_{D}(x_0, r)=B(x_0,r)\cap {D}$.

\begin{thm}\label{t:Dynkin-improve}	
	Suppose that $D\subset\R^d$, $d\ge 2$, is a bounded $C^{1,1}$ open set, $\sB$ satisfies \hyperlink{A1}{{\bf (A1)}}, \hyperlink{A2}{{\bf (A2)}} and \hyperlink{A3}{{\bf (A3)}}, and  $\kappa$ satisfies \hyperlink{K}{{\bf (K)}}. Then	there exist constants $\eta_0\in (0,1/36]$ and $C\ge 1$  such that for any $Q\in \partial D$,  $r\in (0, \eta_0 \wh R]$ 
	and  $x \in B_D(Q,2^{-4} r)$,
	\begin{align*}
		C^{-1}  \bigg(\frac{\delta_D(x)}{r}\bigg)^q\le 	\P_x \big(	Y^\kappa_{\tau^\kappa_{B_D(Q,r)}} 	\in D \big)   \le C   
		\bigg(\frac{\delta_D(x)}{r}\bigg)^q.
	\end{align*}
\end{thm}
\pf Let $Q\in \partial D$. For 
$r \in (0, \wh{R}/2)$, define 
$$
U(r)=\big\{x\in D: x=(\wt{x},x_d)\text{ in CS}_Q \mbox{ with } |\wt{x}|< r,\, 0<\rho_D(x)<r\big\}.
$$	
Here $\rho_D(x):=x_d-\Psi(\wt{x})$, where $\Psi$ is the $C^{1,1}$ function satisfying
\eqref{e:local-coordinate}. 
By  \cite[Theorem 7.4]{CKSV24}, there exist $a_0\in (0,1/24]$ and $c_1\ge 1$ such that for any $r\in (0,a_0 \wh R]$ and $x\in U(r/4)$,
\begin{align}\label{e:Dynkin-improve}	
	c_1^{-1} (\delta_D(x)/r)^q\le	\P_x \big(	Y^\kappa_{\tau^\kappa_{U( r)}}\in D 		\big) 	 \le c_1 (\delta_D(x)/r)^q.	\end{align}

By \cite[(3.15)]{CKSV24}, we have $B_D(Q,2r/3) \subset U(r) \subset B_D(Q,2r)$ for all $r\in (0,a_0\wh R]$. It follows that for all $r\in (0,a_0 \wh R]$ and $x\in B_D(Q,2^{-4}r)$,
\begin{align*}
	\P_x \big(	Y^\kappa_{\tau^\kappa_{U( 3r/2)}}\in D 		\big) &\le 	\P_x \big(	Y^\kappa_{\tau^\kappa_{B_D(Q, r)}}\in D 		\big) \le 	\P_x \big(	Y^\kappa_{\tau^\kappa_{U( r/2)}}\in D 		\big) .
\end{align*}
Thus, \eqref{e:Dynkin-improve}	 yields the desired result with $\eta_0:= 2a_0/3$. \qed

In  the remainder of the paper, we let $\eta_0$ denote the constant from Theorem \ref{t:Dynkin-improve}.

Using Theorem \ref{t:Dynkin-improve}, 
the following result was established in \cite[Proposition 8.6]{CKSV24}.

\begin{prop}\label{p:bound-for-integral-new} 
	Let $Q \in \partial D$ and $\gamma>q-\alpha$.  Then there exists $C\ge 1$ such that  for any $R \in (0, \wh R/24]$,  any Borel set $V$ satisfying 	$B_D(Q,R/4) \subset V \subset B_D(Q,R)$  and any $x \in B_D(Q,R/8)$,
	\begin{equation*} 	
		C^{-1} 	\delta_D(x)^q R^{\alpha+  \gamma   -q}\le 	\E_x \int_0^{\tau_V}\delta_D(		Y^\kappa_t		)^{\gamma }\, dt  \le C 	\delta_D(x)^q R^{\alpha+  \gamma   -q} .
	\end{equation*}
\end{prop}

We recall several results from \cite{CKSV24} that will be used later in this paper. We note that these remain valid under assumptions weaker than those stated in the last paragraph of Subsection \ref{ss:setup}; see \cite{CKSV24} for details.

Under  assumptions \hyperlink{A1}{{\bf (A1)}}, \hyperlink{A2}{{\bf (A2)}} and  \hyperlink{K}{{\bf (K)}}, the process $Y^\kappa$ has a jointly continuous transition density  $p^\kappa(t,x,y)$ defined on 	$(0,\infty) \times D  \times D$; see \cite[Remark 4.19]{CKSV24}. We extend $p^\kappa(t,x,y)$ to $(0,\infty) \times \overline D\times \overline D$ by setting $p^\kappa(t,x,y)=0$ if $x\in \partial D$ or $y\in \partial D$. Recall that $\overline Y$  has a jointly continuous  transition density $\overline p(t,x,y)$
on 	$(0,\infty) \times \overline D  \times \overline D$.

\begin{prop}\label{p:upper-heatkernel} \cite[Propositions 4.2 and 4.15]{CKSV24}
	For any $T>0$, there exists a constant $C=C(T)>0$ such that
	\begin{equation*}
		p^\kappa(t,x,y)\le 		\overline	p(t,x,y) \le C \left( t^{-d/\alpha} \wedge \frac{t}{|x-y|^{d+\alpha}}\right) \quad \text{for all $t\in (0,T]$ and $x,y \in \overline D$}.
	\end{equation*}
\end{prop}

\begin{prop}\label{p:E}
	\cite[Propositions 4.6 and 4.17]{CKSV24} (i)	For any $R_0>0$,  there exists $C=C(R_0)>1$ such that
	\begin{align*}
		C^{-1}r^\alpha \le 	\E_{x_0}[\bar\tau_{	B_{\overline D}(x_0,r)}] \le C r^\alpha 
		\quad \text{for all} \;\, x_0 \in \overline D, \; 0<r\le R_0.
	\end{align*}
	\noindent (ii) For any $R_0>0$, 	there exists $C=C(R_0)>1$ such that	\begin{align*}	C^{-1}r^\alpha \le 	\E_{x_0}[\tau^\kappa_{B(x_0,r)}] \le C r^\alpha \quad \text{for all} \;\, x_0 \in  D, \; 0<r<\delta_D(x_0) \wedge R_0.	\end{align*}
\end{prop}

For a relatively open subset $U$ of $\overline{D}$, we denote
by $\overline{Y}^{U}$  the part of $\overline{Y}$ killed upon exiting $U$, and for any  open set $U\subset {D}$, we denote
by ${Y}^{\kappa, U}$  the part of  $Y^\kappa$ killed upon exiting $U$. Then $\overline Y^U$ and  $Y^{\kappa, U}$ have jointly continuous transition densities $\overline p^{U}(t,x,y)$ and $p^{U}(t,x,y)$  with respect to the Lebesgue measure on $U$, respectively. See \cite[Remarks 4.9 and 4.19]{CKSV24}.

\begin{prop}\label{p:ndl}
	\cite[Propositions 4.5 and 4.16]{CKSV24}	Let $R_0>0$ and  $b\in(0,1)$.
	
	\noindent (i) There exists  $C=C(R_0,b)>0$  such that for any $x_0 \in \overline D$, $0<r\le R_0$ and 
	$0< t\le (b r)^\alpha$, it holds that
	\begin{equation*}
		\overline	p^{	B_{\overline D}(x_0, r)}(t,x,y)  \ge  C t^{-d/\alpha} 
		\quad  \text{for all } x,y \in 	B_{\overline D}(x_0, b t^{1/\alpha}).
	\end{equation*}
	
	\noindent (ii) There exists  $C=C(R_0,b)>0$ such that for any $x_0 \in D$, $0<r<\delta_D(x_0) \wedge R_0$ and  $0< t\le (b r)^\alpha$, it holds that
	\begin{equation*}
		p^{\kappa,B(x_0, r)}(t,x,y)  \ge  C t^{-d/\alpha}\quad  \text{for all } x,y \in 	B(x_0, b t^{1/\alpha}). 
	\end{equation*}
\end{prop}

Since the form $\sE^0$ has jump  measure $\sB(x,y)|x-y|^{-d-\alpha} dxdy$, 
the process $\overline Y$ satisfies  the following  L\'evy system formula on $\overline D\times \overline D$: for  any stopping time $\sigma$ and   any  Borel function $f: \overline D \times  \overline D \to [0,\infty]$ vanishing on the set $\{(x,x):x \in \overline D\}$, the following  holds for all $x\in \overline D$:
\begin{align}
	\E_x\sum_{s\in (0, \sigma]}f(\overline Y_{s-}, \overline Y_s)= \E_x\int^{\sigma}_0\int_{\overline D}f(\overline Y_s, y) \frac{\sB( \overline Y_s, y)}{|\overline Y_s-y|^{d+\alpha}} dy ds.\label{e:levy-system-Y-bar}
\end{align}
Similarly, since the pure-jump part of $\sE^\kappa$ is given by $\sE^0$, the process $Y^\kappa$ satisfies  the following  L\'evy system formula on $ D\times  D$: for  any stopping time $\sigma$ and
any  Borel function $f: (D\cup \{\partial\})\times (D\cup \{\partial\}) \to [0,\infty]$ vanishing on the set $\{(x,x):x \in D\} \cup (\{\partial\}\times D) \cup ( D\times \{\partial\})$, the following  holds for all $x\in D$:
\begin{align}
	\E_x\sum_{s\in (0, \sigma]}f(Y^\kappa_{s-}, Y^\kappa_s)= \E_x\int^{\sigma}_0\int_{D}f(Y^\kappa_s, y) \frac{\sB(Y^\kappa_s, y)}{|Y^\kappa_s-y|^{d+\alpha}} dy ds.\label{e:levy-system-Y-kappa}
\end{align}

The following lemma is a  consequence of the L\'evy system formulas \eqref{e:levy-system-Y-bar} and \eqref{e:levy-system-Y-kappa}, together with the joint continuity and symmetry of the heat kernels; see, e.g., \cite[Lemma 3.15]{CKSV23}. We omit the proof.

\begin{lemma}\label{l:general-upper-2}	
	(i)	Let  $V_1$ and $V_3$ be open subsets of $\overline{D}$ with ${\rm dist}(V_1,V_3)>0$. 
	Set $V_2:=\overline{D} \setminus (V_1 \cup V_3)$. For any $x\in V_1$, $y \in V_3$ and $t>0$, it holds that 	
	\begin{align*}	\bar p(t,x,y) &\le  \P_x(\overline \tau_{V_1}<t) \sup_{s \le t, \, z \in V_2} \bar p(s,z,y)\\			
		& \quad  + {\rm dist}(V_1,V_3)^{-d-\alpha}\int_0^t\int_{V_3} \int_{V_1} \bar p^{V_1}(s, x, u) \sB(u,w) \bar p(t-s, y,w) du dw ds.	
	\end{align*}
	(ii)	Let  $V_1$ and $V_3$ be open subsets of $D$ with ${\rm dist}(V_1,V_3)>0$. Set $V_2:=D \setminus (V_1 \cup V_3)$. 
	For any $x\in V_1$, $y \in V_3$ and $t>0$, it holds that 	
	\begin{align*}	p^\kappa(t,x,y) &\le  \P_x(\tau^\kappa_{V_1}<t< \zeta^\kappa) \sup_{s \le t, \, z \in V_2} p^\kappa(s,z,y)\\			
		&\quad  + {\rm dist}(V_1,V_3)^{-d-\alpha}\int_0^t\int_{V_3} \int_{V_1} p^{\kappa,V_1}(s, x, u) 
		\sB(u,w) p^\kappa(t-s, y,w) du dw ds.		
	\end{align*}
\end{lemma}


\section{Results related to $A_{\Phi_1,\Phi_2,\ell}$}\label{s:AA}

Throughout this section, we assume that $D$ is a \emph{bounded} Lipschitz open set. 
   Recall that,  for  positive functions $f,g,h$, the function $A_{f, g, h}(t,x,y)$ is defined in \eqref{e:def-A(f,g,h)}.
The purpose of this 
section is to establish various inequalities involving the function $A$ and its integrals. The obtained results depend only on the assumptions on the functions $\Phi_1$, $\Phi_2$ and $\ell$,  and geometry of $D$, but not on the assumptions \hyperlink{A1}{{\bf (A1)}}, \hyperlink{A2}{{\bf (A2)}}, \hyperlink{A3}{{\bf (A3)}} and \hyperlink{K}{{\bf (K)}}.

 By using \eqref{e:Phi1-scaling}--\eqref{e:ell-scaling}, 
in the same way as in 
\cite[Lemma 9.2(i)]{CKSV24},  we can show that 
\begin{align}\label{e:upper-bound-A}
	&A_{\Phi_1,\Phi_2,\ell}(t,x,y)\le C \quad \text{for all $t\ge 0$ and $x,y \in D$,} 
\end{align}
and that for any $a \in (0,1)$, there exists a constant $c(a)>0$ such that 
\begin{equation}\label{e:interior-lower-bound-A}
	A_{\Phi_1,\Phi_2,\ell}(t,x, y)\ge c(a)
\end{equation}for all $t\ge 0$ and $x,y\in D$ with  $(\delta_D(x)\wedge\delta_D(y))\vee t^{1/\alpha}\ge a|x-y|$. 

We now give several elementary but useful inequalities for integrals of functions $A_{\Phi_1,\Phi_2, \ell}$. Their proofs are technical,  so the reader may wish to skip the proofs  at the first reading.
 
Recall that for $x\in D$, the point $Q_x\in \partial D$ satisfies  $|x-Q_x|=\delta_D(x)$ and $\mathbf{n}_x=(x-Q_x)/|x-Q_x|$. 
We also recall  that constants $\eta_1$ and $\eta_2$ were introduced in \eqref{e:lifting-property}, $\eps_1$ was defined in \eqref{e:def-eps1},
and the function $\Phi_0$ in \eqref{e:Phi_0}.

\begin{lemma}\label{l:two-jumps-equivalent}
	There exist comparison constants such that  for $(t,x,y)\in (0,\infty) \times D\times D$,
	\begin{align} \label{e:two-jumps-equivalent}
		\begin{split} 
			& \int_{(\delta_D(x) \vee \delta_D(y)\vee t^{1/\alpha}) \wedge (\eps_1|x-y|/2)}^{\eps_1|x-y|}A_{\Phi_1,\Phi_2,\ell}(t,x,x+u\mathbf{n}_x)\,A_{\Phi_1,\Phi_2,\ell}(t,x+u\mathbf{n}_x,y)\frac{du}{u^{\alpha+1}}\\
			&\asymp  \Phi_1\bigg(\frac{\delta_D(y)\vee t^{1/\alpha}}{|x-y|}\bigg) \\
			&\quad \times  \int_{(\delta_D(x)\vee \delta_D(y)\vee t^{1/\alpha})\wedge (|x-y|/4)}^{|x-y|}\Phi_0\bigg(\frac{\delta_D(x)\vee t^{1/\alpha}}{u} \bigg)   \Phi_2\bigg(\frac{u}{|x-y|}\bigg) \ell \bigg(\frac{\delta_D(y)\vee t^{1/\alpha}}{u}\bigg) \frac{du}{u^{\alpha+1}}.
		\end{split} 
	\end{align}
\end{lemma}
\pf 
For $u \in ((\delta_D(x)\vee \delta_D(y) \vee t^{1/\alpha}) \wedge (\eps_1|x-y|/2),\eps_1|x-y|)$, we have
$u\le \eta_1$, and hence $\delta_D(x+u\bn_x)\asymp u \vee \delta_D(x)$ and 
\begin{align}\label{e:two-jumps-equivalent-1}
	A_{\Phi_1,\Phi_2,\ell}(t,x,x+u\mathbf{n}_x)&\asymp \Phi_1\bigg( \frac{\delta_D(x) \vee t^{1/\alpha}}{u}\bigg)  \ell\bigg( \frac{\delta_D(x)\vee t^{1/\alpha}}{ u} \bigg)=\Phi_0\bigg( \frac{\delta_D(x) \vee t^{1/\alpha}}{u}\bigg).
\end{align}
We claim that, for   $u \in ((\delta_D(x)\vee \delta_D(y) \vee t^{1/\alpha}) \wedge (\eps_1|x-y|/2),\eps_1|x-y|)$,
\begin{align}\label{e:two-jumps-equivalent-2}
	A_{\Phi_1,\Phi_2,\ell}(t,x+u\mathbf{n}_x,y)\asymp \Phi_1\bigg( \frac{\delta_D(y) \vee t^{1/\alpha}}{|x-y|}\bigg) \Phi_2\bigg( \frac{u }{|x-y|}\bigg) \ell\bigg( \frac{\delta_D(y) \vee t^{1/\alpha}}{  u} \bigg).
\end{align}
First note that, for   $u \in ((\delta_D(x)\vee \delta_D(y) \vee t^{1/\alpha}) \wedge (\eps_1|x-y|/2),\eps_1|x-y|)$,  we have $(1-\eps_1)|x-y|<|x+u\bn_x-y| < (1+\eps_1)|x-y|$. 

When $\delta_D(x) \vee \delta_D(y)\vee t^{1/\alpha}  < \eps_1|x-y|/2$,  since $\delta_D(x)\vee \delta_D(y) \vee t^{1/\alpha}  \le u \le \eps_1|x-y|$, we have 
$\delta_D(y) \vee t^{1/\alpha} \le  u   \wedge (\eps_1|x-y|) \asymp u   \wedge |x+u\bn_x-y|$ and  $\delta_D(x+u\mathbf{n}_x) \asymp u \vee \delta_D(x) \asymp u \ge \delta_D(y) \vee t^{1/\alpha}$ by \eqref{e:lifting-property}. 
Thus, using \eqref{e:Phi1-scaling}, \eqref{e:Phi2-scaling} and \eqref{e:ell-scaling},  \eqref{e:two-jumps-equivalent-2} holds true.

When $\delta_D(x) \vee \delta_D(y)\vee t^{1/\alpha}   \ge  \eps_1|x-y|/2$, we have  
$\delta_D(x+u\mathbf{n}_x) \vee \delta_D(y)\vee t^{1/\alpha}   \ge  \eta_2\eps_1|x-y|/2$ by \eqref{e:lifting-property}.
Moreover, in this case,  $ u \asymp |x-y|$ and so  $\delta_D(y) \le \delta_D(x) +|x-y| \asymp \delta_D(x) +u \asymp \delta_D(x+u\mathbf{n}_x)$ by \eqref{e:lifting-property}. 
Thus, using  \eqref{e:Phi1-scaling}, \eqref{e:Phi2-scaling} and \eqref{e:ell-scaling},  \eqref{e:two-jumps-equivalent-2} holds true in this case as well.

Note that, by \eqref{e:Phi0-scaling},  \eqref{e:Phi2-scaling} and \eqref{e:ell-scaling}, we have
\begin{align*} 
	&	\int_{\eps_1|x-y|}^{|x-y|}\Phi_0\bigg( \frac{\delta_D(x) \vee t^{1/\alpha}}{u}\bigg)  \Phi_2\bigg( \frac{u }{|x-y|}\bigg) \ell\bigg( \frac{\delta_D(y) \vee t^{1/\alpha}}{u} \bigg) \frac{du}{u^{\alpha+1}} \\
	&\le c_1 \int_{\eps_1|x-y|}^{|x-y|}\Phi_0\bigg( \frac{\delta_D(x) \vee t^{1/\alpha}}{|x-y|}\bigg)  \Phi_2(1) \ell\bigg( \frac{\delta_D(y) \vee t^{1/\alpha}}{|x-y|} \bigg) \frac{du}{|x-y|^{\alpha+1}} \\
	&\le \frac{c_1}{|x-y|^{\alpha}}\Phi_0\bigg( \frac{\delta_D(x) \vee t^{1/\alpha}}{|x-y|}\bigg)   \ell\bigg( \frac{\delta_D(y) \vee t^{1/\alpha}}{|x-y|} \bigg)\\
	&\le  c_2\int_{\eps_1|x-y|/2}^{\eps_1|x-y|}\Phi_0\bigg( \frac{\delta_D(x) \vee t^{1/\alpha}}{u}\bigg) \Phi_2\bigg( \frac{u }{|x-y|}\bigg) \ell\bigg( \frac{\delta_D(y) \vee t^{1/\alpha}}{u} \bigg) \frac{du}{u^{\alpha+1}} 
\end{align*}
and
\begin{align*} 
	&	\int_{\eps_1|x-y|/2}^{|x-y|/4}\Phi_0\bigg( \frac{\delta_D(x) \vee t^{1/\alpha}}{u}\bigg)  \Phi_2\bigg( \frac{u }{|x-y|}\bigg) \ell\bigg( \frac{\delta_D(y) \vee t^{1/\alpha}}{u} \bigg) \frac{du}{u^{\alpha+1}} \\
	&\le \frac{c_3}{|x-y|^{\alpha}}\Phi_0\bigg( \frac{\delta_D(x) \vee t^{1/\alpha}}{|x-y|}\bigg)   \ell\bigg( \frac{\delta_D(y) \vee t^{1/\alpha}}{|x-y|} \bigg)\\
	&\le  c_4\int_{|x-y|/4}^{|x-y|}\Phi_0\bigg( \frac{\delta_D(x) \vee t^{1/\alpha}}{u}\bigg) \Phi_2\bigg( \frac{u }{|x-y|}\bigg) \ell\bigg( \frac{\delta_D(y) \vee t^{1/\alpha}}{u} \bigg) \frac{du}{u^{\alpha+1}} .
\end{align*}
Using these two displays, one can deduce that the right-hand side of \eqref{e:two-jumps-equivalent} is comparable to
\begin{align*}
	&\Phi_1\bigg(\frac{\delta_D(y)\vee t^{1/\alpha}}{|x-y|}\bigg) \\
	&\times \int_{(\delta_D(x)\vee \delta_D(y)\vee t^{1/\alpha})\wedge (\eps_1|x-y|/2)}^{\eps_1|x-y|}\Phi_0\bigg(\frac{\delta_D(x)\vee t^{1/\alpha}}{u} \bigg)   \Phi_2\bigg(\frac{u}{|x-y|}\bigg) \ell \bigg(\frac{\delta_D(y)\vee t^{1/\alpha}}{u}\bigg) \frac{du}{u^{\alpha+1}}.
\end{align*}
Now the assertion follows from \eqref{e:two-jumps-equivalent-1} and \eqref{e:two-jumps-equivalent-2}.\qed

\begin{lemma}\label{l:two-jumps-lower-bound}
	(i)	There exists $C>0$ such that  for all $t>0$ and $x,y \in D$,
	\begin{align}\label{e:two-jumps-lower-bound}
		\begin{split} 
			& |x-y|^\alpha\int_{(\delta_D(x) \vee \delta_D(y)\vee t^{1/\alpha}) \wedge (\eps_1|x-y|/2)}^{\eps_1|x-y|}A_{\Phi_1,\Phi_2,\ell}(t,x,x+u\mathbf{n}_x)\,A_{\Phi_1,\Phi_2,\ell}(t,x+u\mathbf{n}_x,y)\frac{du}{u^{\alpha+1}}\\
			& \ge CA_{\Phi_0, \Phi_0, 1}(t,x,y).
		\end{split}
	\end{align}
	
	\noindent (ii) There exists $C>0$ such that  for all $t>0$ and $x,y \in D$,
	\begin{align}\label{e:two-jumps-lower-bound-t}
		\begin{split} 
			& (t \wedge |x-y|^\alpha)\int_{t^{1/\alpha} \wedge (\eps_1|x-y|/2)}^{\eps_1|x-y|}A_{\Phi_1,\Phi_2,\ell}(t,x,x+u\mathbf{n}_x)\,A_{\Phi_1,\Phi_2,\ell}(t,x+u\mathbf{n}_x,y)\frac{du}{u^{\alpha+1}}\\
			& \ge CA_{\Phi_1,\Phi_2,\ell}(t,x,y).
		\end{split}
	\end{align}
	
\end{lemma}
\pf (i) By \eqref{e:Phi0-scaling}, \eqref{e:Phi2-scaling} (with  $\Phi_2(1)=1$) and \eqref{e:ell-scaling}, we see that for all $u\in (|x-y|/4,|x-y|)$,
\begin{align*} 
	&\Phi_0\bigg(\frac{\delta_D(x)\vee t^{1/\alpha}}{u} \bigg)   \Phi_2\bigg(\frac{u}{|x-y|}\bigg) \ell \bigg(\frac{\delta_D(y)\vee t^{1/\alpha}}{u}\bigg) \frac{1}{u^{\alpha+1}} \\
	&\asymp \Phi_0\bigg(\frac{\delta_D(x)\vee t^{1/\alpha}}{|x-y|} \bigg)  \ell \bigg(\frac{\delta_D(y)\vee t^{1/\alpha}}{|x-y|}\bigg) \frac{1}{|x-y|^{\alpha+1}} .
\end{align*}
Thus, by using Lemma \ref{l:two-jumps-equivalent} and  \eqref{e:def-Phi0}, we deduce that  the left-hand side of \eqref{e:two-jumps-lower-bound} is bounded below by
\begin{align*}
	\frac{c_1}{|x-y|}\Phi_0\bigg(\frac{\delta_D(x)\vee t^{1/\alpha}}{|x-y|} \bigg)  
	\Phi_0 \bigg(\frac{\delta_D(y)\vee t^{1/\alpha}}{|x-y|}\bigg)\int_{|x-y|/4}^{|x-y|} du = c_2A_{\Phi_0,\Phi_0,1}(t,x,y).
\end{align*}

\noindent (ii) If $t^{1/\alpha} \ge \eps_1|x-y|/2$, then by (i), since $\Phi_0$ is almost increasing, the left-hand side of \eqref{e:two-jumps-lower-bound-t} is bounded below by
\begin{align*}
	&(\eps_1/2)^\alpha |x-y|^\alpha\int_{ \eps_1|x-y|/2}^{\eps_1|x-y|}A_{\Phi_1,\Phi_2,\ell}(t,x,x+u\mathbf{n}_x)\,A_{\Phi_1,\Phi_2,\ell}(t,x+u\mathbf{n}_x,y)\frac{du}{u^{\alpha+1}} \\
	&\ge c_3A_{\Phi_0,\Phi_0,1}(t,x,y) \ge  c_4\Phi_0 \bigg( \frac{t^{1/\alpha}}{|x-y|}\bigg)^2 \ge c_5 \Phi_0( \eps_1/2)^2.
\end{align*}
Thus, the result follows from \eqref{e:upper-bound-A}.

Suppose $t^{1/\alpha} < \eps_1|x-y|/2$. By \eqref{e:upper-bound-A} and  \eqref{e:interior-lower-bound-A},  we have
\begin{align}\label{e:two-jumps-lower-bound-t-1}
	A_{\Phi_1,\Phi_2,\ell}(t,x,x+u\mathbf{n}_x)\asymp 1\quad \text{for all $u\in (t^{1/\alpha},2t^{1/\alpha})$.}
\end{align}
Further,   for all $u\in (t^{1/\alpha},2t^{1/\alpha})$, by \eqref{e:lifting-property}, we have
\begin{align*}
	&	\delta_D(x+u\bn_n) \vee \delta_D(y) \vee t^{1/\alpha} \asymp  \delta_D(x)  \vee \delta_D(y)\vee t^{1/\alpha}
\end{align*}
and
\begin{align*}
	&	(\delta_D(x+u\bn_n) \wedge \delta_D(y)) \vee t^{1/\alpha} \asymp \big( (\delta_D(x) \vee t^{1/\alpha}) \wedge \delta_D(y) \big)\vee t^{1/\alpha}\\
	&\asymp \big( (\delta_D(x) \vee t^{1/\alpha}) \wedge (\delta_D(y) \vee t^{1/\alpha} )\big)  \asymp (\delta_D(x) \wedge \delta_D(y))\vee t^{1/\alpha}.
\end{align*}
Note that  $|x+u\bn_x - y| \asymp |x-y|$ for all $u\in (t^{1/\alpha}, 2t^{1/\alpha})$ since $t^{1/\alpha} < \eps_1|x-y|/2$. Thus, by \eqref{e:Phi1-scaling}, \eqref{e:Phi2-scaling}  and \eqref{e:ell-scaling}, we obtain
\begin{align} \label{e:two-jumps-lower-bound-t-2}
	A_{\Phi_1,\Phi_2,\ell}(t,x+u\mathbf{n}_x,y) & \asymp 
	A_{\Phi_1,\Phi_2,\ell}(t,x,y) \quad \text{for all $u\in (t^{1/\alpha},2t^{1/\alpha})$.}
\end{align}
Combining \eqref{e:two-jumps-lower-bound-t-1} and \eqref{e:two-jumps-lower-bound-t-2}, we conclude that the left-hand side of \eqref{e:two-jumps-lower-bound-t} is bounded below by
\begin{align*}
	c_6tA_{\Phi_1,\Phi_2,\ell}(t,x,y)\int_{t^{1/\alpha}}^{2t^{1/\alpha}} \frac{du}{u^{\alpha+1}} =c_7A_{\Phi_1,\Phi_2,\ell}(t,x,y).
\end{align*}
The proof is complete.
\qed 

\begin{lemma}\label{l:equivalences}
Suppose $\beta_1^*<\alpha+\beta_1$. For all $t>0$ and all $x,y\in D$ with $r:=|x-y|$, define 
\begin{align}
I_1&:= t\int_{ t^{1/\alpha}}^{\eps_1r}A_{\Phi_1,\Phi_2,\ell}(t,x,x+u\mathbf{n}_x)\,A_{\Phi_1,\Phi_2,\ell}(t,x+u\mathbf{n}_x,y)\frac{du}{u^{\alpha+1}}\label{e:equivalences-1}\\
		& + t\int_{ t^{1/\alpha}}^{\eps_1r}A_{\Phi_1,\Phi_2,\ell}(t,y,y+u\mathbf{n}_y)\,A_{\Phi_1,\Phi_2,\ell}(t,y+u\mathbf{n}_y,x)\frac{du}{u^{\alpha+1}},\nn
\end{align}
\begin{align}		
		I_2&:=A_{\Phi_1, \Phi_2,  \ell} (t,x,y) \label{e:equivalences-2}\\
		& +  t \int_{(\delta_D(x) \vee \delta_D(y)\vee t^{1/\alpha}) \wedge (\eps_1r/2)}^{\eps_1r}A_{\Phi_1,\Phi_2,\ell}(t,x,x+u\mathbf{n}_x)\,A_{\Phi_1,\Phi_2,\ell}(t,x+u\mathbf{n}_x,y)\frac{du}{u^{\alpha+1}}\nn\\
		&+  t \int_{(\delta_D(x) \vee \delta_D(y)\vee t^{1/\alpha}) \wedge (\eps_1r/2)}^{\eps_1r}A_{\Phi_1,\Phi_2,\ell}(t,y,y+u\mathbf{n}_y)\,A_{\Phi_1,\Phi_2,\ell}(t,y+u\mathbf{n}_y,x)\frac{du}{u^{\alpha+1}}, \nn
\end{align}
\begin{align}
		I_3&:= A_{\Phi_1,\Phi_2,\ell}(t,x,y)\label{e:equivalences-3} \\
		& +  t\Phi_1\bigg(\frac{\delta_D(x)\vee t^{1/\alpha}}{r}\bigg)  \!\!
		\int_{(\delta_D(x)\vee \delta_D(y)\vee t^{1/\alpha})\wedge r}^{r}\!\!\!\Phi_0\bigg(\frac{\delta_D(y)\vee t^{1/\alpha}}{u} \bigg)   \Phi_2\bigg(\frac{u}{r}\bigg) \ell \bigg(\frac{\delta_D(x)\vee t^{1/\alpha}}{u}\bigg) \frac{du}{u^{\alpha+1}}\nn\\
		& +  t \Phi_1\bigg(\frac{\delta_D(y)\vee t^{1/\alpha}}{r}\bigg) \!\!
		\int_{(\delta_D(x)\vee \delta_D(y) \vee  t^{1/\alpha})\wedge r}^{r}\!\!\!\Phi_0\bigg(\frac{\delta_D(x)\vee t^{1/\alpha}}{u} \bigg)   \Phi_2\bigg(\frac{u}{r}\bigg) \ell \bigg(\frac{\delta_D(y)\vee t^{1/\alpha}}{u}\bigg) \frac{du}{u^{\alpha+1}}.\nn
\end{align}
Then there exist comparison constants such that for all $t>0$  satisfying $t^{1/\alpha}<\eps_1r/2$,
$I_1\asymp I_2 \asymp I_3$.
\end{lemma}
\pf 
By symmetry, we can assume $\delta_D(x) \le \delta_D(y)$. 
From Lemmas \ref{l:two-jumps-equivalent} and \ref{l:two-jumps-lower-bound}(ii), we obtain that $I_3\le c_1I_2\le c_2 I_1$.   Thus, it suffices to show that $I_1\le c_3I_3$. Set $\eps:=(\alpha+\beta_1-\beta_1^*)/4>0$.

If $\delta_D(y)>\eps_1r/2$, then  for all $u \in (t^{1/\alpha}, \eps_1 r)$, using  \eqref{e:lifting-property},  \eqref{e:Phi0-scaling}-- \eqref{e:ell-scaling}, the fact that $\Phi_2(u)=1$ for $u\ge 1$, and the boundedness of $\Phi_0$, we obtain
\begin{align*}
	&	A_{\Phi_1,\Phi_2,\ell}(t,x,x+u\mathbf{n}_x)A_{\Phi_1,\Phi_2,\ell}(t,x+u\mathbf{n}_x,y)\\
	& \asymp\Phi_0 \bigg( \frac{\delta_D(x) \vee t^{1/\alpha}}{u} \bigg)  \Phi_2 \bigg( \frac{\delta_D(x) \vee u}{u} \bigg)  \Phi_1 \bigg( \frac{\delta_D(x) \vee u}{r} \bigg)	 \Phi_2\bigg(\frac{\delta_D(y)}{r} \bigg) \ell \bigg( \frac{\delta_D(x) \vee u}{\delta_D(y) \wedge r}\bigg)\\
	& \asymp\Phi_0 \bigg( \frac{\delta_D(x) \vee t^{1/\alpha}}{u} \bigg) \Phi_0 \bigg( \frac{\delta_D(x) \vee u}{r} \bigg)  \\
	&\le c_4 \Phi_0 \bigg( \frac{\delta_D(x) \vee t^{1/\alpha}}{r} \bigg)   \times  \begin{cases}
		\displaystyle  1&\mbox{ if $\delta_D(x)\ge u$},\\
		\displaystyle \bigg( \frac{\delta_D(x)\vee t^{1/\alpha}}{u}\bigg)^{\beta_1-\eps}    \bigg( \frac{u}{\delta_D(x)\vee t^{1/\alpha}}\bigg)^{\beta_1^*+\eps} &\mbox{ if $\delta_D(x)<u$}
	\end{cases}\\
	&\le c_5  \bigg( \frac{u}{ t^{1/\alpha}}\bigg)^{\beta_1^*-\beta_1+2\eps}   \Phi_0 \bigg( \frac{\delta_D(x) \vee t^{1/\alpha}}{r} \bigg)    \asymp    \bigg( \frac{u}{ t^{1/\alpha}}\bigg)^{\beta_1^*-\beta_1+2\eps}    A_{\Phi_1,\Phi_2,\ell}(t,x,y)
\end{align*}
and
\begin{align*}
	&	A_{\Phi_1,\Phi_2,\ell}(t,y,y+u\mathbf{n}_y)\,A_{\Phi_1,\Phi_2,\ell}(t,y+u\mathbf{n}_y,x)   \\
	&\asymp A_{\Phi_1,\Phi_2,\ell}(t,y+u\mathbf{n}_y,x)  \asymp  \Phi_0 \bigg( \frac{\delta_D(x) \vee t^{1/\alpha}}{r} \bigg)   \asymp A_{\Phi_1,\Phi_2,\ell}(t,x,y).
\end{align*}
Since $\alpha+\beta_1-\beta_1^*-2\eps>0$, it follows that
\begin{align*}
I_1\le  c_6t A_{\Phi_1,\Phi_2,\ell}(t,x,y)\int_{t^{1/\alpha}}^{\eps_1r} 
\bigg(  \bigg( \frac{u}{ t^{1/\alpha}}\bigg)^{\beta_1^*-\beta_1+2\eps} +1 \bigg) \frac{du}{u^{\alpha+1}}  
\le  c_7  A_{\Phi_1,\Phi_2,\ell}(t,x,y) \le c_7 I_3.
\end{align*}

Suppose $\delta_D(y) \le \eps_1r/2$. Then by Lemma  \ref{l:two-jumps-equivalent}, we get $I_2 \asymp I_3$. 
If $\delta_D(y)\le t^{1/\alpha}$, then the lower bound in the integrals in $I_2$ is equal to $t^{1/\alpha}$, hence these integrals are equal to the ones in $I_1$. Thus clearly $I_1\le I_2\asymp I_3$.

If $\delta_D(y)>t^{1/\alpha}$, we split $I_1$ into two parts:
\begin{align*}
	I_1'	&:= t\int_{ t^{1/\alpha}}^{ \delta_D(y)}A_{\Phi_1,\Phi_2,\ell}(t,x,x+u\mathbf{n}_x)\,A_{\Phi_1,\Phi_2,\ell}(t,x+u\mathbf{n}_x,y)\frac{du}{u^{\alpha+1}}\\
	& \quad\; + t\int_{ t^{1/\alpha}}^{\delta_D(y)}A_{\Phi_1,\Phi_2,\ell}(t,y,y+u\mathbf{n}_y)\,A_{\Phi_1,\Phi_2,\ell}(t,y+u\mathbf{n}_y,x)\frac{du}{u^{\alpha+1}}
\end{align*}
and
\begin{align*}
	I_1''	&:= t\int_{ \delta_D(y)}^{\eps_1 r}A_{\Phi_1,\Phi_2,\ell}(t,x,x+u\mathbf{n}_x)\,A_{\Phi_1,\Phi_2,\ell}(t,x+u\mathbf{n}_x,y)\frac{du}{u^{\alpha+1}}\\
	& \quad\; + t\int_{ \delta_D(y)}^{\eps_1 r}A_{\Phi_1,\Phi_2,\ell}(t,y,y+u\mathbf{n}_y)\,A_{\Phi_1,\Phi_2,\ell}(t,y+u\mathbf{n}_y,x)\frac{du}{u^{\alpha+1}}
\end{align*}
Since $I_1''$ is equal to  the sum of the integrals in $I_2$, 
in order to obtain $I_1 \le c_3I_3$ it suffices to prove that $I_1'\le c_8A_{\Phi_1,\Phi_2,\ell}(t,x,y)$.
Since $\delta_D(y)>t^{1/\alpha}$, for all $u \in (t^{1/\alpha}, \delta_D(y))$, using   \eqref{e:Phi0-scaling}-- \eqref{e:ell-scaling}  and the boundedness of $\Phi_0$, we obtain
\begin{align*}
	&	A_{\Phi_1,\Phi_2,\ell}(t,x,x+u\mathbf{n}_x)A_{\Phi_1,\Phi_2,\ell}(t,x+u\mathbf{n}_x,y) \\
	& \asymp \Phi_0 \bigg( \frac{\delta_D(x) \vee t^{1/\alpha}}{u} \bigg)\Phi_1 \bigg( \frac{\delta_D(x) \vee u}{r} \bigg) \Phi_2 \bigg( \frac{\delta_D(y) }{r}  \bigg) \ell \bigg( \frac{\delta_D(x) \vee u}{\delta_D(y) }\bigg) \\
	&\le c_9 \Phi_1 \bigg( \frac{\delta_D(x) \vee t^{1/\alpha}}{r} \bigg) \Phi_2 \bigg( \frac{\delta_D(y) }{r}  \bigg) \ell \bigg( \frac{\delta_D(x) \vee t^{1/\alpha}}{\delta_D(y) }\bigg)\\
	&\quad  \times  \begin{cases}
		\displaystyle  1&\mbox{ if $\delta_D(x)\ge u$},\\
		\displaystyle \bigg( \frac{\delta_D(x) \vee t^{1/\alpha}}{ u}\bigg)^{\beta_1-\eps} \bigg( \frac{ u}{\delta_D(x) \vee t^{1/\alpha}}\bigg)^{\beta_1^*+\eps} \bigg( \frac{ u}{\delta_D(x) \vee t^{1/\alpha}}\bigg)^\eps  &\mbox{ if $\delta_D(x)<u$}
	\end{cases}\\
	&\le c_{10} \bigg( \frac{u}{ t^{1/\alpha}}\bigg)^{\beta_1^*-\beta_1+3\eps}  \Phi_1 \bigg( \frac{\delta_D(x) \vee t^{1/\alpha}}{r} \bigg) \Phi_2 \bigg( \frac{\delta_D(y) }{r}  \bigg) \ell \bigg( \frac{\delta_D(x) \vee t^{1/\alpha}}{\delta_D(y) }\bigg)\\
	& \asymp  \bigg( \frac{u}{ t^{1/\alpha}}\bigg)^{\beta_1^*-\beta_1+3\eps} A_{\Phi_1, \Phi_2,\ell} (t,x,y) 
\end{align*}
and
\begin{align*}
	&	A_{\Phi_1,\Phi_2,\ell}(t,y,y+u\mathbf{n}_y)\,A_{\Phi_1,\Phi_2,\ell}(t,y+u\mathbf{n}_y,x) \\
	&\asymp  A_{\Phi_1,\Phi_2,\ell}(t,y+u\mathbf{n}_y,x) \asymp  \Phi_1 \bigg( \frac{\delta_D(x) \vee t^{1/\alpha}}{r} \bigg) \Phi_2 \bigg( \frac{\delta_D(y) }{r}  \bigg) \ell \bigg( \frac{\delta_D(x) \vee t^{1/\alpha}}{\delta_D(y) }\bigg)  \asymp A_{\Phi_1,\Phi_2,\ell}(t,x,y). 
\end{align*}
Since $\alpha+\beta_1-\beta_1^*-3\eps>0$,  we arrive at
\begin{align*}
	I_1'\le  c_{11}t A_{\Phi_1,\Phi_2,\ell}(t,x,y)\int_{t^{1/\alpha}}^{\delta_D(y)} \bigg( \bigg( \frac{u}{ t^{1/\alpha}}\bigg)^{\beta_1^*-\beta_1+3\eps}+1\bigg) \frac{du}{u^{\alpha+1}} \le  c_{12}  A_{\Phi_1,\Phi_2,\ell}(t,x,y) .
\end{align*}
The proof is complete. \qed 

\begin{lemma}\label{l:two-jumps-comparability}
	(i)
	If $\beta_1^*\vee \beta_2^*<\alpha+\beta_1$, then
	there exists $C>0$ such that for all $t>0$ and $x,y\in D$,
	\begin{align}\label{e:two-jumps-comparability-1}
		\begin{split} 
			&\left(t \wedge |x-y|^{\alpha}\right)\!\! \int_{(\delta_D(x) \vee \delta_D(y)\vee t^{1/\alpha}) \wedge (\eps_1|x-y|/2)}^{\eps_1|x-y|}\!\!\!A_{\Phi_1,\Phi_2,\ell}(t,x,x+u\mathbf{n}_x)\,A_{\Phi_1,\Phi_2,\ell}(t,x+u\mathbf{n}_x,y)\frac{du}{u^{\alpha+1}}\\
			&\le C A_{\Phi_1, \Phi_2, \ell}(t,x,y).
		\end{split} 
	\end{align}

	\noindent	(ii) 
	If $\beta_2>\alpha+\beta_1^*$, then
	there exist comparison constants such that for $(t,x,y) \in (0,\infty)\times D\times D$,
	\begin{align}\label{e:two-jumps-comparability-2}
		\begin{split} 
			&\left(t \wedge |x-y|^{\alpha}\right)\!\! \int_{(\delta_D(x) \vee \delta_D(y)\vee t^{1/\alpha}) 
				\wedge (\eps_1|x-y|/2)}^{\eps_1|x-y|}\!\!\!A_{\Phi_1,\Phi_2,\ell}(t,x,x+u\mathbf{n}_x)\,
			A_{\Phi_1,\Phi_2,\ell}(t,x+u\mathbf{n}_x,y)\frac{du}{u^{\alpha+1}}\\
			&\asymp \bigg(1 \wedge \frac{t}{|x-y|^\alpha}\bigg) A_{\Phi_0, \Phi_0, 1}(t,x,y).
		\end{split} 
	\end{align}
\end{lemma}
\pf  	By Lemma \ref{l:two-jumps-equivalent}, the (common) left-hand side
of \eqref{e:two-jumps-comparability-1} and \eqref{e:two-jumps-comparability-2} is comparable to
\begin{align*}I&:= \left(t \wedge |x-y|^{\alpha}\right) \int_{(\delta_D(x)\vee \delta_D(y)\vee t^{1/\alpha})\wedge 
		(|x-y|/4)}^{|x-y|}   f(u) \frac{du}{u^{\alpha+1}},
\end{align*} 
where
\begin{align*}
	f(u):=\Phi_1\bigg(\frac{\delta_D(y)\vee t^{1/\alpha}}{|x-y|}\bigg) \Phi_0\bigg(\frac{\delta_D(x)\vee t^{1/\alpha}}{u} \bigg)   \Phi_2\bigg(\frac{u}{|x-y|}\bigg) \ell \bigg(\frac{\delta_D(y)\vee t^{1/\alpha}}{u}\bigg).
\end{align*}

We first assume that $t^{1/\alpha}\ge |x-y|/4$. By \eqref{e:upper-bound-A}, $I$ is bounded above by 
\begin{align*}
	c_1|x-y|^\alpha \int_{\eps_1|x-y|/2}^{\eps_1|x-y|} \frac{du}{u^{\alpha+1}} =c_2.
\end{align*}
Hence, \eqref{e:two-jumps-comparability-1} follows  from \eqref{e:interior-lower-bound-A}.  Moreover, when $t^{1/\alpha}\ge |x-y|/4$, we see from the boundedness of $\Phi_0$ and  \eqref{e:Phi0-scaling} that
\begin{align*}
	\bigg(1 \wedge \frac{t}{|x-y|^\alpha}\bigg) A_{\Phi_0, \Phi_0, 1}(t,x,y) \asymp 1,
\end{align*}
and from  \eqref{e:Phi0-scaling}, \eqref{e:Phi2-scaling}, \eqref{e:ell-scaling}, and  the definition  \eqref{e:def-Phi0} and  the almost increasing property of $\Phi_0$ that 
\begin{align*}
	I	&\ge 4^{-\alpha}|x-y|^\alpha \int_{ |x-y|/4}^{|x-y|}f(u)\frac{du}{u^{\alpha+1}} \\
	&\ge c_3|x-y|^\alpha \Phi_0\bigg(\frac{\delta_D(x)\vee t^{1/\alpha}}{|x-y|}\bigg) \Phi_0\bigg(\frac{\delta_D(y)\vee t^{1/\alpha}}{|x-y|} \bigg)   \Phi_2(1)\int_{ |x-y|/4}^{|x-y|} \frac{du}{u^{\alpha+1}}  \ge c_4 \Phi_0(1/4)^2 .
\end{align*} 
Thus the assertions are valid if  if  $t^{1/\alpha}\ge |x-y|/4$, without extra assumptions on $\beta_1$ and $\beta_2$.

We now assume that $t^{1/\alpha}< |x-y|/4$.

\smallskip

\noindent (i) 
Let $\eps:=(\alpha+\beta_1- \beta_1^*\vee \beta_2^*)/4$.  If $\delta_D(x) \vee \delta_D(y) \ge |x-y|/4$, then using \eqref{e:Phi0-scaling}, \eqref{e:Phi2-scaling}, \eqref{e:ell-scaling}, we see that
\begin{align}\label{e:two-jumps-comparability-case0}
	\begin{split} 
		I&=t \int_{|x-y|/4}^{|x-y|}   f(u) \frac{du}{u^{\alpha+1}} 
		\asymp
		t \Phi_0\bigg(\frac{\delta_D(x)\vee t^{1/\alpha}}{|x-y|}\bigg) \Phi_0\bigg(\frac{\delta_D(y)\vee t^{1/\alpha}}{|x-y|} \bigg)  \Phi_2(1) \int_{|x-y|/4}^{|x-y|}    \frac{du}{u^{\alpha+1}}\\
		&
		\le \frac{c_5t}{|x-y|^\alpha}
		\Phi_0\bigg(\frac{\delta_D(x)\vee t^{1/\alpha}}{|x-y|} \bigg)\Phi_0\bigg(\frac{\delta_D(y)\vee t^{1/\alpha}}{|x-y|} \bigg).
	\end{split}
\end{align}
We distinguish between two cases.

\smallskip

\noindent
\textbf{Case 1:} $\delta_D(x) \ge \delta_D(y)$.
If $\delta_D(x)\ge |x-y|/4$, we have 
$(\delta_D(x) \vee t^{1/\alpha}) \wedge |x-y| \asymp |x-y|.$
Thus using  \eqref{e:Phi1-scaling}, \eqref{e:Phi2-scaling} and \eqref{e:def-Phi0}, we get
\begin{align*}
	&A_{\Phi_1,\Phi_2,\ell}(t,x,y) \ge c_6  \Phi_1\bigg( \frac{\delta_D(y) \vee t^{1/\alpha}}{|x-y|}\bigg)  \ell \bigg( \frac{\delta_D(y)\vee t^{1/\alpha}}{(\delta_D(x) \vee t^{1/\alpha}) \wedge |x-y|} \bigg) \ge c_7 \Phi_0\bigg( \frac{\delta_D(y) \vee t^{1/\alpha}}{|x-y|}\bigg) . 
\end{align*}
Thus, by \eqref{e:two-jumps-comparability-case0} with  
$t^{1/\alpha}\le |x-y|/4$ 
and the boundedness of $\Phi_0$, the assertion holds in this case.

Suppose $\delta_D(x)<|x-y|/4$. For all $\delta_D(x)\vee t^{1/\alpha}<u<|x-y|$, applying \eqref{e:Phi0-scaling}, \eqref{e:Phi2-scaling} and \eqref{e:ell-scaling},  we obtain
\begin{align*}
	&\Phi_1\bigg(\frac{\delta_D(y)\vee t^{1/\alpha}}{|x-y|}\bigg) \Phi_0\bigg(\frac{\delta_D(x)\vee t^{1/\alpha}}{u} \bigg)   \Phi_2\bigg(\frac{u}{|x-y|}\bigg) \ell \bigg(\frac{\delta_D(y)\vee t^{1/\alpha}}{u}\bigg)\\ &\le c_8\Phi_0(1)  \bigg( \frac{u}{\delta_D(x)\vee t^{1/\alpha}}\bigg)^{-\beta_1+\eps + \beta_2^*+\eps + \eps}	 \Phi_1\bigg(\frac{\delta_D(y)\vee t^{1/\alpha}}{|x-y|}\bigg)  \Phi_2\bigg(\frac{\delta_D(x)\vee t^{1/\alpha}}{|x-y|}\bigg) \ell \bigg(\frac{\delta_D(y)\vee t^{1/\alpha}}{\delta_D(x) \vee t^{1/\alpha}}\bigg)\\
	&=c_8  \bigg( \frac{u}{\delta_D(x)\vee t^{1/\alpha}}\bigg)^{-\beta_1 + \beta_2^*+3\eps} 	A_{\Phi_1,\Phi_2,\ell}(t,x,y).
\end{align*}
Since $\alpha+\beta_1>\beta_2^*+3\eps$, it follows that
\begin{align*}
	I&\le c_8t(\delta_D(x)\vee t^{1/\alpha})^{\beta_1-\beta_2^*-3\eps}A_{\Phi_1,\Phi_2,\ell}(t,x,y) \int_{\delta_D(x)\vee t^{1/\alpha}}^{\infty}  \frac{du}{u^{1+\alpha+\beta_1-\beta_2^*-3\eps }}\\
	&= \frac{c_9t}{(\delta_D(x)\vee t^{1/\alpha})^\alpha} A_{\Phi_1,\Phi_2,\ell}(t,x,y) \le
	c_7 A_{\Phi_1,\Phi_2,\ell}(t,x,y).
\end{align*}

\noindent
\textbf{Case 2:} $\delta_D(x) < \delta_D(y)$. If $\delta_D(y)  \ge |x-y|/4$, then using \eqref{e:Phi0-scaling}, \eqref{e:Phi2-scaling}, \eqref{e:ell-scaling}, and the definition  \eqref{e:def-Phi0} of $\Phi_0$, we see that
\begin{align*}
	A_{\Phi_1,\Phi_2,\ell}(t,x,y) &\ge c_{10}\Phi_1\bigg(\frac{\delta_D(x)\vee t^{1/\alpha}}{|x-y|} \bigg) \ell \bigg(\frac{\delta_D(x)\vee t^{1/\alpha}}{|x-y|}\bigg) = c_{10}\Phi_0\bigg(\frac{\delta_D(x)\vee t^{1/\alpha}}{|x-y|} \bigg).
\end{align*}
Hence, the assertion follows from \eqref{e:two-jumps-comparability-case0} with  $t^{1/\alpha}\le |x-y|$ and the boundedness of $\Phi_0$.

Suppose $\delta_D(y)<|x-y|/4$. For all $\delta_D(y)\vee t^{1/\alpha}<u<|x-y|$, using  \eqref{e:Phi0-scaling}, \eqref{e:Phi2-scaling} and \eqref{e:ell-scaling}, we get that
\begin{align*}
	&\Phi_1\bigg(\frac{\delta_D(y)\vee t^{1/\alpha}}{|x-y|}\bigg) \Phi_0\bigg(\frac{\delta_D(x)\vee t^{1/\alpha}}{u} \bigg)   \Phi_2\bigg(\frac{u}{|x-y|}\bigg) \ell \bigg(\frac{\delta_D(y)\vee t^{1/\alpha}}{u}\bigg)\\ 
	&\le c_{11} \bigg( \frac{\delta_D(x)\vee t^{1/\alpha}}{u}\bigg)^{\beta_1-\eps} \bigg( \frac{\delta_D(y)\vee t^{1/\alpha}}{\delta_D(x)\vee t^{1/\alpha}}\bigg)^{\beta_1^*+\eps } \bigg( \frac{u}{ \delta_D(y)\vee t^{1/\alpha}}\bigg)^{\beta_2^*+\eps } \bigg( \frac{u}{\delta_D(y)\vee t^{1/\alpha}}\bigg)^\eps \\
	&\quad \times 	\Phi_0(1)\ell (1) \Phi_1\bigg(\frac{\delta_D(x)\vee t^{1/\alpha}}{|x-y|}\bigg)  \Phi_2\bigg(\frac{ \delta_D(y) \vee t^{1/\alpha}}{|x-y|}\bigg)  \\
	&\le c_{12}\bigg( \frac{\delta_D(x)\vee t^{1/\alpha}}{u}\bigg)^{\beta_1-\eps} \bigg( \frac{\delta_D(y)\vee t^{1/\alpha}}{\delta_D(x)\vee t^{1/\alpha}}\bigg)^{\beta_1^*+\eps } \bigg( \frac{u}{ \delta_D(y)\vee t^{1/\alpha}}\bigg)^{\beta_2^*+\eps } \bigg( \frac{u}{\delta_D(y)\vee t^{1/\alpha}}\bigg)^\eps \\
	&\quad \times \bigg( \frac{\delta_D(y)\vee t^{1/\alpha}}{\delta_D(x)\vee t^{1/\alpha}}\bigg)^\eps A_{\Phi_1,\Phi_2,\ell}(t,x,y)  \\
	&= c_{12} (\delta_D(x)\vee t^{1/\alpha})^{\beta_1-\beta_1^*-3\eps} (\delta_D(y) \vee t^{1/\alpha})^{\beta_1^*-\beta_2^*} u^{-\beta_1+\beta_2^* +3\eps}  A_{\Phi_1,\Phi_2,\ell}(t,x,y).
\end{align*}
Since $\alpha+\beta_1>\beta_1^*\vee \beta_2^*+3\eps$, it follows that
\begin{align*}
	I&\le c_{12}t(\delta_D(x)\vee t^{1/\alpha})^{\beta_1-\beta_1^*-3\eps} (\delta_D(y) \vee t^{1/\alpha})^{\beta_1^*-\beta_2^*}A_{\Phi_1,\Phi_2,\ell}(t,x,y) \int_{\delta_D(y)\vee t^{1/\alpha}}^{\infty}  \frac{du}{u^{1+\alpha+\beta_1-\beta_2^*-3\eps }}\\
	&= \frac{c_{13}t(\delta_D(x)\vee t^{1/\alpha})^{\beta_1-\beta_1^*-3\eps}}{(\delta_D(y)\vee t^{1/\alpha})^{\alpha+\beta_1-\beta_1^*-3\eps}} A_{\Phi_1,\Phi_2,\ell}(t,x,y) \le\frac{c_{13}(\delta_D(x)\vee t^{1/\alpha})^{\alpha+\beta_1-\beta_1^*-3\eps}}{(\delta_D(y)\vee t^{1/\alpha})^{\alpha+\beta_1-\beta_1^*-3\eps}} A_{\Phi_1,\Phi_2,\ell}(t,x,y)\\
	&\le c_{13}A_{\Phi_1,\Phi_2,\ell}(t,x,y).
\end{align*}

\noindent (ii) Let $\eps':=(\beta_2-\alpha-\beta_1^*)/4$. 
If $(\delta_D(x)\vee \delta_D(y)\vee  t^{1/\alpha}) \wedge (|x-y|/4)<u<|x-y|$, by using \eqref{e:Phi0-scaling}, \eqref{e:Phi2-scaling} and \eqref{e:ell-scaling},   we get that 
\begin{align*}
	&\Phi_1\bigg(\frac{\delta_D(y)\vee t^{1/\alpha}}{|x-y|}\bigg) \Phi_0\bigg(\frac{\delta_D(x)\vee t^{1/\alpha}}{u} \bigg)   \Phi_2\bigg(\frac{u}{|x-y|}\bigg) \ell \bigg(\frac{\delta_D(y)\vee t^{1/\alpha}}{u}\bigg)\\ &\le c_{14}\Phi_2(1)  \bigg( \frac{|x-y|}{u}\bigg)^{\beta_1^*+\eps' - \beta_2+\eps' + \eps'}	 \Phi_1\bigg(\frac{\delta_D(y)\vee t^{1/\alpha}}{|x-y|}\bigg)  \Phi_0\bigg(\frac{\delta_D(x)\vee t^{1/\alpha}}{|x-y|}\bigg) \ell \bigg(\frac{\delta_D(y)\vee t^{1/\alpha}}{|x-y|}\bigg)\\
	&=c_{14} \bigg( \frac{|x-y|}{u}\bigg)^{\beta_1^* - \beta_2+3\eps'}	 	A_{\Phi_0,\Phi_0,1}(t,x,y)
\end{align*}
and
\begin{align*}
	&\Phi_1\bigg(\frac{\delta_D(y)\vee t^{1/\alpha}}{|x-y|}\bigg) \Phi_0\bigg(\frac{\delta_D(x)\vee t^{1/\alpha}}{u} \bigg)   \Phi_2\bigg(\frac{u}{|x-y|}\bigg) \ell \bigg(\frac{\delta_D(y)\vee t^{1/\alpha}}{u}\bigg)\\ &\ge c_{15}\bigg( \frac{|x-y|}{u}\bigg)^{\beta_1 - \beta_2^*-3\eps' }	 	A_{\Phi_0,\Phi_0,1}(t,x,y).
\end{align*}
Since $\beta_2>\alpha+\beta_1^*+3\eps$, it follows that
\begin{align*}
	I &\le c_{14} t |x-y|^{\beta_1^*-\beta_2+3\eps'} 	A_{\Phi_0,\Phi_0,1}(t,x,y) \int_0^{|x-y|} \frac{du}{u^{1+\alpha + \beta_1^*-\beta_2+3\eps'}} = \frac{c_{16}t	A_{\Phi_0,\Phi_0,1}(t,x,y) }{|x-y|^\alpha }
\end{align*}
and
\begin{align*}
	I &\ge c_{15} t |x-y|^{\beta_1-\beta_2^*-3\eps'} 	A_{\Phi_0,\Phi_0,1}(t,x,y) \int_{|x-y|/4}^{|x-y|} \frac{du}{u^{1+\alpha + \beta_1-\beta_2^*-3\eps'}} = \frac{c_{17}t	A_{\Phi_0,\Phi_0,1}(t,x,y) }{|x-y|^\alpha }
\end{align*}

The proof is complete.  
\qed

The next lemma implies that none of the  two terms in the brackets on the right hand side of \eqref{e:main11-case2} dominates the other. Moreover, if $\liminf_{r\to 0}\ell(r)>0$, then the same holds for \eqref{e:main11-case3}.

\begin{lemma}\label{l:two-jumps-non-dominant}
	(i)	If $\beta_2>\alpha+\beta_1^*$, then 
	\begin{align}\label{l:two-jumps-non-dominant-1-1}\lim_{t\to 0}\frac{ t|x-y|^{-\alpha}A_{\Phi_0, \Phi_0, 1}(t,x,y)} 
	{A_{\Phi_1, \Phi_2, \ell}(t,x,y)}   =0 \quad \text{for all $x,y \in D$, $x\neq y$},
	\end{align}
	\begin{align}\label{l:two-jumps-non-dominant-1-2}
		\lim_{t\to 0} \lim_{x\to Q, \, y \to Q'}\frac{ t|x-y|^{-\alpha}A_{\Phi_0, \Phi_0, 1}(t,x,y)} {A_{\Phi_1, \Phi_2, \ell}(t,x,y)} =  \infty \quad \text{for all $Q, Q' \in \partial D$, $Q\neq Q'$}.
	\end{align}

	\noindent (ii)  	Assume the setting of Theorem \ref{t:main11}(iii) and suppose in addition that $\liminf_{s\to 0} \ell(s) >0$. 
	Let
	\begin{align*}
		\sI(t,x,y)&:=  A_{\Phi_0,\Phi_0,1}(t,x,y)\int_{(\delta_D(x)\vee \delta_D(y)\vee t^{1/\alpha})\wedge (|x-y|/4)}^{|x-y|}\frac{du}{u}\\
		&\qquad \times \bigg( \frac{\ell((\delta_D(x) \vee t^{1/\alpha})/u) \,\ell((\delta_D(y) \vee t^{1/\alpha})/u)\, \phi(u/|x-y|)}{\ell((\delta_D(x) \vee t^{1/\alpha})/|x-y|)\, \ell((\delta_D(y) \vee t^{1/\alpha})/|x-y|) }  \bigg) .
	\end{align*} 
	Then
	\begin{align}\label{l:two-jumps-non-dominant-2-1}\lim_{t\to 0}\frac{ t|x-y|^{-\alpha}\sI(t,x,y)} {A_{\Phi_1, \Phi_2, \ell}(t,x,y)}   =0 \quad \text{for all $x,y \in D$, $x\neq y$},
	\end{align}
	\begin{align}\label{l:two-jumps-non-dominant-2-2}
		\limsup_{t\to 0} \lim_{x\to Q, \, y \to Q'}\frac{ t|x-y|^{-\alpha}\sI(t,x,y)} {A_{\Phi_1, \Phi_2, \ell}(t,x,y)} =  \infty \quad \text{for all $Q, Q' \in \partial D$, $Q\neq Q'$}.
	\end{align}
\end{lemma}
\pf 
Equalities \eqref{l:two-jumps-non-dominant-1-1} and \eqref{l:two-jumps-non-dominant-2-1} are evident since $\lim_{t\to 0} t|x-y|^{-\alpha}A_{\Phi_0,\Phi_0,1}(t,x,y)=\lim_{t\to 0} t|x-y|^{-\alpha}\sI(t,x,y)=0$ and $\lim_{t\to 0} A_{\Phi_1,\Phi_2,\ell}(t,x,y)=A_{\Phi_1,\Phi_2,\ell}(0,x,y)>0$  for all $x,y\in D$, $x\neq y$.

Let $Q,Q' \in \partial D$ be distinct and write $r:= |Q-Q'|>0$.

\noindent (i) Let $\eps:=(\beta_2-\alpha-\beta_1^*)/5>0$. For any  $t\in (0, r^\alpha)$, using \eqref{e:Phi1-scaling}--\eqref{e:ell-scaling}, we obtain
\begin{align*}
	&	\lim_{x\to Q, \, y \to Q'}\frac{ t|x-y|^{-\alpha}A_{\Phi_0, \Phi_0, 1}(t,x,y)} {A_{\Phi_1, \Phi_2, \ell}(t,x,y)} = \frac{t\Phi_0(t^{1/\alpha}/r)^2}{r^\alpha \Phi_1(t^{1/\alpha}/r)  \Phi_2(t^{1/\alpha}/r)} = \frac{t\Phi_1(t^{1/\alpha}/r) \ell(t^{1/\alpha}/r)^2 }{r^\alpha \Phi_2(t^{1/\alpha}/r)}\\
	&\ge  \frac{c_1\Phi_1(1) \ell(1)^2 }{\Phi_2(1)} \bigg( \frac{t^{1/\alpha}}{r}\bigg)^{\alpha+\beta_1^* + \eps + 2\eps - (\beta_2-\eps)} = \frac{c_1\Phi_1(1) \ell(1)^2 }{\Phi_2(1)} \bigg( \frac{t^{1/\alpha}}{r}\bigg)^{-\eps}.
\end{align*}
This yields \eqref{l:two-jumps-non-dominant-1-2}.

\noindent (ii) Since $\liminf_{s\to 0} \ell(s)>0$, we have  $c_2:=\inf_{s\in (0,1)} \ell(s)>0$. For all $t\in (0,r^\alpha)$, we have
\begin{align*}
	\lim_{x\to Q, \, y \to Q'}\frac{ t|x-y|^{-\alpha}\sI(t,x,y)} {A_{\Phi_1, \Phi_2, \ell}(t,x,y)}&=\frac{ (t^{1/\alpha}/r)^{\alpha+ 2\beta_1} \ell(t^{1/\alpha}/r)^2}{ (t^{1/\alpha}/r)^{\alpha+2\beta_1} \phi(t^{1/\alpha}/r) } \int_{t^{1/\alpha}}^r \frac{\ell(t^{1/\alpha}/u)^2 \phi(u/r)}{\ell(t^{1/\alpha}/r)^2} \frac{du}{u}\\
	&\ge \frac{c_2^2}{\phi(t^{1/\alpha}/r)} \int_{t^{1/\alpha}}^r   \frac{\phi(u/r)}{u} du = \frac{c_2^2}{\phi(t^{1/\alpha}/r)} \int_{1/r}^{1/t^{1/\alpha}}   \frac{\phi(1/(ra))}{a} da,
\end{align*}
where we used the change of the variables $u=1/a$ in the last equality. Define $f(a):=\phi(1/(ra))$. Since $\phi$ has lower and upper Matuszewska indices both equal to $0$, by \cite[Theorem 2.6.1(f)]{BGT87}, 
\begin{align*}
	\liminf_{t\to 0} \frac{\phi(t^{1/\alpha}/r)}{\int_{1/r}^{1/t^{1/\alpha}}  a^{-1} \phi(1/(ra))  da} = \liminf_{t\to 0} \frac{f(1/t^{1/\alpha})}{\int_{1/r}^{1/t^{1/\alpha}}  a^{-1} f(a)  da}=0.
\end{align*}
This leads to \eqref{l:two-jumps-non-dominant-2-2}.
\qed


\section{Preliminary upper estimates}\label{s:pue}

Recall that we work under the framework introduced at the end of Subsection \ref{ss:setup}.

\medskip

 \emph{The following notational convention will be used throughout Sections \ref{s:pue}--\ref{s:ub}.	
 	When we consider $\overline{Y}$, 
 	we  write $Y$, $Y^U$,  $p(t, x, y)$, $p^U(t, x, y)$, $\tau_U$ and 	$\zeta$ instead of  
 	$\overline{Y}$, $\overline Y^U$, $\bar{p}(t, x, y)$, $\bar{p}^{U}(t, x, y)$, $\bar{\tau}_U$ and $\infty$.	
 	When we consider $Y^\kappa$, 
 	we  write $Y$, $Y^U$,  $p(t, x, y)$, $p^U(t, x, y)$, $\tau_U$ and	$\zeta$ instead of  
 	$Y^\kappa$, $Y^{\kappa,U}$, $p^\kappa(t, x,y)$, $p^{\kappa, U}(t, x, y)$, $\tau^\kappa_U$ and $\zeta^\kappa$.	
 	Further, when we consider $\overline{Y}$, we let $q=0$, and when  we consider $Y^\kappa$, 
 	we let $q\in [(\alpha-1)_+,\alpha+\beta_1)$ be the 	(strictly)  positive constant satisfying \eqref{e:C(alpha,p,F)}.} 

\medskip 

The following proposition is the main result of this section. 
\begin{prop}\label{p:UHK-rough} 
		For any $T>0$, there exists  $C=C(T)>0$ such that
	\begin{equation}\label{e:UHK-rough}		
		p(t,x,y) \le C \left(1 \wedge \frac{\delta_D(x)}{t^{1/\alpha}}\right)^{q}\left(1 \wedge \frac{\delta_D(y)}{t^{1/\alpha}} \right)^{q} \left( t^{-d/\alpha} \wedge \frac{t}{|x-y|^{d+\alpha}}\right)
	\quad 	\text{for all 
		$t\le T$ and $x,y \in D$.}
	\end{equation}
\end{prop}

For $q=0$, \eqref{e:UHK-rough}		 is established in  Proposition \ref{p:upper-heatkernel}. Hence, we assume 
$q \in [(\alpha-1)_+,\alpha+\beta_1) \cap (0, \infty)$
 in the remainder of this section and so $D$ is a bounded $C^{1,1}$ open set and $p(t,x,y)=p^\kappa(t,x,y)$.

Note that, since $t^{-d/\alpha} \wedge \frac{t}{|x-y|^{d+\alpha}}$ is comparable to the transition density of 
the isotropic $\alpha$-stable process in $\R^d$, 
there exists  $C>0$ such that for all $t,s>0$ and  $x,y \in D$,
\begin{align}\label{e:stableu1} 
	\int_{D} \left(t^{-d/\alpha} \wedge \frac{t}{|x-z|^{d+\alpha}}\right) dz \le C 
\end{align}
and
\begin{align} \label{e:stableu2} 
	\int_{D} \left( t^{-d/\alpha} \wedge \frac{t}{|x-z|^{d+\alpha}}\right) \left( s^{-d/\alpha} \wedge \frac{s}{|y-z|^{d+\alpha}}\right)dz \le C \left( (t+s)^{-d/\alpha} \wedge 
	\frac{t+s}{|x-y|^{d+\alpha}}\right),\end{align}
where in \eqref{e:stableu2} we used the semigroup property.

Recall that $\eta_0\in (0,1/36]$ denotes the constant in Theorem \ref{t:Dynkin-improve}.
The next lemma shows that  \eqref{e:UHK-rough} 
(hence Proposition \ref{p:UHK-rough}) is 
equivalent to the following weaker inequality: There exists $C>0$ such that
\begin{equation}\label{e:UHKDn}
	p(t,x,y) \le C \left(1 \wedge \frac{\delta_D(x)}{t^{1/\alpha}}\right)^{q}
	t^{-d/\alpha} \quad \text{for all 
		$t\le(\eta_0\wh{R})^{\alpha}$ and $x,y \in D$}.
\end{equation}

\begin{lemma}\label{l:(5.6)implies(5.1)}
	If \eqref{e:UHKDn} holds true,  then,  
	for any $T>0$, there exists  $C=C(T)>0$ such that \eqref{e:UHK-rough} also holds.
\end{lemma}
\pf
Without loss of generality, we assume $T >(\eta_0\wh{R})^{\alpha}$. 

\noindent
\textbf{Step 1:}
We claim that there exists a constant $c_1>0$ such that for all 
		$t\le (\eta_0\wh{R})^{\alpha}$ and $x,y \in D$,
\begin{equation}\label{e:UHK-claim}
	p(t,x,y) \le c_1 \left(1 \wedge \frac{\delta_D(x)}{t^{1/\alpha}}\right)^{q} \left( t^{-d/\alpha} \wedge \frac{t}{|x-y|^{d+\alpha}}\right).
\end{equation} 

If 
$\delta_D(x) \ge 2^{-4} t^{1/\alpha}$ or $|x-y| \le 8t^{1/\alpha}$, then \eqref{e:UHK-claim} follows from Proposition \ref{p:upper-heatkernel} or the assumption \eqref{e:UHKDn} respectively.

Assume now that  
$\delta_D(x)<2^{-4}t^{1/\alpha}$ and $|x-y|>8t^{1/\alpha}$. For \eqref{e:UHK-claim}, it suffices to show that
\begin{equation}\label{e:UHK-claim1}
	p(t,x,y) \le
	c_2\left(\frac{\delta_D(x)}{t^{1/\alpha}}\right)^q \frac{t}{|x-y|^{d+\alpha}}.
\end{equation}

Let $Q\in \partial D$ satisfy $|x-Q|=\delta_D(x)$,  let $V_1=B_D(Q, t^{1/\alpha})$, $V_3=\{w \in  D:|w-y| < |x-y|/2\}$ and $V_2=D \setminus (V_1 \cup V_3)$. 
By the 
triangle inequality, for any $u \in V_1$ and $w \in V_3$, we have
\begin{equation}\label{e:UHK-rough-1} 
	|u-w| \ge |x-y| - |x-u| - |y-w| \ge |x-y| - 2t^{1/\alpha} - \frac{|x-y|}{2} \ge \frac{|x-y|}{4} >2t^{1/\alpha}. 
\end{equation}
Let $\eps>0$ be such that $q-\alpha< \beta_1  - \eps$. By \eqref{e:B4-a}, \eqref{e:UHK-rough-1} and \eqref{e:Phi0-scaling},  for all $u \in V_1$ and $w \in V_3$,
\begin{equation*}
	\sB(u,w) \le 
	c_5 \Phi_0\left(\frac{\delta_D(u)\wedge \delta_D(w)}{|u-w|}\right)\le c_6\Phi_0\left(\frac{\delta_D(u)}{t^{1/\alpha}}\right)\le c_7\left(\frac{\delta_D(u)}{t^{1/\alpha}}\right)^{\beta_1-\eps}.
\end{equation*}
Thus, by Proposition \ref{p:bound-for-integral-new}, we obtain that
\begin{align}
	&\int_0^{ t }\int_{V_3} \int_{V_1} p^{V_1}(s, x, u) \sB(u,w) p(t-s, y,w) du dw ds 
	\label{e:UHK-rough111}\\
	&\le c_7
	\int_0^{t}\int_{V_1} p^{V_1}(s, x, u)\left(\frac{\delta_D(u)}{t^{1/\alpha}}\right)^{\beta_1-\eps} 
	\left( \int_{V_3}p(t-s, y,w) dw\right) ds \nn\\
	& \le \frac{c_7}{t^{(\beta_1-\eps)/\alpha}}
	\int_0^\infty\int_{V_1} p^{V_1}(s, x, u) \delta_D(u)^{\beta_1-\eps}  du  ds \nn\\
	&=
	\frac{c_7}{t^{(\beta_1-\eps)/\alpha}}\E_x	\int_0^{\tau_{V_1} } \delta_D(Y_s)^{\beta_1-\eps}  ds
	\le c_8 \left(\frac{\delta_D(x)}{t^{1/\alpha}}\right)^q t. \nn
\end{align}

It follows from  Theorem \ref{t:Dynkin-improve} (with $r=t^{1/\alpha}$) that
\begin{align}\label{l:lemma-upper-1-1}
	\P_x(\tau_{V_1}<t< \zeta)  \le \P_x(Y_{\tau_{V_1}} \in D)  
	\le c_3 \left(\frac{\delta_D(x)}{t^{1/\alpha}}\right)^q.
\end{align}
Note also that Proposition \ref{p:upper-heatkernel} implies that
\begin{equation}\label{e:UHK-rough-0}
	\sup_{s \le t, \, z \in V_2} p(s,z,y) \le c_4 \sup_{s \le t, \, z \in D, |z-y|\ge |x-y|/2} \frac{s}{|z-y|^{d+\alpha}} = 2^{d+\alpha}c_4 \frac{t}{|x-y|^{d+\alpha}}.
\end{equation}
Now combining \eqref{e:UHK-rough-1}--\eqref{e:UHK-rough-0} and  Lemma \ref{l:general-upper-2}, we get \eqref{e:UHK-claim1} and hence \eqref{e:UHK-claim}.

\noindent
\textbf{Step 2:} We first note that,      when $t_0:=(\eta_0\wh{R})^{\alpha}< t\le T$, 
by the semigroup property, symmetry and \eqref{e:UHKDn},
\begin{align*}
	&p(t,x,y)=\int_D p(t_0,x,z)p(t-t_0,y,z)\, dz\le C  \left(1 \wedge \frac{\delta_D(x)}{t_0^{1/\alpha}}\right)^{q}  t_0^{-d/\alpha}  \int_D p(t-t_0,y,z)\, dz\\
	&
	\le c_5  \left(1 \wedge \frac{\delta_D(x)}{t^{1/\alpha}}\right)^{q}  
	\le c_6  \left(1 \wedge \frac{\delta_D(x)}{t^{1/\alpha}}\right)^{q}  \left( t^{-d/\alpha} \wedge \frac{t}{|x-y|^{d+\alpha}}\right).
	\end{align*}
We have used the boundedness of $D$ in 	the last inequality. Therefore,   \eqref{e:UHK-claim} holds for all 
$t \in (0,  T]$. Using this and,  by the semigroup property, symmetry and \eqref{e:stableu2}, we conclude that for all $t \in (0,  T]$
\begin{align*}
	&p(t,x,y) = \int_{D} p(t/2,x,z)p(t/2,y,z)dz \\
	&\le c_2^2 \left(1 \wedge \frac{\delta_D(x)}{t^{1/\alpha}}\right)^{q}  \left(1 \wedge \frac{\delta_D(y)}{t^{1/\alpha}}\right)^{q} \int_{D} \left( (t/2)^{-d/\alpha} \wedge \frac{t/2}{|x-z|^{d+\alpha}}\right) \left( (t/2)^{-d/\alpha} \wedge \frac{t/2}{|y-z|^{d+\alpha}}\right)dz\nn\\
	&\le c_9 \left(1 \wedge \frac{\delta_D(x)}{t^{1/\alpha}}\right)^{q}  \left(1 \wedge \frac{\delta_D(y)}{t^{1/\alpha}}\right)^{q} \left( t^{-d/\alpha} \wedge \frac{t}{|x-y|^{d+\alpha}}\right).\nn
\end{align*}
This completes the proof.

\qed

To prove \eqref{e:UHKDn}, we will need the following simple lemma.
\begin{lemma}\label{l:lemma-upper-3}	Let  $B:=B_D(Q,r)$ for 	$Q\in \partial D$ and $r>0$. For all $\gamma\ge 0$, $a>0$, $t>0$ 
and $x \in B$, it holds that
	\begin{equation}\label{e:lemma-upper-3-1}
	  	\int_D p(t,x,z)\left(1 \wedge \frac{\delta_D(z)}{a}\right)^\gamma dz 
	  	\le \E_x \left[ \bigg(1 \wedge \frac{\delta_D(Y^{ B}_{t})}{a}\bigg)^\gamma : \tau_{B}>t\right] + \P_x(Y_{\tau_{B}} \in D).
	\end{equation}
	In particular, it holds that
	\begin{equation}\label{e:lemma-upper-3-2}
		\P_x(\zeta>t) \le t^{-1}\E_x\tau_{B} + \P_x(Y_{\tau_{B}}\in D).
	\end{equation}
\end{lemma}
\pf Since $Y^{B}_t=Y_t$ for $t<\tau_{B}$, we have
\begin{align*}
	&\int_D p(t,x,z)\left(1 \wedge \frac{\delta_D(z)}{a}\right)^\gamma dz 
	= \E_x \left[ \left(1 \wedge \frac{\delta_D(Y_{t})}{a}\right)^\gamma : t<\zeta \right]\\
	&=  \E_x \left[ \bigg(1 \wedge \frac{\delta_D(Y^{ B}_{t})}{a}\bigg)^\gamma : \tau_{B}>t \right] 
	+ \E_x \left[ \left(1 \wedge \frac{\delta_D(Y_{t})}{a}\right)^\gamma : \tau_{B}\le t<\zeta \right].
\end{align*}
Combining this with
\begin{align*}
	\E_x \left[ \left(1 \wedge \frac{\delta_D(Y_{t})}{a}\right)^\gamma : \tau_{B}\le t<\zeta \right] \le  \P_x( \tau_B \le t<\zeta) \le \P_x(Y_{\tau_B} \in D),
\end{align*}
we obtain \eqref{e:lemma-upper-3-1}.

For \eqref{e:lemma-upper-3-2}, by taking $\gamma=0$ in \eqref{e:lemma-upper-3-1} and using Markov's inequality, we get 
\begin{equation*}
	\P_x(\zeta>t)  \le \P_x(\tau_{B}>t) + \P_x(Y_{\tau_{B}}\in D)\le t^{-1}\E_x\tau_{B}+ \P_x (Y_{\tau_{B}} \in D).
\end{equation*}
\qed

We first prove that \eqref{e:UHKDn} holds in the case $q\in (0,\alpha)$.  

\begin{lemma}\label{l:UHKD_n} 
	If $q<\alpha$, then  \eqref{e:UHKDn} holds true. 
\end{lemma}
\pf 
By the semigroup property and Proposition \ref{p:upper-heatkernel},  for all $t\le (\eta_0\wh{R})^{\alpha}$ and $x,y\in D$,
\begin{align*}		&p(t,x,y) = \int_D p(t/2,x,z)p(t/2,z,y) dz\le c_1t^{-d/\alpha} \int_D p(t/2,x,z)dz   = c_1t^{-d/\alpha} \P_x(\zeta>t/2).\end{align*}
Thus, to obtain  \eqref{e:UHKDn}, it suffices to show that 
\begin{align}\label{e:UHKD_n-claim}
	\P_x(\zeta>t/2) \le c_1\left(1\wedge \frac{\delta_D(x)}{t^{1/\alpha}}\right)^q\quad \text{for all $t\le (\eta_0\wh{R})^{\alpha}$ and $x\in D$}.
\end{align}
Let $t\le (\eta_0\wh{R})^{\alpha}$ and $x\in D$. If $\delta_D(x)\ge 2^{-4}t^{1/\alpha}$, then
$$
\P_x(\zeta>t/2)\le 1\le 2^{4q}\left(1\wedge \frac{\delta_D(x)}{t^{1/\alpha}}\right)^q.
$$
Assume  $\delta_D(x)<2^{-4}t^{1/\alpha}$ and let $Q\in \partial D$ be such that $|x-Q|=\delta_D(x)$. Using \eqref{e:lemma-upper-3-2} (with $r=t^{1/\alpha}$), Proposition \ref{p:bound-for-integral-new} (with $\gamma=0>q-\alpha$),
and Theorem \ref{t:Dynkin-improve}, we obtain\begin{align*}
	\P_x(\zeta>t/2)&\le 2t^{-1}\E_x[\tau_{B_D(Q,  t^{1/\alpha})}]+\P_x(Y_{\tau_{B_D(Q, t^{1/\alpha})}}\in D)\\ & \le 2c_2 t^{-1}(t^{1/\alpha})^{\alpha-q} \delta_D(x)^q +c_2 \left(\frac{\delta_D(x)}{t^{1/\alpha}}\right)^q=c_3 \left(\frac{\delta_D(x)}{t^{1/\alpha}}\right)^q,
\end{align*} 
proving that \eqref{e:UHKD_n-claim} holds.  
The proof is complete.
\qed

In the next lemma, we remove the restriction $q<\alpha$ 
from Lemma \ref{l:UHKD_n}.

\begin{lemma}\label{l:UHKD}
	Inequality \eqref{e:UHKDn} holds true.
\end{lemma}
\pf Let $t\le (\eta_0\wh{R})^{\alpha}$.
By Lemma \ref{l:UHKD_n},  \eqref{e:UHKDn}  holds when  $q<\alpha$. Now, assume that \eqref{e:UHKDn} holds  
if $q <k\alpha $ for some $k\in \mathbb{N}$. We now show \eqref{e:UHKDn} also holds 
if $q\in [k\alpha, (k+1)\alpha)$.  By induction, it then follows that  \eqref{e:UHKDn} holds,
 without restriction.

Suppose $q\in [k\alpha,  (k+1)\alpha)$ and  fix $\eps \in ((2 \alpha-q-1)_+, (k+1)\alpha-q) \subset (0, \alpha)$. Since $(\alpha-1)_+ <q-\alpha + \eps<(k\alpha)\wedge (\alpha+\beta_1)$, using 
the strict-increasing property and positivity  of $p\mapsto C(p;\alpha, \F)$ on  $((\alpha-1)_+,\alpha+\beta_1)$, 
we have that $0<C( q-\alpha+\eps ;\alpha, \F)<C(q;\alpha,  \F)$.  Define 
$$
\wt \kappa(x):= \frac{C( q-\alpha+\eps ;\alpha, \F)}{C(q;\alpha, \F)}\kappa(x), \quad x\in D. 
$$
It is easy to see that $\wt \kappa$ satisfies  \hyperlink{K}{{\bf (K)}} with the same $C_1$ and $\kappa_0 = C( q-\alpha+\eps;\alpha, \F)$. Further, since $\wt \kappa(x) < \kappa(x)$ for all $x\in D$, we  get   $p(s,z,w)\le p^{\wt \kappa}(s,z,w)$ for all $s>0$ and $z,w\in D$. Applying the induction hypothesis to $p^{\wt \kappa}$, we obtain for all $t/2\le s\le t$ and $z,w \in D$,
\begin{align*}	p(s,z,w) &\le p^{\wt \kappa}(s,z,w) \le c_1 \left(1\wedge \frac{\delta_D(z)}{s^{1/\alpha}}\right)^{q-\alpha+\eps} s^{-d/\alpha} \le c_2 \left(1\wedge \frac{\delta_D(z)}{t^{1/\alpha}}\right)^{q-\alpha+\eps} t^{-d/\alpha} .\end{align*}
Therefore, by the semigroup property and symmetry, we get for all $x,y \in D$,
\begin{align*}
	&p(t,x,y) = 16t^{-2}\int_0^{t/4} \int_0^{t/4} p(t,x,y) dsdu\\
	&= 16t^{-2}\int_0^{t/4} \int_0^{t/4} \int_D\int_D p(s,x,z) p(t-s-u,z,w) p(u,y,w)dzdw dsdu\\
	& \le 16c_2 t^{-2-d/\alpha}\bigg(  \int_0^{t/4} \int_D p(s,x,z) \left(1\wedge \frac{\delta_D(z)}{t^{1/\alpha}}\right)^{q-\alpha+\eps}  dz ds \bigg)  \bigg(  \int_0^{t/4} \int_D p(u,y,w)   dw du \bigg)\\
	& \le 4c_2 t^{-1-d/\alpha} \int_0^{t/4} \int_D p(s,x,z) \left(1\wedge \frac{\delta_D(z)}{t^{1/\alpha}}\right)^{q-\alpha+\eps}  dz ds.
\end{align*}
Thus, to conclude \eqref{e:UHKDn} by induction, it suffices to show that there exists  $c_3>0$ such that 
\begin{equation}\label{e:UHKD-claim}
	\int_0^{t/4} \int_D p(s,v,z) \left(1\wedge \frac{\delta_D(z)}{t^{1/\alpha}}\right)^{q-\alpha+\eps}  dz ds 
	\le c_3t \left(1\wedge \frac{\delta_D(v)}{t^{1/\alpha}}\right)^q, \quad v \in D.
\end{equation}
Let $v\in D$. If $\delta_D(v)\ge 2^{-4}t^{1/\alpha}$, then
\begin{align*}
	&\int_0^{t/4} \int_D p(s,v,z) \left(1\wedge \frac{\delta_D(z)}{t^{1/\alpha}}\right)^{q-\alpha+\eps} dz ds\le \int_0^{t/4} \int_D p(s,v,z) dz ds\le \frac{t}{4}
	\le c_4 t \left(1\wedge \frac{\delta_D(v)}{t^{1/\alpha}}\right)^q.
\end{align*}
Assume $\delta_D(v)<2^{-4}t^{1/\alpha}$ and let $Q\in \partial D$ satisfy $|v-Q|=\delta_D(v)$. Write $B:=B_D(Q, t^{1/\alpha})$. Using \eqref{e:lemma-upper-3-1} in the first inequality below,  Fubini's theorem in the second, and Theorem \ref{t:Dynkin-improve} and Proposition \ref{p:bound-for-integral-new} in the fourth, we obtain
\begin{align*}
	&\int_0^{t/4} \int_D p(s,v,z) \left(1\wedge \frac{\delta_D(z)}{t^{1/\alpha}}\right)^{q-\alpha+\eps} dz ds\\
	& \le \int_0^{t/4} \E_v \left[ \bigg(1 \wedge \frac{\delta_D(Y^{B}_{s})}{t^{1/\alpha}}\bigg)^{q-\alpha+\eps} : \tau_{B}>s\right]  ds
	+ \int_0^{t/4} \P_v(Y_{\tau_{B}} \in D) ds\\
	& \le \E_v\int_0^{\tau_{B}}\bigg(1 \wedge \frac{\delta_D(Y^{B}_{s})}{t^{1/\alpha}}\bigg)^{q-\alpha+\eps}ds + \frac{t}{4}\P_v(Y_{\tau_{B}} \in D)\\ 
	&\le t^{-(q-\alpha+\eps)/\alpha}\E_v\int_0^{\tau_{B}}
	\delta_D(Y_{s})^{q-\alpha+\eps}ds+ \frac{t}{4}\P_v(Y_{\tau_{B}} \in D)\\
	&\le c_5 t \left(\frac{\delta_D(v)}{t^{1/\alpha}}\right)^q = c_5 t \left(1\wedge \frac{\delta_D(v)}{t^{1/\alpha}}\right)^q.
\end{align*}
This completes the proof of \eqref{e:UHKD-claim}, and thus of \eqref{e:UHKDn}. \qed

The proof of Proposition \ref{p:UHK-rough}  is now immediate.

\medskip

\noindent \textbf{Proof of Proposition \ref{p:UHK-rough}.} The assertion follows from Lemmas \ref{l:(5.6)implies(5.1)} and \ref{l:UHKD}. \qed

As a direct consequence of Proposition \ref{p:UHK-rough}, we obtain the following corollary.
\begin{corollary}\label{c:life}
	For any $T>0$,  there exists  $C=C(T) >0$ such that 
	$$
	\P_x (\zeta>t) \le  C \left(1 \wedge \frac{\delta_D(x)}{t^{1/\alpha}}\right)^q \quad \text{for all $t\in (0,T]$ and $x\in D$}.
	$$
\end{corollary}
\pf 
By Proposition \ref{p:UHK-rough} and \eqref{e:stableu1},
\begin{align*}
	& \P_x (\zeta>t) =\int_D p(t, x,y)dy \\
	&\le c_1  \left(1 \wedge \frac{\delta_D(x)}{t^{1/\alpha}}\right)^q \int_D \left( t^{-d/\alpha} \wedge \frac{t}{|x-y|^{d+\alpha}}\right)dy \le c_2 \left(1 \wedge \frac{\delta_D(x)}{t^{1/\alpha}}\right)^q.
\end{align*} 
\qed

We end this section with the following survival probability estimate.

\begin{lemma}\label{l:survival-exponential}
	There exist constants $C,C'>0$ such that for all $x\in D$, $r\in (0,1]$ and 
	$t\ge r^\alpha$, 
	\begin{align*}
		\P_x ( \tau_{B_{\overline D}(x,r)}\ge t) \le  C \left( 1 \wedge \frac{\delta_D(x)}{r}\right)^q e^{-C't/r^\alpha}.
	\end{align*}
\end{lemma}
\pf Let $x\in D$ and $r\in (0,1]$. Since $t\mapsto	\P_x ( \tau_{B_{\overline D}(x,r)}>t)$ is non-increasing, it suffices to show that there exist constants $c_1,c_2\ge 1$ independent of $x$ and $r$ such that 
\begin{align}\label{e:survival-exponential-claim}
	\P_x ( \tau_{B_{\overline D}(x,r)}\ge  (c_1 n+1)  r^\alpha) \le  c_2 \left( 1 \wedge \frac{\delta_D(x)}{r}\right)^q 2^{-n} \quad \text{for all $n\ge 0$}.
\end{align} 
It follows from Proposition \ref{p:E}(i) that $\sup_{z\in B_{\overline D}(x,r)} \E_z [\overline \tau_{B_{\overline D}(z,2r)} ]\le c_3r^\alpha$ for some constant $c_3>0$. Choose $c_1$ so that $c_3/c_1\le 1/2$. Then by Markov's inequality, 
\begin{align*}
	&\sup_{z\in B_{\overline D}(x,r)}	\P_z ( \tau_{B_{\overline D}(x,r)}\ge c_1r^\alpha) \le \sup_{z\in B_{\overline D}(x,r)}	\P_z ( \overline \tau_{B_{\overline D}(z,2r)}\ge c_1r^\alpha) \\
	&\le  (c_1r^\alpha)^{-1} 	\sup_{z\in B_{\overline D}(x,r)} \E_z [\overline \tau_{B_{\overline D}(z,2r)} ]  \le c_3/c_1 \le 1/2.
\end{align*}
Using this, and applying the Markov property at time $r^\alpha$ in the first step and at time $c_1r^\alpha$ in the subsequent steps, we get that for all $n\ge 0$, 
\begin{align*}
	& \P_x ( \tau_{B_{\overline D}(x,r)}\ge (c_1 n+1)r^\alpha)\le  \E_x\left[	\P_{Y_{r^\alpha}} ( \tau_{B_{\overline D}(x,r)}\ge c_1nr^\alpha) ;\,  \tau_{B_{\overline D}(x,r)}\ge r^\alpha \right]\\
	&\le  \P_x ( \tau_{B_{\overline D}(x,r)}\ge r^\alpha)  \sup_{z\in B_{\overline D}(x,r)}	\P_z ( \tau_{B_{\overline D}(x,r)}\ge c_1nr^\alpha) \\
	&\le \cdots \le  \P_x ( \tau_{B_{\overline D}(x,r)}\ge r^\alpha) 	\Big( \sup_{z\in B_{\overline D}(x,r)}	\P_z ( \tau_{B_{\overline D}(x,r)}\ge c_1r^\alpha)  \Big)^{n} \le 2^{-n} \P_x ( \zeta \ge r^\alpha) .
\end{align*}
Combining this with  Corollary \ref{c:life}, we obtain \eqref{e:survival-exponential-claim}. \qed


\section{Lower  bound}\label{s:lb}

In this section, we establish sharp lower bounds for the heat kernel. 
We emphasize that the assumption $\beta_1^*<\alpha +  \beta_1$ will not be used in this section.
We begin with the following interior diagonal estimate, which refines Proposition \ref{p:ndl}.

\begin{lemma}\label{l:ndl-refined}
	Let   $a \in (0,1]$. There exists  $C =C(a)>0$ 
	such that  for all $t>0$  and  $x,y \in D$ satisfying 
	$\delta_D(x) \ge a t^{1/\alpha}$ and $|x-y|< (t^{1/\alpha}/2)\wedge(\delta_D(x)/2^{3+1/\alpha}) $,
	\begin{equation}\label{e:ndl-refined-claim}
		p(t,x,y)  \ge  C t^{-d/\alpha}.
	\end{equation}
\end{lemma}
\pf Let $r:=\delta_D(x)/2$. Applying Proposition \ref{p:ndl} with $R_0=\diam(D)$ and $b=1/2$,  we get 
\begin{align}\label{e:ndl-refined-1}
	p(s,z,y) \ge 	p^{B(x,r)}(s,z,y) \ge  c_1s^{-d/\alpha} \quad \text{for all $s\le (r/2)^\alpha$ and $z,y\in B(x, s^{1/\alpha}/2)$.} 
\end{align}
If $t \le (r/2)^\alpha$, then  \eqref{e:ndl-refined-claim} follows from \eqref{e:ndl-refined-1}.

Assume $t>(r/2)^\alpha$. Set $t_0:=(r/2)^\alpha/2$ and let $n\ge 2$ be such that  $n < t/t_0 \le n+1$. Note that $|x-y| <r/2^{2+1/\alpha} = t_0^{1/\alpha}/2$ and
$n<t/t_0 \le (2r/a)^\alpha/t_0 =2 (4/a)^\alpha$. Using the semigroup property and \eqref{e:ndl-refined-1}, we obtain
\begin{align*}
	p(t,x,y) &= \int_{B(x,t_0^{1/\alpha}/2)\times \cdots B(x,t_0^{1/\alpha}/2)} p(t_0,x,z_1) \cdots p(t_0,z_{n-2},z_{n-1})  \\
	&\qquad \qquad \qquad \qquad\qquad \qquad  \times p(t - (n-1)t_0, z_{n-1}, y )\,  dz_1\cdots dz_{n-1}\\
	&\ge c_1(2t_0)^{-d/\alpha}(c_1 t_0^{-d/\alpha})^{n-1}  (c_2(t_0^{1/\alpha}/2)^d)^{n-1} =c_3 c_4^n t_0^{-d/\alpha} \ge c_3 (c_4\wedge 1)^{2(4/a)^\alpha} t^{-d/\alpha}.
\end{align*}
\qed

\subsection{Some estimates related to lifted points}\label{ss:lifted-points}

Recall that $\eta_1$ and $\eta_2$ are the constants in \eqref{e:lifting-property}.  
For $x\in D$ and $t>0$, we define the lifted point\begin{equation}\label{e:def-x(t)}	x(t):= x+ (t^{1/\alpha}/2)\bn_x.\end{equation}
Note that, by \eqref{e:lifting-property}, we have  $x(t)\in D$ for all $x\in D$ and $t\in (0, (2\eta_1)^\alpha]$.

The following simple consequence of \eqref{e:lifting-property} will be used several times in the sequel: 
For all  $x\in D$ and $t\in (0,\eta_1^\alpha]$, it holds that
\begin{align}\label{e:lifting-standard}
	\delta_D(y)  >2^{-2}\eta_2(t/4)^{1/\alpha} \quad \text{for all  $y \in B(x(t/4), 2^{-2}\eta_2(t/4)^{1/\alpha})$.} 
\end{align}
In fact, 
$$
\delta_D(y) \ge \delta_D(x(t/4))-|x(t/4)-y|  > 2^{-1}\eta_2(t/4)^{1/\alpha}-2^{-2}\eta_2(t/4)^{1/\alpha}=2^{-2}\eta_2(t/4)^{1/\alpha}.
$$
In particular, for all $x\in D$ and $t\in (0,\eta_1^\alpha]$, we have 
\begin{align}\label{e:lifting-standard1}
 B(x(t/4), 2^{-2}\eta_2(t/4)^{1/\alpha}) \subset D.
 \end{align}

\begin{lemma}\label{l:lift-interior-estimate}
	Let $t\in (0,\eta_1^\alpha]$ and   $x,y \in D$. 
	There exists $C>0$ independent of $t,x,y$ such that
	\begin{align*}
		p(t/2, z,w) \ge  A_{\Phi_1,\Phi_2,\ell}(t,x,y) 
		\left( t^{-d/\alpha} \wedge \frac{t}{|x-y|^{d+\alpha}}\right) 
	\end{align*}
	for all	$z\in B(x(t/4), 2^{-6}\eta_2(t/8)^{1/\alpha})$ and 	$w\in B(y(t/4), 2^{-6}\eta_2(t/8)^{1/\alpha})$.
\end{lemma}
\pf 
Without loss of generality, we  assume $\delta_D(x)\le \delta_D(y)$. Write $r:=2^{-6}\eta_2(t/8)^{1/\alpha}$. Let  $z\in B(x(t/4),r)$ and 	$w\in B(y(t/4), r)$.
Recall from \eqref{e:lifting-standard} that  $\delta_D(z) \wedge \delta_D(w) > 2^{-2}\eta_2(t/4)^{1/\alpha}$, which implies that, for
$s\in [t/6,t/2]$, we have  $\delta_D(z) \wedge \delta_D(w) > 2^{-2-1/\alpha}\eta_2s^{1/\alpha}$ and 
$
2r =2^{-5}\eta_2(t/8)^{1/\alpha}
< 2^{-5}\eta_2s^{1/\alpha} \wedge 2^{-3-1/\alpha} \delta_D(y).$
Thus, by Lemma \ref{l:ndl-refined}, it holds that 
\begin{align}\label{e:lift-interior-estimate-1}
	\bigg(\inf_{s\in [t/6,t/2],\, |z-z'|<2r}	p(s, z, z') \bigg) \wedge  \bigg( \inf_{s\in [t/6,t/2],\, |w-w'|<2r}	p(s,  w,w') \bigg) \ge c_1 t^{-d/\alpha}.
\end{align}
If $|z-w|<2r$, then the result follows from \eqref{e:lift-interior-estimate-1} and \eqref{e:upper-bound-A}.

Suppose $|z-w|\ge 2r$. For all $z'\in B(z,r/2)$, by \eqref{e:lifting-property}, 
$$
\delta_D(z') \ge \delta_D(x(t/4)) -3r/2 > \eta_2\delta_D(x) + 2^{-2}\eta_2(t/4)^{1/\alpha} \ge 2^{-2-2/\alpha}\eta_2(\delta_D(x) \vee t^{1/\alpha})
$$ 
and  $\delta_D(z') \le \delta_D(x(t/4)) + 3r/2 \le \delta_D(x) + (t/4)^{1/\alpha} \le \delta_D(x) \vee t^{1/\alpha}$. Similarly,  we have that $2^{-2-2/\alpha}\eta_2(\delta_D(y) \vee t^{1/\alpha})\le \delta_D(w') \le  \delta_D(y)\vee t^{1/\alpha}$ for all 
$w'\in B(w,r/2)$. Further, for all $z'\in B(z,r/2)$ and $w'\in B(w,r/2)$,
\begin{align*}
	|z'-w'| <|x(t/4)-y(t/4)|+3r\le  |x-y| + |x-x(t/4)|+ |y-y(t/4)| +3r \le |x-y|+c_2t^{1/\alpha}
\end{align*}
and $|z'-w'| >|x-y|- c_3t^{1/\alpha}$. 
Hence, $|z'-w'|\asymp t^{1/\alpha}$ if $|x-y| \le 2c_3t^{1/\alpha}$, and $|z'-w'|\asymp |x-y|$ if $|x-y| > 2c_3t^{1/\alpha}$. Thus, for all  $z'\in B(z,r/2)$ and $w'\in B(w,r/2)$, using  \hyperlink{A1}{{\bf (A1)}},   \eqref{e:Phi1-scaling}--\eqref{e:ell-scaling} and \eqref{e:upper-bound-A}, we obtain that  if $|x-y| \le 2c_3t^{1/\alpha}$, then
\begin{align}\label{e:lift-interior-estimate-2}
	\frac{\sB(z',w')}{|z'-w'|^{d+\alpha}}& \ge \frac{c_4}{t^{(d+\alpha)/\alpha}} \Phi_1 \bigg( \frac{\delta_D(x) \vee t^{1/\alpha}}{t^{1/\alpha}}\bigg)  \Phi_2 \bigg( \frac{\delta_D(y) \vee t^{1/\alpha}}{t^{1/\alpha}}\bigg) \ell \bigg( \frac{\delta_D(x) \vee t^{1/\alpha}}{t^{1/\alpha}}\bigg)  \nn\\
	&=  \frac{c_4}{t^{(d+\alpha)/\alpha}} \ge  \frac{c_5A_{\Phi_1,\Phi_2,\ell}(t,x,y)}{t^{(d+\alpha)/\alpha}},
\end{align}
and if $|x-y| > 2c_3t^{1/\alpha}$, then
\begin{align}\label{e:lift-interior-estimate-3}
	\frac{\sB(z',w')}{|z'-w'|^{d+\alpha}} \ge \frac{c_6A_{\Phi_1,\Phi_2,\ell}(t,x,y)}{|x-y|^{d+\alpha}}.
\end{align}
Using the strong Markov property in the first line below, \eqref{e:lift-interior-estimate-1} in the second,  the L\'evy system formulas \eqref{e:levy-system-Y-bar}--\eqref{e:levy-system-Y-kappa} in the third,  and  \eqref{e:lift-interior-estimate-2}--\eqref{e:lift-interior-estimate-3} in the fourth, we deduce that 
\begin{align}\label{e:lift-interior-estimate-4}
	p(t/2,z,w) &\ge \E_z \bigg[ p(t/2- \tau_{B(z,r/2)}, Y_{\tau_{B(z,r/2)}}, w);\,  \tau_{B(z,r/2)} \le t/3, \, Y_{\tau_{B(z,r/2)}} \in B(w,r/2)  \bigg] \\
	&\ge c_1t^{-d/\alpha} \P_z \Big( \tau_{B(z,r/2)} \le t/3, \, Y_{\tau_{B(z,r/2)}} \in B(w,r/2)  \Big)\nn\\
	&\ge c_1t^{-d/\alpha} \E_z \bigg[ \int_{0}^{\tau_{B(z,r/2)} \wedge (t/3)} \int_{B(w,r/2)}\frac{\sB(Y_s, w')}{|Y_s-w'|^{d+\alpha}} dw\,ds \bigg]\nn\\
	&\ge \frac{c_7 A_{\Phi_1,\Phi_2,\ell}(t,x,y)}{(|x-y| \vee t^{1/\alpha})^{d+\alpha}} \E_z [\tau_{B(z,r/2)} \wedge (t/3)].\nn
\end{align}
By applying Proposition \ref{p:ndl}(ii) with $b=1/2$, we get that 
\begin{align}\label{e:lift-interior-estimate-5}
	&	\E_z [\tau_{B(z,r/2)} \wedge (t/3)] \ge (r/4)^\alpha \P_z(\tau_{B(z,r/2)} \ge (r/4)^\alpha ) \ge (r/4)^\alpha \P_z\big(Y^{B(z,r/2)}_{(r/4)^\alpha} \in B(z, r/8) \big) \\
	&= (r/4)^\alpha \int_{B(z,r/8)} p^{B(z,r/2)}((r/4)^\alpha, z, z')dz' \ge c_8 (r/4)^{\alpha-d}\int_{B(z,r/8)}dz'  = c_9r^\alpha = c_{10}t.\nn
\end{align}
Combining \eqref{e:lift-interior-estimate-4} and  \eqref{e:lift-interior-estimate-5} with the fact that $t/(|x-y|\vee t^{1/\alpha})^{d+\alpha} = t^{-d/\alpha} \wedge (t|x-y|^{-d-\alpha})$, we obtain the desired result. \qed

\begin{lemma}\label{l:lift-probability-pre}
	If $q>0$, then there
	exists $C>0$  such that for all $x\in D$, $t\in (0,1]$ and $\eps\in (0,\wh R/48]$,
	\begin{align*}
		\E_x \left[\int_0^{\tau_{B_{ D}(x, \eps t^{1/\alpha})}} \Phi_0\bigg(\frac{\delta_D(Y_s)}{t^{1/\alpha}}\bigg)\, ds\right] \ge C \eps^{\alpha+\beta_1^*+1}t\left( 1 \wedge \frac{\delta_D(x)}{\eps t^{1/\alpha}}\right)^q.
	\end{align*}
\end{lemma}
\pf Assume first that $\delta_D(x)\ge \eps t^{1/\alpha}/4$. Then,  for all $z\in B_D(x,\eps t^{1/\alpha}/8)$,  by the almost increasing property of $\Phi_0$ and  \eqref{e:Phi0-scaling}, we have $\Phi_0(\delta_D(z)/t^{1/\alpha}) \ge c_1\Phi_0 (\eps/8) \ge c_2(\eps/8)^{\beta_1^*+1}$. Hence, using Proposition \ref{p:E}(ii), we obtain
\begin{align*}
	&\E_x \left[\int_0^{\tau_{B_{ D}(x, \eps t^{1/\alpha})}} \Phi_0\bigg(\frac{\delta_D(Y_s)}{t^{1/\alpha}}\bigg)\, ds\right] \ge 	\E_x \left[\int_0^{\tau_{B(x, \eps t^{1/\alpha}/8)}} \Phi_0\bigg(\frac{\delta_D(Y_s)}{t^{1/\alpha}}\bigg)\, ds\right]\\
	&\ge  c_2(\eps/8)^{\beta_1^*+1} \E_x[\tau_{B(x, \eps t^{1/\alpha}/8)}] \ge c_3 \eps^{\alpha+\beta_1^*+1}t.
\end{align*}

Suppose $\delta_D(x)<\eps t^{1/\alpha}/4$ and let $Q_x\in \partial D$ satisfy $|x-Q_x|=\delta_D(x)$. Note that $2\eps t^{1/\alpha} \le \wh R/24$, $B_D(Q_x, \eps t^{1/\alpha}/2)\subset B_D(x,\eps t^{1/\alpha}) \subset B_D(Q_x, 2\eps t^{1/\alpha})$ and $x\in B_D(Q_x, \eps t^{1/\alpha}/4)$. Hence, using \eqref{e:Phi0-scaling} and Proposition \ref{p:bound-for-integral-new} with $\gamma=\beta_1^*+1>q-\alpha$, we get that
\begin{align*}
	\E_x \bigg[\int_0^{\tau_{B_{ D}(x, \eps t^{1/\alpha})}} \Phi_0\bigg(\frac{\delta_D(Y_s)}{t^{1/\alpha}}\bigg)\, ds\bigg] & \ge c_4	\E_x \bigg[\int_0^{\tau_{B_{ D}(x, \eps t^{1/\alpha})}} \bigg(\frac{\delta_D(Y_s)}{t^{1/\alpha}}\bigg)^{\beta_1^*+1}\, ds\bigg] \\
	&\ge \frac{c_5 (\eps t^{1/\alpha})^{\alpha+\beta_1^*+1-q} \delta_D(x)^q}{t^{(\beta_1^*+1)/\alpha}}   = c_5 \eps^{\alpha+\beta_1^*+1}t\bigg(\frac{\delta_D(x)}{\eps t^{1/\alpha}}\bigg)^q.
\end{align*}
\qed

\begin{lemma}\label{l:lift-probability}
If $q>0$, then there
 exist constants $\lambda \in (0,
2^{-2-3/\alpha} \wedge (\wh R/48)]$ and  $C>0$ such that for all $x\in D$ and  $t\in (0, \eta_1^\alpha]$,
\begin{align*}
	\P_x \Big( Y_{\tau_{B_{ D}(x, \lambda t^{1/\alpha})} \wedge (t/8)} \in B(x(t/4), 2^{-6}\eta_2(t/8)^{1/\alpha}) \Big)  \ge C\left( 1 \wedge \frac{\delta_D(x)}{t^{1/\alpha}}\right)^q.
\end{align*}	
\end{lemma}
\pf  
Let 
$\lambda\in (0,2^{-1-1/\alpha} \wedge (\wh R/48)]$ be a constant to be chosen later.   
Since $\lambda t^{1/\alpha}  + 2^{-6}\eta_2(t/8)^{1/\alpha} <(17/64)(t/4)^{1/\alpha} =(17/32) |x-x(t/4)| $, it holds that
\begin{align}\label{e:lift-probability-1}
		|z-y| \asymp t^{1/\alpha}  \quad \text{for  $z\in B_{ D}(x, \lambda t^{1/\alpha})$ and $y \in B(x(t/4), 2^{-6}\eta_2(t/8)^{1/\alpha})$.} 
\end{align}
 Thus, 
combining  \hyperlink{A1}{{\bf (A1)}}, \eqref{e:lifting-standard},  \eqref{e:Phi1-scaling}--\eqref{e:ell-scaling}  and \eqref{e:def-Phi0} with the fact that $\Phi_0(r)=1$ for all $r\ge 1$, we get that for  $z\in B_{ D}(x, \lambda t^{1/\alpha})$ and $y \in B(x(t/4), 2^{-6}\eta_2(t/8)^{1/\alpha})$,
\begin{align}\label{e:lift-probability-2}
\sB(z,y) &\asymp \begin{cases}
	\Phi_0(\delta_D(z)/t^{1/\alpha}) &\mbox{ if $\delta_D(z)\le 2^{-2}\eta_2(t/4)^{1/\alpha}$},\\
	1 &\mbox{ if $\delta_D(z)> 2^{-2}\eta_2(t/4)^{1/\alpha}$}
\end{cases} \, \asymp  \Phi_0\bigg(\frac{\delta_D(z)}{t^{1/\alpha}}\bigg).
\end{align}
Using the L\'evy system formulas \eqref{e:levy-system-Y-bar}--\eqref{e:levy-system-Y-kappa},   \eqref{e:lift-probability-1} and \eqref{e:lift-probability-2},  we obtain
\begin{align}\label{e:lift-probability-3}
		&\P_x \Big( Y_{\tau_{B_{ D}(x, \lambda t^{1/\alpha})} \wedge (t/8)} \in B(x(t/4), 2^{-6}\eta_2(t/8)^{1/\alpha})  \Big)\\
		  &= \E_x \bigg[\int_0^{\tau_{B_{ D}(x, \lambda t^{1/\alpha})} \wedge (t/8) }\int_{B(x(t/4), 2^{-6}\eta_2(t/8)^{1/\alpha}) }  \frac{\sB(Y_s, y)}{|Y_s-y|^{d+\alpha}}\,dy\, ds\bigg] \nn\\
		  &\ge c_1 t^{-d/\alpha-1}\E_x \bigg[\int_0^{\tau_{B_{ D}(x, \lambda t^{1/\alpha})} \wedge (t/8) }\Phi_0\bigg(\frac{\delta_D(Y_s)}{t^{1/\alpha}}\bigg) ds\int_{B(x(t/4), 2^{-6}\eta_2(t/8)^{1/\alpha}) } \,dy\bigg] \nn\\
		  & \ge c_2 t^{-1}\E_x \bigg[\int_0^{\tau_{B_{ D}(x, \lambda t^{1/\alpha})}  } \Phi_0\bigg(\frac{\delta_D(Y_s)}{t^{1/\alpha}}\bigg)\, ds ;\, \tau_{B_{ D}(x, \lambda t^{1/\alpha})}\le t/8 \bigg] \nn\\
		  &\ge c_2 t^{-1}\E_x \left[\int_0^{\tau_{B_{ D}(x, \lambda t^{1/\alpha})}} \Phi_0\bigg(\frac{\delta_D(Y_s)}{t^{1/\alpha}}\bigg)\, ds\right] - c_3t^{-1} \E_x \big[ \tau_{B_{ D}(x, \lambda t^{1/\alpha})}  ;\, \tau_{B_{ D}(x, \lambda t^{1/\alpha})} > t/8 \big].\nn
\end{align}
We used the boundedness of $\Phi_0$ in the last inequality.  Note that  $\lambda^\alpha<1/8$.    Using Lemma \ref{l:survival-exponential}  and the fact that $ \sup_{a\ge 0}a^ke^{-a} = (k/e)^k$ for all $k>0$, we obtain
\begin{align}\label{e:lift-probability-4}
	&t^{-1}\E_x \big[ \tau_{B_{ D}(x, \lambda t^{1/\alpha})}  ;\, \tau_{B_{ D}(x, \lambda t^{1/\alpha})} > t/8 \big]\\
	& \le \sum_{n=1}^\infty \frac{n}{4}\P_x \big( \tau_{B_{ D}(x, \lambda t^{1/\alpha})} \in (nt/8, (n+1)t/8]\big)
\le c_4  \left( 1 \wedge \frac{\delta_D(x)}{\lambda t^{1/\alpha}}\right)^q\sum_{n=1}^\infty n e^{-c_5n/\lambda^\alpha} \nn\\
	&\le c_6 \lambda^{2\alpha+\beta_1^*+1} \left( 1 \wedge \frac{\delta_D(x)}{\lambda t^{1/\alpha}}\right)^q \sum_{n=1}^\infty n^{-(\alpha+\beta_1^*+1)/\alpha} = c_7 \lambda^{2\alpha+\beta_1^*+1} \left( 1 \wedge \frac{\delta_D(x)}{\lambda t^{1/\alpha}}\right)^q,\nn
\end{align} 
where $c_7$ is independent of $\lambda$. Combining \eqref{e:lift-probability-3} with Lemma \ref{l:lift-probability-pre} and  \eqref{e:lift-probability-4}, we arrive at
\begin{align*}
	\P_x \Big( Y_{\tau_{B_{ D}(x, \lambda t^{1/\alpha})} \wedge (t/8)} \in B(x(t/4), 2^{-6}\eta_2(t/8)^{1/\alpha})  \Big) \ge  (c_8\lambda^{\alpha+\beta_1^*+1}-c_7 \lambda^{2\alpha+\beta_1^*+1})\left( 1 \wedge \frac{\delta_D(x)}{\lambda t^{1/\alpha}}\right)^q.
\end{align*}
By choosing $\lambda$ sufficiently small, this yields the desired result.
\qed

\begin{lemma}\label{l:lift}
	There exists $C>0$ such that  for all $x\in D$,  $t\in (0, \eta_1^\alpha]$ and $y \in B(x(t/4),$ $ 2^{-6}\eta_2(t/8)^{1/\alpha})$, 
	$$ 	p(t/4,x,y) \ge C  \left(1 \wedge \frac{\delta_D(x)}{t^{1/\alpha}}\right)^{q}t^{-d/\alpha}.
	$$
\end{lemma}
\pf 
When $q=0$, since $|x-y|
\le |x-x(t/4)|+|x(t/4)-y|<
(t/4)^{1/\alpha}/2 + 2^{-6}\eta_2(t/8)^{1/\alpha} <(3/4)(t/4)^{1/\alpha}$,  applying Proposition \ref{p:ndl} with $r=R_0=4^{1+1/\alpha}\eta_1$ and $b=3/4$, we get the result.

Suppose $q>0$. 
 Note that,  by \eqref{e:lifting-standard}, 
for all $s\in [t/8,t/4]$ and $z\in B(x(t/4),2^{-6}\eta_2(t/8)^{1/\alpha})$, 
$\delta_D(y) > 2^{-2}\eta_2(t/4)^{1/\alpha} \ge 2^{-2}\eta_2s^{1/\alpha}$ and 
$$
|y-z| \le |y- x(t/4)|+|x(t/4)-z| <2^{-5}\eta_2(t/8)^{1/\alpha} 
\le 2^{-5}\eta_2s^{1/\alpha} \wedge 2^{-3-1/\alpha} \delta_D(y)$$
Thus,  we can apply Lemma \ref{l:ndl-refined} with $a=2^{-2}\eta_2$ to obtain 
\begin{align}\label{e:lift-on-diagonal-1}
	p(s,z,y) \ge  c_2 t^{-d/\alpha}
\end{align}
for all $s\in [t/8,t/4]$ and $z\in B(x(t/4),2^{-6}\eta_2(t/8)^{1/\alpha})$.  Let $\lambda$ be the constant in Lemma \ref{l:lift-probability} and write $B:=B_D(x,\lambda t^{1/\alpha})$.
Using the strong Markov property, \eqref{e:lift-on-diagonal-1} and Lemma \ref{l:lift-probability}, we obtain
\begin{align*}
&	p(t/4,x,y) \ge 	\E_x\left[p(t/4- \tau_{B}, Y_{\tau_B}, y);\, \tau_{B} \le t/8,\, Y_{\tau_B} \in B(x(t/4), 2^{-6}\eta_2(t/8)^{1/\alpha})\right]\nn\\
	&\ge \bigg(
	\inf_{  s\in [t/8,t/4], \, z \in B(x(t/4), 2^{-6}\eta_2(t/8)^{1/\alpha})} p (s, z,y) 	\bigg)	\P_x \Big( Y_{\tau_{B} \wedge (t/8)} \in B(x(t/4), 2^{-6}\eta_2(t/8)^{1/\alpha}) \Big)  \nn\\
	&\ge c_3  \left(1 \wedge \frac{\delta_D(x)}{t^{1/\alpha}}\right)^{q}t^{-d/\alpha}.
\end{align*}
\qed

\subsection{Sharp lower bound}\label{ss:slb}

\begin{lemma}\label{l:HKE-lower-0}
There exists $C>0$ such that for all $t\in (0, \eta_1^\alpha]$ and  $x,y \in D$, 
	\begin{align*}
		p(t,x,y) \ge C \left(1 \wedge \frac{\delta_D(x)}{t^{1/\alpha}}\right)^q 
		\left(1 \wedge \frac{\delta_D(y)}{t^{1/\alpha}}\right)^q A_{\Phi_1,\Phi_2,\ell}(t,x,y) 
		\left( t^{-d/\alpha} \wedge \frac{t}{|x-y|^{d+\alpha}}\right).
	\end{align*}
\end{lemma}
\pf By the semigroup property, Lemmas \ref{l:lift} and \ref{l:lift-interior-estimate}, we obtain
\begin{align*}
	p(t,x,y)& \ge \int_{B(x(t/4), 2^{-6}\eta_2(t/8)^{1/\alpha})}  \int_{B(y(t/4), 2^{-6}\eta_2(t/8)^{1/\alpha})} p(t/4,x,z)p(t/2,z,w)p(t/4,w,y) dwdz\\
	&\ge c_1t^{-2d/\alpha} \left(1 \wedge \frac{\delta_D(x)}{t^{1/\alpha}}\right)^q 
	\left(1 \wedge \frac{\delta_D(y)}{t^{1/\alpha}}\right)^q  A_{\Phi_1,\Phi_2,\ell}(t,x,y)
	\left( t^{-d/\alpha} \wedge \frac{t}{|x-y|^{d+\alpha}}\right)\\
	&\quad \times \int_{B(x(t/4), 2^{-6}\eta_2(t/8)^{1/\alpha})} dz \int_{B(y(t/4), 2^{-6}\eta_2(t/8)^{1/\alpha})} dw\\
	&= c_2 \left(1 \wedge \frac{\delta_D(x)}{t^{1/\alpha}}\right)^q
	\left(1 \wedge \frac{\delta_D(y)}{t^{1/\alpha}}\right)^q A_{\Phi_1,\Phi_2,\ell}(t,x,y)
	\left( t^{-d/\alpha} \wedge \frac{t}{|x-y|^{d+\alpha}}\right).
\end{align*}
\qed

\begin{lemma}\label{l:HKE-lower-1}
	For any $T>0$ and $\eps>0$, there exists $C=C(T,\eps)>0$ such that for all $t\in (0,T]$ and $x,y\in D$ with $(\delta_D(x) \wedge \delta_D(y)) \vee t^{1/\alpha} \ge \eps|x-y|$,
	\begin{align*}
		p(t,x,y) \ge C \left(1 \wedge \frac{\delta_D(x)}{t^{1/\alpha}}\right)^q 
		\left(1 \wedge \frac{\delta_D(y)}{t^{1/\alpha}}\right)^q \bigg( t^{-d/\alpha} \wedge \frac{t}{|x-y|^{d+\alpha}}\bigg).
	\end{align*}
\end{lemma}
\pf If $t\le \eta_1^\alpha$, then the assertion follows from 
Lemma \ref{l:HKE-lower-0} and \eqref{e:interior-lower-bound-A}. 

  Suppose $T>\eta_1^\alpha$ and $t\in (\eta_1^\alpha,T]$. 
  Then  $ t^{-d/\alpha} \wedge \frac{t}{|z-w|^{d+\alpha}} \asymp 1$ for all $z, w\in D$. 
  Let $B:=B(x(\eta_1^\alpha),  2^{-6-1/\alpha}\eta_1\eta_2)$. By \eqref{e:lifting-property}, we have $\delta_D(z)\ge \eta_2\delta_D(x)+ \eta_1\eta_2/4$ for all $z\in B$. Applying Lemma \ref{l:ndl-refined} with $a=\eta_1\eta_2/(4T^{1/\alpha})$, we get
 \begin{align}\label{e:HKE-lower-1}
 	p(t-\eta_1^\alpha/2,z,w) \ge c_1 (t-\eta_1^\alpha/2)^{-d/\alpha} \ge c_1T^{-d/\alpha} \quad \text{for all $z,w\in B$}.
 \end{align}
 On the other hand, for all $z\in B$,
 we have  $|y-z|\le |x-y|+|x-z|\le |x-y|+\eta_1$,
 \begin{align*}
 	(\delta_D(x) \wedge \delta_D(z)) \vee \frac{\eta_1}{4^{1/\alpha}}& \ge \frac{\eta_1}{4^{1/\alpha}} \ge \frac{|x-z|}{4^{1/\alpha}}
 \end{align*}
 and
 \begin{align*}
 	(\delta_D(y) \wedge \delta_D(z)) \vee \frac{\eta_1}{4^{1/\alpha}} & \ge  (\delta_D(y) \wedge  (\eta_2\delta_D(x) + \frac{\eta_1\eta_2}{4}))  \vee \frac{\eta_1}{4^{1/\alpha}} \\
 	&\ge \frac{\eta_1\eta_2}{4^{1/\alpha}T^{1/\alpha}}((\delta_D(x) \wedge \delta_D(y)) \vee t^{1/\alpha})\ge \frac{\eta_1\eta_2  }{4^{1/\alpha}T^{1/\alpha}} ((\eps|x-y|) \vee \eta_1)\\
 	& \ge \frac{\eta_1\eta_2 \eps }{2\cdot 4^{1/\alpha}T^{1/\alpha}} (|x-y| + \eta_1) \ge \frac{\eta_1\eta_2 \eps }{2\cdot 4^{1/\alpha}T^{1/\alpha}} |y-z|.
 \end{align*}
Hence, for all $z\in B$, by 
Lemma \ref{l:HKE-lower-0} and \eqref{e:interior-lower-bound-A}  with the fact that  $\delta_D(z) \ge \eta_1\eta_2/4$,  we get
 \begin{align}\label{e:HKE-lower-2}
 	p(\eta_1^\alpha/4,x,z) \ge c_2\bigg(1\wedge \frac{\delta_D(x)}{\eta_1}\bigg)^q 
	\quad \text{and} \quad
	p(\eta_1^\alpha/4,z,y) &\ge  c_2\bigg(1\wedge \frac{\delta_D(y)}{\eta_1}\bigg)^q.
 \end{align} 
 Using the semigroup property, \eqref{e:HKE-lower-1} and \eqref{e:HKE-lower-2}, we arrive at
 \begin{align*}
 	p(t,x,y) &\ge \int_B \int_B 	p(\eta_1^\alpha/4,x,z) 	p(t-\eta_1^\alpha/2,z,w) p(\eta_1^\alpha/4,w,y) dzdw\\
 	&\ge 	c_2^2
	\bigg(1\wedge \frac{\delta_D(x)}{\eta_1}\bigg)^q\bigg(1\wedge \frac{\delta_D(y)}{\eta_1}\bigg)^q  \int_B \int_B dzdw \\
 	&\asymp 
	\left(1 \wedge \frac{\delta_D(x)}{t^{1/\alpha}}\right)^q
 	\left(1 \wedge \frac{\delta_D(y)}{t^{1/\alpha}}\right)^q \bigg( t^{-d/\alpha} \wedge \frac{t}{|x-y|^{d+\alpha}}\bigg),
 \end{align*} 
 where the comparison constants in the last line depend only on $T$ and $\eps$. 
 The proof is complete.
 
  \qed

 Recall that the constant $\eps_1$ is defined in \eqref{e:def-eps1}.
\begin{lemma}\label{l:HKE-lower-two-jumps}
	There exists  $C>0$ such that for all $t\in (0,\eta_1^\alpha]$ and $x,y \in D$,  
	\begin{align*}
		p(t,x,y) &\ge C \left(1 \wedge \frac{\delta_D(x)}{t^{1/\alpha}}\right)^q \left(1 \wedge \frac{\delta_D(y)}{t^{1/\alpha}}\right)^q
		\left( t^{-d/\alpha} \wedge \frac{t}{|x-y|^{d+\alpha}}\right)\\
		&\quad \times \big(t\wedge  |x-y|^{\alpha} \big) 
		\int_{ t^{1/\alpha}\wedge (\eps_1|x-y|/2)}^{\eps_1|x-y|}A_{\Phi_1,\Phi_2,\ell}(t,x,x+u\mathbf{n}_x)\,A_{\Phi_1,\Phi_2,\ell}(t,x+u\mathbf{n}_x,y)\frac{du}{u^{\alpha+1}}.
	\end{align*}
\end{lemma}
\pf  Let $r:=\eps_1|x-y|$. 
Note that, by the definition of $\eps_1$,  we have $r\le \epsilon_1 {\rm diam} (D) \le \eta_1$.
If $ t^{1/\alpha}\ge r/2$, then by Lemma \ref{l:HKE-lower-1}, 
$$
p(t,x,y) \ge  c_1 \left(1 \wedge \frac{\delta_D(x)}{t^{1/\alpha}}\right)^q \left(1 \wedge \frac{\delta_D(y)}{t^{1/\alpha}}\right)^q t^{-d/\alpha}.
$$
It follows from \eqref{e:upper-bound-A} that
\begin{align*}
	&\big(t\wedge  |x-y|^{\alpha} \big) 
	\int_{ t^{1/\alpha}\wedge (r/2)}^{r}A_{\Phi_1,\Phi_2,\ell}(t,x,x+u\mathbf{n}_x)\,A_{\Phi_1,\Phi_2,\ell}(t,x+u\mathbf{n}_x,y)\frac{du}{u^{\alpha+1}}\\
	&\le  c_2|x-y|^\alpha 	\int_{r/2}^{r}\frac{du}{u^{\alpha+1}} = c_3.
\end{align*}
Thus the assertion is valid in this case.

Suppose  $t^{1/\alpha}<r/2$ and let $N\ge 2$ be such that $2^{ N-1} t^{1/\alpha} < r \le 2^{ N}t^{1/\alpha}$. For all $1\le n\le N$, we see that 
$2^{\alpha n}t<(2r)^\alpha \le (2\eta_1)^\alpha$ and  $|x(2^{\alpha n}t) -x|  = 2^{n-1}t^{1/\alpha} <r\le |x-y|/4$. Hence,  by \eqref{e:lifting-property},   for all $1\le n\le N$ and $z\in B(x(2^{\alpha n}t), 2^{n-3}\eta_2 t^{1/\alpha})$,
\begin{align}\label{e:lower-two-jumps-1}
	\begin{split} 
	& \delta_D(x(2^{\alpha n}t)) \ge 2^{n-1} \eta_2 t^{1/\alpha }, \quad \frac34\le \frac{\delta_D(z)}{\delta_D(x(2^{\alpha n}t))} \le \frac54,\\
	& \frac14\le \frac{|x-z|}{2^nt^{1/\alpha}} \le 1 \quad \text{and} \quad \frac12\le \frac{|y-z|}{|x-y|} \le  \frac32.
	\end{split}
\end{align}
Moreover, $B(x(2^{\alpha n}t), 2^{n-3}\eta_2 t^{1/\alpha})$, $1\le n\le N$, are disjoint.  By
Lemma \ref{l:HKE-lower-0}, \eqref{e:lower-two-jumps-1} and \eqref{e:Phi1-scaling}--\eqref{e:ell-scaling}, we have for all $1\le n \le N$ and  $z\in B(x(2^{\alpha n}t), 2^{n-3}\eta_2 t^{1/\alpha})$,
\begin{align}\label{e:lower-two-jumps-2}
	p(t/2,x,z) \ge \frac{c_4}{2^{(d+\alpha)n} t^{d/\alpha}}  \left(1 \wedge \frac{\delta_D(x)}{t^{1/\alpha}}\right)^q  A_{\Phi_1,\Phi_2,\ell}(t,x,x(2^{\alpha n}t)) 
\end{align}
and 
\begin{align}\label{e:lower-two-jumps-3}
	p(t/2,z,y) \ge \frac{c_4t}{|x-y|^{d+\alpha}}  \left(1 \wedge \frac{\delta_D(y)}{t^{1/\alpha}}\right)^q  A_{\Phi_1,\Phi_2,\ell}(t,x(2^{\alpha n}t),y). 
\end{align}
Using the semigroup property, \eqref{e:lower-two-jumps-2} and  \eqref{e:lower-two-jumps-3}, we obtain
\begin{align}\label{e:lower-two-jumps-4}
&	p(t,x,y) \ge \sum_{n=1}^N \int_{B(x(2^{\alpha n}t),2^{n-3} \eta_2 t^{1/\alpha})}  p(t/2, x,z) p(t/2, z, y) dz\\
	&\ge \frac{c_5t}{|x-y|^{d+\alpha}} \left(1 \wedge \frac{\delta_D(x)}{t^{1/\alpha}}\right)^q \left(1 \wedge \frac{\delta_D(y)}{t^{1/\alpha}}\right)^q  \nn\\
	&\quad \times \sum_{n=1}^N  \int_{B(x(2^{\alpha n}t),2^{n-3} \eta_2 t^{1/\alpha})} \frac{A_{\Phi_1,\Phi_2,\ell}(t,x, x(2^{\alpha n}t))A_{\Phi_1,\Phi_2,\ell}(t,x(2^{\alpha n}t),y)}{2^{(d+\alpha)n} t^{d/\alpha}} dz \nn\\
	&= \frac{c_6t}{|x-y|^{d+\alpha}} \left(1 \wedge \frac{\delta_D(x)}{t^{1/\alpha}}\right)^q \left(1 \wedge \frac{\delta_D(y)}{t^{1/\alpha}}\right)^q \sum_{n=1}^N   \frac{A_{\Phi_1,\Phi_2,\ell}(t,x, x(2^{\alpha n}t))A_{\Phi_1,\Phi_2,\ell}(t,x(2^{\alpha n}t),y)}{2^{\alpha n} }.\nn
\end{align}
On the other hand, 
by using \eqref{e:lifting-property} and the fact that  $2^N t^{1/\alpha}<2r \le |x-y|/2$, we see that for all $1\le n\le N$ and $u \in [2^{n-1}t^{1/\alpha},2^{n}t^{1/\alpha}]$,
\begin{align*}
\eta_2 &\le 	\frac{\delta_D(x+u\bn_x)}{\delta_D(x(2^{\alpha n}t))} \le \frac{2}{\eta_2}, \quad 1\le \frac{|x+u\bn_x - x|}{|x(2^{\alpha n}t) -x|} = \frac{u}{2^{n-1}t^{1/\alpha}} \le 2,\\
&\frac13= \frac{(1/2)|x-y|}{(3/2)|x-y|}\le  \frac{|x+u\bn_x- y|}{|x(2^{\alpha n}t)-y|}\le  \frac{(3/2)|x-y|}{(1/2)|x-y|}=3.
\end{align*}
Thus, by \eqref{e:Phi1-scaling}--\eqref{e:ell-scaling},  we have for all $1\le n\le N$ and $u \in [2^{n-1}t^{1/\alpha},2^{n}t^{1/\alpha}]$,
\begin{align*}
	A_{\Phi_1,\Phi_2,\ell}(t,x,x(2^{\alpha n}t))A_{\Phi_1,\Phi_2,\ell}(t,x(2^{\alpha n}t),y)\asymp A_{\Phi_1,\Phi_2,\ell}(t,x,x +u\bn_x)A_{\Phi_1,\Phi_2,\ell}(t,x+u\bn_x,y).
\end{align*}
Hence, since $r\le 2^Nt^{1/\alpha}$, it holds that  
\begin{align}\label{e:lower-two-jumps-5}
& \sum_{n=1}^N   \frac{A_{\Phi_1,\Phi_2,\ell}(t,x,x(2^{\alpha n}t))A_{\Phi_1,\Phi_2,\ell}(t,x(2^{\alpha n}t),y)}{2^{\alpha n} }\\
	&\ge c_7t\sum_{n=1}^N \int_{2^{n-1}t^{1/\alpha}}^{2^{n}t^{1/\alpha}} A_{\Phi_1,\Phi_2,\ell}(t,x,x +u\bn_x)A_{\Phi_1,\Phi_2,\ell}(t,x+u\bn_x,y)\frac{du}{u^{\alpha+1}} \nn\\
		&\ge c_7t \int_{t^{1/\alpha}}^{r} A_{\Phi_1,\Phi_2,\ell}(t,x,x +u\bn_x)A_{\Phi_1,\Phi_2,\ell}(t,x+u\bn_x,y)\frac{du}{u^{\alpha+1}} .\nn
\end{align}
Combining \eqref{e:lower-two-jumps-4} with \eqref{e:lower-two-jumps-5} completes the proof. \qed

\begin{thm}\label{t:LHK}
For any $T>0$, there exists  $C=C(T)>0$ such that 
for all $t\in(0,T]$ and $x,y \in D$,
\begin{align*}		
	\begin{split} 			p(t,x,y)&\ge C \left(1 \wedge \frac{\delta_D(x)}{t^{1/\alpha}}\right)^q \left(1 \wedge \frac{\delta_D(y)}{t^{1/\alpha}}\right)^q \left( t^{-d/\alpha} \wedge \frac{t}{|x-y|^{d+\alpha}}\right)  \\			&  \times \bigg[ \big( t \wedge |x-y|^\alpha\big)  \int_{		t^{1/\alpha} \wedge (\eps_1|x-y|/2)}^{\eps_1|x-y|}A_{\Phi_1,\Phi_2,\ell}(t,x,x+u\mathbf{n}_x)\,A_{\Phi_1,\Phi_2,\ell}(t,x+u\mathbf{n}_x,y)\frac{du}{u^{\alpha+1}}\\			& \quad + \big( t \wedge |x-y|^\alpha\big) \int_{
			t^{1/\alpha}\wedge (\eps_1|x-y|/2)}^{\eps_1|x-y|}A_{\Phi_1,\Phi_2,\ell}(t,y,y+u\mathbf{n}_y)\,A_{\Phi_1,\Phi_2,\ell}(t,y+u\mathbf{n}_y,x)\frac{du}{u^{\alpha+1}} \bigg] .		\end{split} 	\end{align*}
\end{thm} 
\pf Let $\eps_1$ be  as defined in  \eqref{e:def-eps1}. If $(\delta_D(x) \wedge \delta_D(y))\vee t^{1/\alpha}\ge \eps_1|x-y|/2$, then  the  result follows from Lemma \ref{l:HKE-lower-1} and \eqref{e:upper-bound-A}.  
Assume $(\delta_D(x) \wedge \delta_D(y))\vee  t^{1/\alpha}<\eps_1|x-y|/2$. Then  $t^{1/\alpha}<\eta_1/2$.
Applying Lemma \ref{l:HKE-lower-two-jumps} together with  symmetry of $p(t,\cdot,\cdot)$ and the fact that $a\ge c_1$ and $a\ge c_2$ imply $a\ge (c_1+c_2)/2$, we obtain the desired result. 
\qed

\begin{corollary}\label{c:life-2}
For any $T>0$,	there exists $C>0$ such that 
	\begin{align*}
		\P_x (\zeta>t) \asymp   \left(1 \wedge \frac{\delta_D(x)}{t^{1/\alpha}}\right)^q  
		\quad \text{for all $t\in (0,T]$ and $x\in D$}.
	\end{align*}
\end{corollary}
\pf 
Let $z:=x(\eta_1^\alpha t/T)$.  By \eqref{e:lifting-property}, it holds that for all $y\in B(z, 2^{-2}\eta_1\eta_2(t/T)^{1/\alpha})$,
\begin{align*}
	|x-y| \le |x-z| + 2^{-2}\eta_1\eta_2(t/T)^{1/\alpha} <\eta_1(t/T)^{1/\alpha} \quad \text{and} \quad 
		\delta_D(y)\ge 2^{-2}\eta_1\eta_2 (t/T)^{1/\alpha}.\end{align*} Hence, by Theorem \ref{t:LHK}, Lemma \ref{l:two-jumps-lower-bound}(ii)  and \eqref{e:interior-lower-bound-A}, we obtain for all $y\in B(z, 2^{-2}\eta_1\eta_2(t/T)^{1/\alpha})$,
		\begin{align*}	
			p(t,x,y) \ge c_1 \bigg( 1 \wedge \frac{\delta_D(x)}{t^{1/\alpha}} \bigg)^q A_{\Phi_1,\Phi_2,\ell}(t,x,y) \,t^{-d/\alpha} \ge c_2 \bigg( 1 \wedge \frac{\delta_D(x)}{t^{1/\alpha}} \bigg)^q t^{-d/\alpha}.\end{align*}
It follows that
\begin{align*}
		\P_x (\zeta>t)& \ge \int_{B(z, 2^{-2}\eta_1\eta_2(t/T)^{1/\alpha})} p(t,x,y)dy\\
			&  \ge c_2  \left(1 \wedge \frac{\delta_D(x)}{t^{1/\alpha}}\right)^q t^{-d/\alpha} \int_{B(x, 2^{-2}\eta_1\eta_2(t/T)^{1/\alpha})}dy  =  c_3 \left(1 \wedge \frac{\delta_D(x)}{t^{1/\alpha}}\right)^q .\end{align*}Combining the above with Corollary \ref{c:life}, we arrive at the result.  \qed


\section{Upper bound}\label{s:ub}
In this section, we establish the
sharp heat kernel upper bounds using a  self-improving argument. 
A related approach was used in \cite{CKSV23}.  
However, since the present setting involves general  functions $\Phi_1$, $\Phi_2$ and $\ell$, rather than power functions as in \cite{CKSV23}, a more refined method is required.  We begin with the following lemma.
Recall that  the function $A_{f, g, h}(t,x,y)$ is defined in \eqref{e:def-A(f,g,h)}, and the class $\sM^\uparrow(\gamma,\gamma^*)$ is defined at the beginning of Section \ref{s:preliminary}.

\begin{lemma}\label{l:first-term}
Let $F\in \sM^\uparrow(\gamma,\gamma^*)$ for some $\gamma^*\ge \gamma \ge 0$. 	Suppose that there exist constants $T>0$ and $c_0>0$  such that 	for all $t\in (0, T]$ and $y,z \in D$,
 	\begin{equation}\label{e:first-term-ass}
 		p(t,y,z) \le c_0  \left(1 \wedge \frac{\delta_D(y)}{t^{1/\alpha}}\right)^{q}
 		A_{F, 1,1}(t,y,z)
 		\left( t^{-d/\alpha} \wedge \frac{t}{|y-z|^{d+\alpha}}\right).
 	\end{equation}
 	Then there exists $C >0$  such that for all $t\in (0,T]$, 
 	$r\in (0,\text{\rm diam}(D))$  and $x\in D$,
 	\begin{equation}\label{e:first-term}
 		\P_x(\tau_{B_{\overline D}(x, r)}<t< \zeta)\le \frac{Ct}{r^\alpha} \left(1 \wedge \frac{\delta_D(x)}{t^{1/\alpha}}\right)^{q}   F\bigg(\frac{\delta_D(x)\vee t^{1/\alpha}}{r}\bigg).
 	\end{equation}
 \end{lemma}
 \pf  We first note that  \eqref{e:first-term-ass} implies that for all $t\in (0,T]$ and $y \in D$,
 \begin{align}\label{e:first-term-pre}
 	&\P_y(|Y_t-y|>r/2, \, t<\zeta)= \int_{z \in D:|y-z| > r/2} p(t,y,z)dz	 \\
 & \le c_0 t\left(1 \wedge \frac{\delta_D(y)}{t^{1/\alpha}}\right)^{q}	\int_{z \in D: |y-z| > r/2} F\bigg(\frac{(\delta_D(y) \wedge \delta_D(z))\vee t^{1/\alpha}}{|y-z|}\bigg)\frac{dz}{|y-z|^{d+\alpha}}\nn\\
 & \le c_1 t\left(1 \wedge \frac{\delta_D(y)}{t^{1/\alpha}}\right)^{q} F\bigg(\frac{\delta_D(y)\vee t^{1/\alpha}}{r}\bigg)	\int_{z \in D: |y-z| > r/2} \frac{dz}{|y-z|^{d+\alpha}} \nn\\
 	&\le  \frac{c_2t}{r^{\alpha}}\left(1 \wedge \frac{\delta_D(y)}{t^{1/\alpha}}\right)^{q} 	 F\bigg(\frac{\delta_D(y)\vee t^{1/\alpha}}{r}\bigg).\nn
 \end{align}
 We used the almost increasing and upper  scaling properties  of $F$ in the second inequality above. On the other hand, by Proposition \ref{p:upper-heatkernel}, we have
 \begin{equation}\label{e:regular-path}
 	\sup_{s \le t,\, y \in D}	\P_y\big( |Y_s - y|\ge r/2, \, s<\zeta\big) \le c_3 \sup_{s \le t, \, y \in D} \int_{z \in D: |y-z| \ge r/2} \frac{s}{|y-z|^{d+\alpha}}dz \le \frac{c_4t}{r^\alpha}.
 \end{equation}
 Set $c_5:=1/(2c_4)$. If $tr^{-\alpha}\ge c_5$, then
  by the almost increasing property of $F$, 
  $$t r^{-\alpha}F((\delta_D(x)\vee t^{1/\alpha})/r) \ge c_5 c_6F(c_5^{1/\alpha}),$$
hence \eqref{e:first-term} follows from Corollary \ref{c:life}. 
 
 Assume $t r^{-\alpha} <c_5$. 
 Set $B:=B_{\overline D}(x, r)$. By the strong Markov property and \eqref{e:regular-path},  
 \begin{align*}
 &\P_x\big(\tau_{B}<t< \zeta, \, |Y_t-Y_{\tau_{B}}|\ge r/2\big)	=  \E_x\left[ \P_{Y_{\tau_{B}}}  	\left( |Y_{t- \tau_{B} } - Y_0 |  \ge r/2\right) :  	\tau_{B}<t< \zeta \right] \nn\\	& \le \P_x(\tau_{B}<t< \zeta)  \sup_{s \le t,\, y \in D}	\P_y( |Y_s - y|\ge r/2, \, s<\zeta) 	\le 2^{-1}\P_x(\tau_{B}<t< \zeta).
 \end{align*}
 Thus,  we get
 \begin{align}\label{e:first-term-1}
 	\P_x(\tau_{B}<t< \zeta)&=2\big(\P_x(\tau_{B}<t< \zeta)-2^{-1}\P_x(\tau_{B}<t< \zeta)\big)
 	\\
 	&\le 2 \big(\P_x(\tau_{B}<t< \zeta)-\P_x(\tau_{B}<t< \zeta, \, |Y_t-Y_{\tau_{B}}|\ge r/2) \big)\nn\\
 	&=  2\P_x(\tau_{B}<t< \zeta, \, |Y_t-Y_{\tau_{B}}|<r/2).
 	\nn
 \end{align}
By the triangle inequality, for any $y \in D\setminus B$ and 
  $z\in B(y,r/2)$, we have  $|z-x|>r/2$. 
 Therefore, by \eqref{e:first-term-1} and \eqref{e:first-term-pre},  we obtain that
 \begin{align*}
 	&\P_x(\tau_{B}<t< \zeta)\le  2\P_x(\tau_{B}<t< \zeta, \, |Y_t-Y_{\tau_{B}}|<r/2)\\
 	& \le 2\P_x(|Y_t-x|>r/2, \, t<\zeta)
 	\le \frac{2c_2t}{r^\alpha} \left(1 \wedge \frac{\delta_D(x)}{t^{1/\alpha}}\right)^{q}   
 	F\bigg(\frac{\delta_D(x)\vee t^{1/\alpha}}{r}\bigg).
 \end{align*} 
 The proof is complete.\qed
 
 For any $\beta \ge 0$, define
 \begin{align*}
 	F_\beta(r):=(r^\beta \wedge 1), \quad\;\; r>0.
 \end{align*}
 
Recall that $\beta_1$ and $\beta_1^*$ denote the lower and upper Matuszewska indices of $\Phi_1$ respectively.

 \begin{lemma}\label{l:UHK-case1-induction-pre}  
 	Let $\underline \beta_1 \ge 0$ be such that $\underline \beta_1 \in ((q-\alpha)_+,\beta_1)$ if $\beta_1>0$, 
 	and $\underline \beta_1 =0$ if $\beta_1=0$.
 	For any $T>0$, there exists a constant $C=C(\underline \beta_1, T)>0$ such that for all $t\in (0,T]$ and $x,y \in D$,
 	\begin{equation*}
 		p(t,x,y) \le C \left(1 \wedge \frac{\delta_D(x) \wedge \delta_D(y)}{t^{1/\alpha}}\right)^{q}
 		A_{F_{\underline \beta_1}, 1, 1}(t,x,y)\left( t^{-d/\alpha} \wedge \frac{t}{|x-y|^{d+\alpha}}\right).
 	\end{equation*}
 \end{lemma}
 \pf   Let  $T>0$ be fixed.   
 Let $N:= \lfloor \underline \beta_1/\alpha \rfloor+ 1$ and set $\theta:=1/N$. 
 For any  integer $0\le n \le N$, we define
 $$
 f_{n}(r):=F_{n\theta\underline \beta_1}(r), \quad r>0.
 $$ 
 Note that, 
 by \eqref{e:Phi0-scaling}, 
  there exists $c_1\ge 1$ such that
 \begin{align}\label{e:f-n-prop0}
 	\Phi_{0}(r) \le c_1 F_{\underline \beta_1}(r)   \le c_1 f_n(r) \quad \text{for all $r>0$ and $0\le n \le N$.}
 \end{align}
 
We will prove below by induction that for any $0\le n\le N$, there exists  $C=C(\underline \beta_1, T,n)>0$ such that  for all $t\in (0,T]$ and $x,y \in D$,
 \begin{align}\label{e:UHK-f-induction-0}
 	p(t,x,y) \le C \left(1 \wedge \frac{\delta_D(x) \wedge \delta_D(y)}{t^{1/\alpha}}\right)^{q}
 	A_{f_{n},1, 1}(t,x,y)\left( t^{-d/\alpha} \wedge \frac{t}{|x-y|^{d+\alpha}}\right).
 \end{align} 
 Then the assertion of the lemma follows immediately from \eqref{e:UHK-f-induction-0} by taking taking $n=N$.
 
 When $n=0$, we have  $A_{f_{0}, 1,1}(t,x,y)= 1$. Hence, \eqref{e:UHK-f-induction-0} holds by Proposition \ref{p:UHK-rough}.
 
 Suppose \eqref{e:UHK-f-induction-0} holds for $n-1$,  
 $1\le n \le N$. By symmetry, we can assume without loss of generality that $\delta_D(x) \le \delta_D(y)$. 
 Set $a:=\wh{R}/(48\,\mathrm{diam}(D))$ and $r:=a|x-y| \le (|x-y| \wedge \wh R)/48$.

 If $\delta_D(x) \vee t^{1/\alpha} \ge 
 r/4$, then 
 \begin{align*}
 		A_{f_{n},1, 1}(t,x,y)=\bigg(\frac{\delta_D(x)\vee t^{1/\alpha}}{|x-y|} \wedge 1\bigg)^{n \theta \underline \beta_1} \ge (a/4)^{n \theta \underline \beta_1}.
 \end{align*}
 Thus,  \eqref{e:UHK-f-induction-0} 
 follows from Proposition \ref{p:UHK-rough}.

 Assume $ \delta_D(x)\vee t^{1/\alpha} < 
 r/4$.  
Our goal is to show that there exists $c_2>0$ independent of $t,x,y$ such that 
 \begin{equation}\label{e:UHK-f-case1}
 	p(t,x,y) \le \frac{c_2t}{r^{d+\alpha}}\left(1 \wedge \frac{\delta_D(x) }{t^{1/\alpha}}\right)^{q} f_{n}\bigg(\frac{\delta_D(x)\vee t^{1/\alpha}}{r}\bigg).
 \end{equation}
Set $V_1:=B_{\overline D}(x,r)$,  $V_3:=B_{\overline D}(y,r)$ and $V_2:=\overline{D} \setminus (V_1 \cup V_3)$.  By Proposition \ref{p:upper-heatkernel} and the triangle inequality,
 \begin{equation}\label{e:UHK-case1-2}
 	\sup_{s \le t, \, z \in V_2} p(s,z,y) \le  \sup_{s \le t, \, z \in D,\, |z-y|\ge 
 	r} \frac{c_3s}{|z-y|^{d+\alpha}} = \frac{c_3t}{r^{d+\alpha}}
 \end{equation}
 and
 \begin{equation}\label{e:UHK-case1-3}
 	{\rm dist}(V_1, V_3) \ge \sup_{u\in V_1, \, w \in V_3} (|x-y|-|x-u|-|y-w|)> 46r.
 \end{equation}  
 By  the induction hypothesis,  condition \eqref{e:first-term-ass} in Lemma  \ref{l:first-term} holds for   $F=f_{n-1}$. 
 Thus,  combining
 Lemma  \ref{l:first-term} and \eqref{e:UHK-case1-2} with the facts that $\theta \underline \beta_1<\alpha$ and $\delta_D(x)\vee t^{1/\alpha}<r/4$, we obtain 
 \begin{align}\label{e:UHK-f-case1-term1}
 	&\P_x(\tau_{V_1}<t< \zeta) \sup_{s \le t, \, z \in V_2} p(s,z,y) \\
 	&\le 
 	\frac{c_4t^2}{r^{d+2\alpha}}\left(1 \wedge \frac{\delta_D(x)}{t^{1/\alpha}}\right)^{q}
 	f_{n-1}\bigg(\frac{\delta_D(x)\vee t^{1/\alpha}}{r}\bigg) \nn\\
 &=  \frac{c_4t}{r^{d+\alpha}}  \left(1 \wedge \frac{\delta_D(x)}{t^{1/\alpha}}\right)^q  	f_{n-1}\bigg(\frac{\delta_D(x)\vee t^{1/\alpha}}{r}\bigg)  \bigg(\frac{t^{1/\alpha}}{r}\bigg)^{\theta \underline \beta_1} \bigg(\frac{t^{1/\alpha}}{r}\bigg)^{\alpha-\theta \underline \beta_1} \,\nn\\
 	&\le  \frac{c_4t}{r^{d+\alpha}}   \left(1 \wedge \frac{\delta_D(x)}{t^{1/\alpha}}\right)^q  	f_{n-1}\bigg(\frac{\delta_D(x)\vee t^{1/\alpha}}{r}\bigg)  \bigg(\frac{\delta_D(x)\vee t^{1/\alpha}}{r}\bigg)^{\theta \underline \beta_1} \nn\\
 	&\le  \frac{c_4t}{r^{d+\alpha}} \left(1 \wedge \frac{\delta_D(x)}{t^{1/\alpha}}\right)^q  
 	f_{n}\bigg(\frac{\delta_D(x)\vee t^{1/\alpha}}{r}\bigg)  .\,\nn
 \end{align}

 By   Lemma \ref{l:general-upper-2} and  \eqref{e:UHK-f-case1-term1}, to get the inequality \eqref{e:UHK-f-case1}, it suffices to show that
 \begin{align}\label{e:UHK-case1-f-term2-claim}
 	&\int_0^t\int_{V_3} \int_{V_1} p^{V_1}(s, x, u) \sB(u,w) p(t-s, y,w) du dw ds
 	\le c_5 t\left(1\wedge \frac{\delta_D(x)}{t^{1/\alpha}}\right)^q f_{n}\bigg(\frac{\delta_D(x)\vee t^{1/\alpha}}{r}\bigg). 
 \end{align}
 By \eqref{e:B4-a}, \eqref{e:UHK-case1-3} and  the almost   increasing property of $\Phi_0$, we have
 for all $u \in V_1$ and $w \in V_3$,
 \begin{align}\label{e:case1-B}
 	\sB(u,w) \le c_6 \Phi_0\left(\frac{\delta_D(u)\wedge \delta_D(w)}{|u-w|}\right)
 	\le c_{7}\Phi_0 \left(\frac{\delta_D(u)}{r}\right).
 \end{align}
For \eqref{e:UHK-case1-f-term2-claim},  we consider the following two cases separately.
 
 \smallskip

 \noindent
 \textbf{Case 1:}  $q-\alpha \ge n \theta \underline \beta_1$  and $\delta_D(x)< t^{1/\alpha}$.  
 Let $Q\in \partial D$ satisfy $|x-Q|=\delta_D(x)$. Then we have  $B_D(Q,r/2 ) \subset V_1 \subset B_D(Q,2r)$ and $x\in B_D(Q, r/4)$.  
 By  \eqref{e:case1-B},  \eqref{e:f-n-prop0} and Proposition \ref{p:bound-for-integral-new}, since  $\underline \beta_1>q-\alpha$, we get
 \begin{align*}
 	&\int_0^t\int_{V_3} \int_{V_1} p^{V_1}(s, x, u) \sB(u,w) p(t-s, y,w) du dw ds\\
	&\le c_{8}\int_0^t \int_{V_1} p^{V_1}(s, x, u) \left(\frac{\delta_D(u)}{r}\right)^{\underline\beta_1} du   \int_{V_3} p(t-s, y,w)  dw ds\nn \\
 	&\le c_{8}r^{-\underline \beta_1}\int_0^\infty \int_{V_1} p^{V_1}(s, x, u) \delta_D(u)^{\underline\beta_1} du ds \le c_{9}  r^{\alpha-q} \delta_D(x)^q.
 \end{align*}
 Thus, using
  $\delta_D(x)<t^{1/\alpha}<r/8$
  and  $q-\alpha \ge n \theta \underline \beta_1$, we obtain
 \begin{align*}
 	&\int_0^t\int_{V_3} \int_{V_1} p^{V_1}(s, x, u) \sB(u,w) p(t-s, y,w) du dw ds\\
 	&\le 	 c_{9} t \left(1\wedge \frac{\delta_D(x)}{t^{1/\alpha}}\right)^q \bigg( \frac{\delta_D(x) \vee  t^{1/\alpha}}{r}\bigg)^{q-\alpha}	\le c_{9} t\left(1\wedge \frac{\delta_D(x)}{t^{1/\alpha}}\right)^q f_{n}\bigg(\frac{\delta_D(x)\vee t^{1/\alpha}}{r}\bigg), \nn
 \end{align*}
 proving that \eqref{e:UHK-f-case1} holds in this case.
 
 \smallskip
 
 \noindent
 \textbf{Case 2:} $q- \alpha < n\theta \underline \beta_1$ or $\delta_D(x)\ge t^{1/\alpha}$.  By \eqref{e:case1-B} and \eqref{e:f-n-prop0}, we have 
 \begin{align}\label{e:case2-B}
 	&\int_0^t\int_{V_3} \int_{V_1} p^{V_1}(s, x, u) \sB(u,w) p(t-s, y,w) du dw ds \\
 	&\le c_{10}\int_0^t \int_{V_3} p(t-s, y,w)   \int_{V_1} p(s, x, u) f_n \left(\frac{\delta_D(u)}{r}\right) du   dw ds.\nn
 \end{align}
Using  Corollary \ref{c:life},  we get that for any $0<s < t$, 
 \begin{align}\label{e:f-case2-1}
 	& \int_{u \in V_1:\delta_D(u) < \delta_D(x)} p(s, x, u) f_n\left(\frac{\delta_D(u)}{r}\right) du \le  f_n\bigg(\frac{\delta_D(x)}{r}\bigg) \int_{u \in V_1:\delta_D(u) < \delta_D(x)} p(s, x, u) du \\
 	&\le \P_x(\zeta>s)  f_n\bigg(\frac{\delta_D(x) \vee s^{1/\alpha}}{r}\bigg)
 		\le  c_{11}   \left(1 \wedge \frac{\delta_D(x)}{s^{1/\alpha}} \right)^qf_{n}\bigg(\frac{\delta_D(x) \vee s^{1/\alpha}}{r}\bigg) \nn.
 \end{align}

 Next, we note that the lower Matuszewska index of $f_{n-1}$ is $(n-1)\theta \underline \beta_1$ and the upper Matuszewska index of $f_n$ is $n\theta \underline \beta_1$.  Since $\theta \underline \beta_1<\alpha$,  by applying 
 Lemma \ref{l:analog-of-10.10}(i)  with  $\Phi= f_{n-1}$, $\Psi=f_n$ and $\ell=1$,  we have that for any $0<s<t$,
  \begin{align}\label{e:f-case2-2}	&\int_{u \in V_1:\delta_D(u) \ge \delta_D(x)}  	 \left(s^{-d/\alpha} \wedge \frac{s}{|x-u|^{d+\alpha}}\right)f_{n-1}\bigg(\frac{\delta_D(x)\vee s^{1/\alpha}}{|x-u|} \bigg) f_{n}\left(\frac{\delta_D(u)}{r}\right) du  \\
  &\le c_{12} f_{n} \bigg(\frac{\delta_D(x) \vee s^{1/\alpha}}{r}\bigg) . \nn\end{align}
 Using the induction hypothesis, \eqref{e:f-n-prop0} and \eqref{e:f-case2-2}, we get  that for any $0<s< t$,
 \begin{align}\label{e:f-case2-3}
 	& \int_{u \in V_1:\delta_D(u) \ge \delta_D(x)} p(s, x, u) f_n \left(\frac{\delta_D(u)}{r}\right) du  \\
 	&\le c_{13} \left(1 \wedge \frac{\delta_D(x)}{s^{1/\alpha}} \right)^q  \nn\\
 	&\quad \times \int_{u \in V_1:\delta_D(u) \ge \delta_D(x)}  
  \left(s^{-d/\alpha} \wedge \frac{s}{|x-u|^{d+\alpha}}\right)f_{n-1}\bigg(\frac{\delta_D(x)\vee s^{1/\alpha}}{|x-u|} \bigg)  f_n\left(\frac{\delta_D(u)}{r}\right)du \nn\\
 	&\le c_{14}  \left(1 \wedge \frac{\delta_D(x)}{s^{1/\alpha}} \right)^qf_{n} \bigg(\frac{\delta_D(x) \vee s^{1/\alpha}}{r}\bigg). \nn
 \end{align}
 Combining \eqref{e:case2-B}, \eqref{e:f-case2-1} and  \eqref{e:f-case2-3}, and using  Lemma \ref{cal:00}, we conclude that
 \begin{align*}
 	&\int_0^t\int_{V_3} \int_{V_1} p^{V_1}(s, x, u) \sB(u,w) p(t-s, y,w) du dw ds \\
 	&\le c_{15}   \int_0^t \left(1 \wedge \frac{\delta_D(x)}{s^{1/\alpha}} \right)^q  f_{n}\bigg(\frac{\delta_D(x)\vee s^{1/\alpha}}{r}\bigg) \int_{V_3} p(t-s, y,w)    dw ds \nn\\
 	&\le c_{15}   \int_0^t 
 	\left(1 \wedge \frac{\delta_D(x)}{s^{1/\alpha}} \right)^q f_{n}\bigg(\frac{\delta_D(x)\vee s^{1/\alpha}}{r}\bigg)  ds \le c_{17}	t\left(1 \wedge \frac{\delta_D(x)}{t^{1/\alpha}} \right)^q f_{n}\bigg(\frac{\delta_D(x)\vee t^{1/\alpha}}{r}\bigg) .
 \end{align*}
 This finishes the proof of \eqref{e:UHK-case1-f-term2-claim} and hence \eqref{e:UHK-f-case1} holds for $n$. The proof is complete.
 \qed
 
 \begin{lemma}\label{l:UHK-case1-induction}
	If $\beta_1^*<\alpha +  \beta_1$, then for 
	any $T>0$, there exists a constant $C=C(T)>0$ such that for all $t\in (0,T]$ and $x,y \in D$,
 	\begin{equation*}
 		p(t,x,y) \le C \left(1 \wedge \frac{\delta_D(x) \wedge \delta_D(y)}{t^{1/\alpha}}\right)^q 
 		A_{\Phi_0, 1, 1}(t,x,y)\left( t^{-d/\alpha} \wedge \frac{t}{|x-y|^{d+\alpha}}\right).
 	\end{equation*}
 \end{lemma}
 \pf   Let $T>0$ be fixed. 
Choose $\eps>0$ such that $q \vee \beta_1^*<\alpha+\beta_1-2\eps$, and let $\lb_1:=(\beta_1-\eps)_+$ and $\ub_1:=\beta_1^*+\eps$. Note that $q \vee \ub_1<\alpha+\lb_1$. Set $N:=\lfloor \ub_1/\alpha \rfloor + 1$ and $\theta:=1/N$.
For each integer $0\le n\le N$, we define
 $$
 \Phi_{0,n}(r):= \Phi_0(r)^{n\theta}F_{\underline \beta_1}(r)^{1-n\theta}, \quad r>0.
 $$ 
 It is easy to see that $\Phi_{0,n}$ satisfies the upper and lower scaling conditions with indices $n\theta\overline{\beta}_1 + (1-n\theta)\underline \beta_1$ and $\underline{\beta}_1$ respectively.  
 By \eqref{e:Phi0-scaling}, there exists $c_1\ge 1$ such that
 \begin{align*}
 	c_1^{-1} F_{\overline \beta_1}(r)\le 	\Phi_{0}(r) \le c_1  F_{\underline \beta_1}(r) \quad \text{for all $r>0$.}
 \end{align*}
 This implies that for all $0\le n \le N$ and $r\in (0,1]$,
 \begin{gather}
 	\Phi_{0,n}(r)  \ge c_1^{-1}\Phi_{0}(r),\label{e:Phi-n-prop0}\\
 		\Phi_{0,n+1}(r) 
 	 \ge c_1^{-1} \Phi_{0,n}(r) r^{\theta (\overline \beta_1 -\underline \beta_1)  }  \ge c_1^{-1} \Phi_{0,n}(r) r^{\theta \overline \beta_1  } ,\label{e:Phi-n-prop1}\\
 		\Phi_{0,n}(r) \ge  c_1^{-1} r^{n\theta \overline \beta_1 + (1-n\theta) \underline \beta_1}.\label{e:Phi-n-prop2}
 \end{gather}

 We will prove below by induction that there exists a constant $C=C(T)>0$ such that for all $0\le n\le N$, $t\in (0,T]$ and $x,y \in D$,
 \begin{align}\label{e:UHK-induction-0}
 	p(t,x,y) \le C \left(1 \wedge \frac{\delta_D(x) \wedge \delta_D(y)}{t^{1/\alpha}}\right)^{q}
	A_{\Phi_{0,n}, 1, 1}
	(t,x,y)\left( t^{-d/\alpha} \wedge \frac{t}{|x-y|^{d+\alpha}}\right).
 \end{align}  
 Then the assertion of the lemma follows from \eqref{e:UHK-induction-0} directly upon taking $n=N$.
 
 When $n=0$,  \eqref{e:UHK-induction-0} holds by Lemma \ref{l:UHK-case1-induction-pre}.
 
 Suppose \eqref{e:UHK-induction-0} holds for $n-1$, $1\le n \le N$. By symmetry, we can assume without loss of generality that $\delta_D(x) \le \delta_D(y)$. Set $a:=\wh{R}/(48\,\mathrm{diam}(D))$ 
 and $r:=a|x-y|$.
As in the proof of Lemma \ref{l:UHK-case1-induction-pre},  by Proposition \ref{p:UHK-rough}, it suffices to consider the case $\delta_D(x) \vee t^{1/\alpha} < r/4$ and show that the following holds:
 \begin{equation}\label{e:UHK-case1}
 	p(t,x,y) \le \frac{c_2t}{r^{d+\alpha}}\left(1 \wedge \frac{\delta_D(x) }{t^{1/\alpha}}\right)^{q} \Phi_{0,n}\bigg(\frac{\delta_D(x)\vee t^{1/\alpha}}{r}\bigg),
 \end{equation}
 where $c_2>0$ is a constnat independent of $t,x,y$.
 
 Assume $\delta_D(x)\vee t^{1/\alpha}<r/4$. Set $V_1:=B_{\overline D}(x,r)$, 
 $V_3:=B_{\overline D}(y,r)$ and $V_2:=\overline{D} \setminus (V_1 \cup V_3)$. Note that \eqref{e:UHK-case1-2} and \eqref{e:UHK-case1-3} are valid. By  the induction hypothesis,  condition \eqref{e:first-term-ass} in Lemma  \ref{l:first-term} holds for $F=\Phi_{0,n-1}$. 
 Thus, using Lemma  \ref{l:first-term}, \eqref{e:UHK-case1-2} and \eqref{e:Phi-n-prop1} with the facts that $\theta \overline \beta_1<\alpha$ 
 and $\delta_D(x) \vee t^{1/\alpha}<r/4$, we obtain
 \begin{align}\label{e:UHK-case1-term1}
 	&\P_x(\tau_{V_1}<t< \zeta) \sup_{s \le t, \, z \in V_2} p(s,z,y) \\
 	&\le 
 	\frac{c_3t^2}{r^{d+2\alpha}}\left(1 \wedge \frac{\delta_D(x)}{t^{1/\alpha}}\right)^{q}
 	\Phi_{0,n-1}\bigg(\frac{\delta_D(x)\vee t^{1/\alpha}}{r}\bigg) \nn\\
 	&=  \frac{c_3t}{r^{d+\alpha}}  \left(1 \wedge \frac{\delta_D(x)}{t^{1/\alpha}}\right)^q  
 	\Phi_{0,n-1}\bigg(\frac{\delta_D(x)\vee t^{1/\alpha}}{r}\bigg)  \bigg(\frac{t^{1/\alpha}}{r}\bigg)^{\theta \overline \beta_1} \bigg(\frac{t^{1/\alpha}}{r}\bigg)^{\alpha-\theta \overline \beta_1} \,\nn\\
 	&\le  \frac{c_3t}{r^{d+\alpha}}   \left(1 \wedge \frac{\delta_D(x)}{t^{1/\alpha}}\right)^q  
 	\Phi_{0,n-1}\bigg(\frac{\delta_D(x)\vee t^{1/\alpha}}{r}\bigg)  \bigg(\frac{\delta_D(x)\vee t^{1/\alpha}}{r}\bigg)^{\theta \overline \beta_1} \nn\\
 	&\le  \frac{c_4t}{r^{d+\alpha}} \left(1 \wedge \frac{\delta_D(x)}{t^{1/\alpha}}\right)^q  
 	\Phi_{0,n}\bigg(\frac{\delta_D(x)\vee t^{1/\alpha}}{r}\bigg)  .\,\nn
 \end{align}

 Observe that \eqref{e:case1-B} is still valid.
 Using  Corollary \ref{c:life}, the almost increasing property of  $\Phi_0$  and \eqref{e:Phi-n-prop0}, we get that for any $0<s < t$, 
 \begin{align}\label{e:case2-1}
 	& \int_{u \in V_1:\delta_D(u) < \delta_D(x)} p(s, x, u) \Phi_0 \bigg(\frac{\delta_D(u)}{r}\bigg) du \\
 	&\le c_{5}\P_x(\zeta>s)\Phi_0 \bigg(\frac{\delta_D(x) \vee s^{1/\alpha}}{r}\bigg) 
 	\le  c_{6}   \left(1 \wedge \frac{\delta_D(x)}{s^{1/\alpha}} \right)^q\Phi_{0,n}\bigg(\frac{\delta_D(x) \vee s^{1/\alpha}}{r}\bigg) \nn.
 \end{align}
 Next,  recall that $\Phi_{0,n-1}$ satisfies the lower  scaling condition with index 
 $ \underline \beta_1$ and $\Phi_{0,n}$ satisfies the upper  scaling condition with index 
 $n\theta \overline{\beta}_1 + (1-n\theta)\underline \beta_1$. Since
 \begin{align*}
 	&\alpha+ \underline \beta_1 - (n\theta \overline{\beta}_1 + (1-n\theta)\underline \beta_1)  = \alpha -   n \theta  (\overline \beta_1- \underline \beta_1)\ge \alpha - (\overline \beta_1- \underline \beta_1)   >0,
 \end{align*}
 we get from Lemma \ref{l:analog-of-10.10}(i) (with $\Phi=\Phi_{0,n-1}$, $\Psi=\Phi_{0,n}$ and $\ell=1$)  that
 \begin{align}\label{e:case2-2}
 	&\int_{u \in V_1:\delta_D(u) \ge \delta_D(x)}   \left(s^{-d/\alpha} \wedge \frac{s}{|x-u|^{d+\alpha}}\right) 
 	\Phi_{0,n-1}\bigg(\frac{\delta_D(x)\vee s^{1/\alpha}}{|x-u|} \bigg) \Phi_{0,n}\bigg(\frac{\delta_D(u)}{r}\bigg) du  \\
 	&\le c_{7} \Phi_{0,n} \bigg(\frac{\delta_D(x) \vee s^{1/\alpha}}{r}\bigg).\nn
 \end{align}
 Using the induction hypothesis, \eqref{e:Phi-n-prop0} and \eqref{e:case2-2}, we get  that for any $0<s< t$,
 \begin{align}\label{e:case2-3}
 	& \int_{u \in V_1:\delta_D(u) \ge \delta_D(x)} p(s, x, u) \Phi_0\left(\frac{\delta_D(u)}{r}\right) du  \\
 	&\le c_{8} \left(1 \wedge \frac{\delta_D(x)}{s^{1/\alpha}} \right)^q \nn\\
 	&\quad \times  \int_{u \in V_1:\delta_D(u) \ge \delta_D(x)}  
 	 \left(s^{-d/\alpha} \wedge \frac{s}{|x-u|^{d+\alpha}}\right) \Phi_{0,n-1}\bigg(\frac{\delta_D(x)\vee s^{1/\alpha}}{|x-u|} \bigg)  \Phi_{0,n}\bigg(\frac{\delta_D(u)}{r}\bigg) du\nn\\
 	&\le c_{9}  \left(1 \wedge \frac{\delta_D(x)}{s^{1/\alpha}} \right)^q\Phi_{0,n} \bigg(\frac{\delta_D(x) \vee s^{1/\alpha}}{r}\bigg). \nn
 \end{align}
 Combining \eqref{e:case1-B}, \eqref{e:case2-1} and  \eqref{e:case2-3}, 
 using Lemma \ref{cal:00} (with the fact that $\Phi_{0,n}$ satisfies the lower scaling condition with index $\underline \beta_1>q-\alpha$), we arrive at 
 \begin{align}\label{e:UHK-case1-term2}
 	&\int_0^t\int_{V_3} \int_{V_1} p^{V_1}(s, x, u) \sB(u,w) p(t-s, y,w) du dw ds \\
 	&\le c_{10}   \int_0^t \left(1 \wedge \frac{\delta_D(x)}{s^{1/\alpha}} \right)^q  \Phi_{0,n}\bigg(\frac{\delta_D(x)\vee s^{1/\alpha}}{r}\bigg) \int_{V_3} p(t-s, y,w)    dw ds \nn\\
 	&\le c_{11}   \int_0^t 
 	\left(1 \wedge \frac{\delta_D(x)}{s^{1/\alpha}} \right)^q \Phi_{0,n}\bigg(\frac{\delta_D(x)\vee s^{1/\alpha}}{r}\bigg)  ds\nn\\
 	& \le c_{12} t
 	\left(1 \wedge \frac{\delta_D(x)}{t^{1/\alpha}} \right)^q \Phi_{0,n}\bigg(\frac{\delta_D(x)\vee t^{1/\alpha}}{r}\bigg). \nn	
 \end{align}
 Applying Lemma \ref{l:general-upper-2}, from \eqref{e:UHK-case1-term1} and \eqref{e:UHK-case1-term2}, we deduce that  \eqref{e:UHK-case1} holds for $n$. The proof is complete.\qed

 \begin{lemma}\label{l:firstpart} 
Let $r\in (0, \text{\rm diam}(D))$, $t>0$ and  $x,y\in D$ be such that $\delta_D(x)\le \delta_D(y)$ and  $\delta_D(x) \vee t^{1/\alpha} \vee(|x-y|/4)  \le r$. If $\beta_1^*<\alpha+\beta_1$, then
there exists $C>0$ independent of $r,t,x$ and $y$ such that 
 	\begin{align*}	&\P_x(\tau_{B_{\overline D}(x,r)}<t< \zeta) \sup_{s \le t, \, z \in D,  |z-y|\ge r} p(s,z,y) \nn \\	&\le \frac{Ct^2}{r^{d+2\alpha}} \bigg(1 \wedge \frac{\delta_D(x)}{t^{1/\alpha}}\bigg)^q \bigg( 1 \wedge \frac{\delta_D(y)}{t^{1/\alpha}}\bigg)^q  A_{\Phi_0,\Phi_0,1}(t,x,y).
	\end{align*}
 \end{lemma}
 \pf  By Lemma \ref{l:UHK-case1-induction}, \eqref{e:first-term-ass} holds with $F=\Phi_0$. Thus, by  Lemma \ref{l:first-term} and  \eqref{e:Phi0-scaling} (with $|x-y|\le 4r$), to get the result, it suffices to show that  
 \begin{align}\label{e:firstpart-claim}
 	\sup_{s \le t, \, z \in D,  |z-y|\ge r} p(s,z,y)& \le  \frac{c_1t}{r^{d+\alpha}} \left( 1 \wedge \frac{\delta_D(y)}{t^{1/\alpha}}\right)^q \Phi_0\bigg(\frac{\delta_D(y)\vee t^{1/\alpha}}{r}\bigg) .
 \end{align}
 By  Lemma \ref{l:UHK-case1-induction} and the almost increasing property of $\Phi_0$, we have 
 \begin{align}\label{e:firstpart-1}
 		\sup_{s \le t, \, z \in D,  |z-y|\ge r} p(s,z,y) &\le c_2 \sup_{s \le t, \, z \in D, |z-y|\ge  r} s \left( 1 \wedge \frac{\delta_D(z) \wedge \delta_D(y)}{s^{1/\alpha}}\right)^q  \frac{ 
 		A_{\Phi_0,1,1}(s,z,y)}{|z-y|^{d+\alpha}} \\
 	&\le  c_3 \sup_{s\le t}	\frac{ s}{r^{d+\alpha}} \left( 1 \wedge \frac{\delta_D(y)}{s^{1/\alpha}}\right)^q 
 	\Phi_0\bigg(\frac{\delta_D(y)\vee s^{1/\alpha}}{r}\bigg) .\nn
 \end{align}
 Let $\eps:= \alpha+\beta_1-q>0$. Note that $\delta_D(y)\le \delta_D(x)+|x-y|\le 5r$.
 Using \eqref{e:Phi0-scaling} (with $\eta=1/5$), we see that for all $s\in (0,t]$, 
 \begin{align*}
 	&s\left( 1 \wedge \frac{\delta_D(y)}{s^{1/\alpha}}\right)^q  
 	\Phi_0\bigg(\frac{\delta_D(y)\vee s^{1/\alpha}}{r}\bigg) =	s\left( \frac{\delta_D(y)}{\delta_D(y) \vee s^{1/\alpha}}\right)^q  
 	\Phi_0\bigg(\frac{\delta_D(y)\vee s^{1/\alpha}}{r}\bigg) \\
 	&\le c_4 	s\left( \frac{\delta_D(y)}{\delta_D(y) \vee t^{1/\alpha}}\right)^q  
 	\Phi_0\bigg(\frac{\delta_D(y)\vee t^{1/\alpha}}{r}\bigg)  \bigg( \frac{\delta_D(y) \vee s^{1/\alpha}}{\delta_D(y) \vee t^{1/\alpha}}\bigg)^{\beta_1-\eps-q}\\
 	&= c_4 	s\left( \frac{\delta_D(y)}{\delta_D(y) \vee t^{1/\alpha}}\right)^q 
 	\Phi_0\bigg(\frac{\delta_D(y)\vee t^{1/\alpha}}{r}\bigg)  \bigg( \frac{\delta_D(y) \vee t^{1/\alpha}}{\delta_D(y) \vee s^{1/\alpha}}\bigg)^{\alpha}\\
 	&= c_4t \left( \frac{\delta_D(y)}{\delta_D(y) \vee t^{1/\alpha}}\right)^q  
 	\Phi_0\bigg(\frac{\delta_D(y)\vee t^{1/\alpha}}{r}\bigg)  \times \begin{cases}
 		1 &\mbox{ if $\delta_D(y)\le s^{1/\alpha}$},\\
 		s /\delta_D(y)^{\alpha}&\mbox{ if $ s^{1/\alpha}<\delta_D(y)\le t^{1/\alpha}$},\\
 		s/t&\mbox{ if $\delta_D(y)> t^{1/\alpha}$}
 	\end{cases}\\
 	&\le c_4t \left( \frac{\delta_D(y)}{\delta_D(y) \vee t^{1/\alpha}}\right)^q  
 	\Phi_0\bigg(\frac{\delta_D(y)\vee t^{1/\alpha}}{r}\bigg) = c_4t \left( 1\wedge  \frac{\delta_D(y)}{ t^{1/\alpha}}\right)^q 
 	\Phi_0\bigg(\frac{\delta_D(y)\vee t^{1/\alpha}}{r}\bigg) .
 \end{align*}
 Combining this with \eqref{e:firstpart-1}, we obtain \eqref{e:firstpart-claim}. The proof is complete.
 \qed

\begin{lemma}\label{l:UHK-case1-main1}
Assume that $\beta_1^*<\alpha+\beta_1$. Let  $r\in (0,\text{\rm diam}(D))$ and  $x\in D$ satisfy $\delta_D(x)<5r$. 
 
 \noindent (i) There exists $C>0$  such that  for all  $t\in (0,r^\alpha]$ and $k>0$,
\begin{align*}
&	\int_{B_{ D}(x,r)}p(t,x,z)\Phi_1\left(\frac{\delta_D(z)}{r}\right)  \ell\left(\frac{\delta_D(z)}{k}\right) \, dz\\
		&\le C \left(1\wedge \frac{\delta_D(x)}{t^{1/\alpha}}\right)^q \Phi_1\bigg(\frac{\delta_D(x)\vee t^{1/\alpha}}{r}\bigg)\ell\bigg(\frac{\delta_D(x)\vee t^{1/\alpha}}{k}\bigg).
\end{align*}

\noindent (ii) There exists $C>0$  such that  for all  $t\in (0,r^\alpha]$ and $k>0$,
\begin{align*}
	&	\int_{B_{ D}(x,r)}p(t,x,z)\Phi_2\left(\frac{\delta_D(z)}{r}\right)  \ell\left(\frac{k}{\delta_D(z)}\right) \, dz\\
	&\le C \left(1\wedge \frac{\delta_D(x)}{t^{1/\alpha}}\right)^q \Phi_2\bigg(\frac{\delta_D(x)\vee t^{1/\alpha}}{r}\bigg)\ell\bigg(\frac{k}{\delta_D(x)\vee t^{1/\alpha}}\bigg)\\
&\quad + Ct \left(1\wedge \frac{\delta_D(x)}{t^{1/\alpha}}\right)^q   \int_{(\delta_D(x)\vee t^{1/\alpha})\wedge r}^{r}  \Phi_0\bigg(\frac{\delta_D(x)\vee t^{1/\alpha}}{l} \bigg)\Phi_2\bigg(\frac{l}{r}\bigg)\ell\bigg(\frac{k}{l}\bigg)\frac{dl}{l^{\alpha+1}}.
\end{align*}
 \end{lemma}
 \pf By Lemma \ref{l:UHK-case1-induction} and the almost increasing property of $\Phi_0$, we get for all $z\in B_D(x,r)$,
 \begin{align*}
 	p(t,x,z) &\le c_1 \bigg( 1 \wedge \frac{\delta_D(x) \wedge \delta_D(z)}{t^{1/\alpha}} \bigg)^q\Phi_0 \bigg( \frac{(\delta_D(x) \wedge \delta_D(z)) \vee t^{1/\alpha}}{|x-z|}\bigg) \bigg(t^{-d/\alpha} \wedge \frac{t}{|x-z|^{d+\alpha}} \bigg) \\
 	&\le c_2 \bigg( 1 \wedge \frac{\delta_D(x)}{t^{1/\alpha}} \bigg)^q\Phi_0 \bigg( \frac{\delta_D(x) \vee t^{1/\alpha}}{|x-z|}\bigg) \bigg(t^{-d/\alpha} \wedge \frac{t}{|x-z|^{d+\alpha}} \bigg) .
 \end{align*}
 Thus, applying Lemma \ref{l:analog-of-10.10}  (with  $\beta_1^*<\alpha+\beta_1$),  we obtain the desired results.
 \qed

For $r,t>0$ and $x,y \in D$, we define
 \begin{align}\label{e:sI-1}
 \begin{split}
 \mathcal{I}_{r,t}^{(1)}(x,y)=\mathcal{I}_{r,t}^{(1)}(x,y;\Phi_1,\Phi_2,\ell)&:=\int_0^{t/2}\int_{B_D(x,r)}\int_{B_D(y,r)}\1_{\{\delta_D(u)\le \delta_D(w)\}}p(s,x,u)p(t-s,y,w)\\ &\qquad \qquad \qquad \times \Phi_1\left(\frac{\delta_D(u)}{r}\right) \Phi_2\left(\frac{\delta_D(w)}{r}\right)	\ell\left(\frac{\delta_D(u)}{\delta_D(w)}\right) dw  du  ds
 \end{split}
 \end{align}
 and
 \begin{align}\label{e:sI-2}
 	\begin{split} 
   \mathcal{I}_{r,t}^{(2)}(x,y)= \mathcal{I}_{r,t}^{(2)}(x,y;\Phi_1,\Phi_2,\ell)&:=\int_0^{t/2}\int_{B_D(x,r)}\int_{B_D(y,r)}\1_{\{\delta_D(u)\ge \delta_D(w)\}}p(s,x,u)p(t-s,y,w)\\ &\qquad \qquad \qquad  \times \Phi_1\left(\frac{\delta_D(w)}{r}\right) \Phi_2\left(\frac{\delta_D(u)}{r}\right)	\ell\left(\frac{\delta_D(w)}{\delta_D(u)}\right) dw  du  ds.
\end{split}
 \end{align}

 \begin{lemma}\label{l:newlemma1} 
If $\beta_1^*<\alpha+\beta_1$, then
there exists $C>0$ such that   for all $r\in (0, \text{\rm diam}(D))$, $t\in (0,r^\alpha]$ and $x,y \in D$ with $\delta_D(x) \vee \delta_D(y)<5r$,
 	\begin{align*}
 	&	\sI_{r,t}^{(1)}(x,y)
 	 \le Ct \bigg(1 \wedge \frac{\delta_D(x)}{t^{1/\alpha}} \bigg)^q \bigg(1 \wedge \frac{\delta_D(y)}{t^{1/\alpha}} \bigg)^q	\\
 	 &\times \bigg[ \Phi_1\bigg(\frac{\delta_D(x)\vee t^{1/\alpha}}{r}\bigg) \Phi_2\bigg(\frac{\delta_D(y)\vee t^{1/\alpha}}{r}\bigg)	\ell\bigg(\frac{\delta_D(x)\vee t^{1/\alpha}}{\delta_D(y)\vee t^{1/\alpha}}\bigg)\\
 	&\qquad +  t \Phi_1\bigg(\frac{\delta_D(x)\vee t^{1/\alpha}}{r}\bigg) \int_{(\delta_D(y)\vee t^{1/\alpha})\wedge r}^{r}\Phi_0\bigg(\frac{\delta_D(y)\vee t^{1/\alpha}}{l} \bigg)  \Phi_2\bigg(\frac{l}{r}\bigg) \ell \bigg(\frac{\delta_D(x)\vee t^{1/\alpha}}{l}\bigg) \frac{dl}{l^{\alpha+1}}\bigg].
 	\end{align*}
 \end{lemma}
 \pf  By using   Lemma \ref{l:UHK-case1-main1}(ii), \eqref{e:Phi2-scaling} and \eqref{e:ell-scaling}  together with the fact that  
 $t-s\asymp t$  if  $0\le s\le t/2$,  we see that  for $0\le s\le t/2$ and $u\in B_D(x,r)$,  
 \begin{align*}
 	&\int_{B_D(y,r)} p(t-s, y,w)   \Phi_2 \left(\frac{\delta_D(w)}{r}\right) \ell\left(\frac{\delta_D(u)}{\delta_D(w)}\right) dw\\		&\le c_1 \left(1\wedge \frac{\delta_D(y)}{t^{1/\alpha}}\right)^q\Phi_2\bigg(\frac{\delta_D(y)\vee t^{1/\alpha}}{r}\bigg)\ell\bigg(\frac{\delta_D(u)}{\delta_D(y)\vee t^{1/\alpha}}\bigg) \\
 	&\quad  + c_1 t \left(1\wedge \frac{\delta_D(y)}{t^{1/\alpha}}\right)^q \int_{(\delta_D(y)\vee (t/2)^{1/\alpha})\wedge r}^{r}  \Phi_0\bigg(\frac{\delta_D(y)\vee t^{1/\alpha}}{l} \bigg)\Phi_2\bigg(\frac{l}{r}\bigg)\ell\bigg(\frac{\delta_D(u)}{l}\bigg)\frac{dl}{l^{\alpha+1}}.
 \end{align*}
 Notice that, by \eqref{e:Phi0-scaling},  \eqref{e:Phi2-scaling} and \eqref{e:ell-scaling}, 
 \begin{align}\label{e:newlemma1-1} 
 	\begin{split}
			&t \int_{(\delta_D(y)\vee (t/2)^{1/\alpha})\wedge r}^{(\delta_D(y)\vee t^{1/\alpha})\wedge r}  \Phi_0\bigg(\frac{\delta_D(y)\vee t^{1/\alpha}}{l} \bigg)\Phi_2\bigg(\frac{l}{r}\bigg)\ell\bigg(\frac{\delta_D(u)}{l}\bigg)\frac{dl}{l^{\alpha+1}}\\
 		&\le t \int_{\delta_D(y)\vee (t/2)^{1/\alpha}}^{\delta_D(y)\vee t^{1/\alpha}}  \Phi_0\bigg(\frac{\delta_D(y)\vee t^{1/\alpha}}{l} \bigg)\Phi_2\bigg(\frac{l}{r}\bigg)\ell\bigg(\frac{\delta_D(u)}{l}\bigg)\frac{dl}{l^{\alpha+1}}\\
 	 &\le c_2t \Phi_2\bigg(\frac{\delta_D(y)\vee t^{1/\alpha}}{r}\bigg)\ell\bigg(\frac{\delta_D(u)}{\delta_D(y) \vee t^{1/\alpha}}\bigg)\int_{\delta_D(y)\vee (t/2)^{1/\alpha}}^{\delta_D(y)\vee t^{1/\alpha}} 
	  \frac{dl}{l^{\alpha+1}}\\
 	 &\le \frac{c_3t}{(\delta_D(y)\vee t^{1/\alpha})^\alpha} \Phi_2\bigg(\frac{\delta_D(y)\vee t^{1/\alpha}}{r}\bigg)\ell\bigg(\frac{\delta_D(u)}{\delta_D(y) \vee t^{1/\alpha}}\bigg)\\
 	 &\le c_3 \Phi_2\bigg(\frac{\delta_D(y)\vee t^{1/\alpha}}{r}\bigg)\ell\bigg(\frac{\delta_D(u)}{\delta_D(y) \vee t^{1/\alpha}}\bigg).
 	\end{split}
 \end{align}
  Thus, by   Fubini's theorem,  Lemma  \ref{l:UHK-case1-main1}(ii),
  we get that for $0\le s\le t/2$,  
 \begin{align*}
 	&\int_{B_D(x,r)} \int_{B_D(y,r)}
 	p(s, x,u) 	p(t-s,y,w)  \Phi_1\bigg(\frac{\delta_D(u)}{r}\bigg) \Phi_2\bigg(\frac{\delta_D(w)}{r}\bigg) \ell\bigg(\frac{\delta_D(u)}{\delta_D(w)}\bigg) \, dw\, du\nn\\
 	&\le c_4 \bigg(1 \wedge \frac{\delta_D(y)}{t^{1/\alpha}} \bigg)^q \Phi_2\bigg(\frac{\delta_D(y)\vee t^{1/\alpha}}{r}\bigg)	\int_{B_D(x,r)}  	p(s, x,u) \Phi_1\bigg(\frac{\delta_D(u)}{r}\bigg)\ell\bigg(\frac{\delta_D(u)}{\delta_D(y) \vee t^{1/\alpha}}\bigg)\, du \nn \\ 
 	&\quad +c_4 t \left(1\wedge \frac{\delta_D(y)}{t^{1/\alpha}}\right)^q \int_{(\delta_D(y)\vee t^{1/\alpha})\wedge r}^{r}\Phi_0\bigg(\frac{\delta_D(y)\vee t^{1/\alpha}}{l} \bigg) \Phi_2\bigg(\frac{l}{r}\bigg)\\
 	&\qquad  \times   \int_{B_D(x,r)} p(s,x,u) \Phi_1\bigg( \frac{\delta_D(u)}{r}\bigg)\ell\bigg(\frac{\delta_D(u)}{l}\bigg) du \,\frac{dl}{l^{\alpha+1}} \\
 &\le c_5 \bigg(1 \wedge \frac{\delta_D(x)}{s^{1/\alpha}} \bigg)^q \bigg(1 \wedge \frac{\delta_D(y)}{t^{1/\alpha}} \bigg)^q	\Phi_1\bigg(\frac{\delta_D(x)\vee s^{1/\alpha}}{r}\bigg) \Phi_2\bigg(\frac{\delta_D(y)\vee t^{1/\alpha}}{r}\bigg)	\ell\bigg(\frac{\delta_D(x)\vee s^{1/\alpha}}{\delta_D(y)\vee t^{1/\alpha}}\bigg)\\
 &\quad +  c_5 t \bigg(1 \wedge \frac{\delta_D(x)}{s^{1/\alpha}} \bigg)^q \bigg(1 \wedge \frac{\delta_D(y)}{t^{1/\alpha}} \bigg)^q \Phi_1\bigg( \frac{\delta_D(x) \vee s^{1/\alpha}}{r} \bigg)\\
 &\qquad \times  \int_{(\delta_D(y)\vee t^{1/\alpha})\wedge r}^{r}\Phi_0\bigg(\frac{\delta_D(y)\vee t^{1/\alpha}}{l} \bigg)  \Phi_2\bigg(\frac{l}{r}\bigg) \ell \bigg(\frac{\delta_D(x)\vee s^{1/\alpha}}{l}\bigg) \frac{dl}{l^{\alpha+1}}.
 \end{align*}
Since $q<\alpha+\beta_1$, by applying Lemma \ref{cal:00}, we arrive at the desired result. \qed

 \begin{lemma}\label{l:newlemma1-case1} 
Let  $x,y \in D$ with $r:=|x-y|/4>\delta_D(x)\wedge  \delta_D(y)$, and let $t\in (0,r^\alpha]$. If $\beta_1^*<\alpha+\beta_1$, then  there exists $C>0$ independent of $x,y$ and $t$ such that 
 	\begin{align*}
 		&	\sI_{r,t}^{(1)}(x,y) + 	\sI_{r,t}^{(2)}(x,y)\\
 		& \le Ct \bigg(1 \wedge \frac{\delta_D(x)}{t^{1/\alpha}} \bigg)^q \bigg(1 \wedge \frac{\delta_D(y)}{t^{1/\alpha}}\bigg)^q \bigg[  A_{\Phi_1,\Phi_2,\ell}(t,x,y)\\
 		&\quad +  t\Phi_1\bigg(\frac{\delta_D(x)\vee t^{1/\alpha}}{r}\bigg)  \int_{(\delta_D(x)\vee \delta_D(y)\vee t^{1/\alpha})\wedge r}^{r}\Phi_0\bigg(\frac{\delta_D(y)\vee t^{1/\alpha}}{l} \bigg)   \Phi_2\bigg(\frac{l}{r}\bigg) \ell \bigg(\frac{\delta_D(x)\vee t^{1/\alpha}}{l}\bigg) \frac{dl}{l^{\alpha+1}}\\
 			&\quad +  t \Phi_1\bigg(\frac{\delta_D(y)\vee t^{1/\alpha}}{r}\bigg) \int_{(\delta_D(x)\vee \delta_D(y) \vee  t^{1/\alpha})\wedge r}^{r}\Phi_0\bigg(\frac{\delta_D(x)\vee t^{1/\alpha}}{l} \bigg)   \Phi_2\bigg(\frac{l}{r}\bigg) \ell \bigg(\frac{\delta_D(y)\vee t^{1/\alpha}}{l}\bigg) \frac{dl}{l^{\alpha+1}}\bigg].
 	\end{align*}
 \end{lemma}
 \pf  Note that $\delta_D(x) \vee \delta_D(y)\le \delta_D(x) \wedge  \delta_D(y)+ |x-y| <\frac54|x-y|$. Hence,
by \eqref{e:Phi1-scaling}--\eqref{e:ell-scaling},  there exist comparison constants independent of $x,y$ and $t$ such that
\begin{align}\label{e:newlemma1-case1-comparability} 
	\begin{split} 
&A_{\Phi_1,\Phi_2,\ell}(t,x,y)\\
&\asymp \Phi_1\bigg( \frac{(\delta_D(x)\wedge \delta_D(y)) \vee t^{1/\alpha}}{r}\bigg)   \Phi_2\bigg( \frac{\delta_D(x)\vee \delta_D(y)\vee t^{1/\alpha}}{r}\bigg) \ell\bigg(\frac{(\delta_D(x)\wedge\delta_D(y))\vee t^{1/\alpha}}{\delta_D(x)\vee \delta_D(y)\vee t^{1/\alpha}}\bigg).
\end{split} 
\end{align}
We distinguish between two cases.
 
 \smallskip
 
 \noindent
 \textbf{Case 1:} $\delta_D(x) \le \delta_D(y) \vee t^{1/\alpha}$. In this case, we have $(\delta_D(x)\wedge \delta_D(y)) \vee t^{1/\alpha} = \delta_D(x) \vee t^{1/\alpha}$ and $(\delta_D(x)\vee \delta_D(y)) \vee t^{1/\alpha} = \delta_D(y) \vee t^{1/\alpha}$. Hence, the upper bound for  $\sI_{r,t}^{(1)}(x,y; \Phi_1,\Phi_2,\ell)$
 follows from Lemma \ref{l:newlemma1} and \eqref{e:newlemma1-case1-comparability}. 
 
We decompose $\sI_{r,t}^{(2)}(x,y; \Phi_1,\Phi_2,\ell)$ as
$	\sI_{r,t}^{(2)}(x,y; \Phi_1,\Phi_2,\ell) = I_1+I_2$, 
where
\begin{align*}
	&	I_1:=\int_0^{t/2}\int_{B_D(x,r)}\int_{B_D(y,r)}\1_{\{\delta_D(w)\le \delta_D(u)< 2(\delta_D(y)\vee t^{1/\alpha})\}}  \, \cdots \,  dw\, du\, ds,\\
	&	I_2:=\int_0^{t/2}\int_{B_D(x,r)}\int_{B_D(y,r)}\1_{\{\delta_D(w) \le \delta_D(u), \,2(\delta_D(y)\vee t^{1/\alpha}) \le \delta_D(u)\}}  \, \cdots \,  dw\, du\, ds.
\end{align*}
Using Lemma \ref{l:Phi-ell-monotonicity}(ii) in the first inequality below, Lemma \ref{l:Phi-ell-monotonicity}(i) in the second, Lemma \ref{l:UHK-case1-main1}(i) and Corollary \ref{c:life} (with the fact that $t-s\asymp t$ for $s\in (0,t/2)$) in the third,  Lemma \ref{cal:00} (with $p<\alpha+\beta_1$) in the fourth, and \eqref{e:newlemma1-case1-comparability}  in the last,  we obtain
\begin{align*}
	I_1&\le c_1 \int_0^{t/2}\int_{B_D(x,r)}\int_{B_D(y,r)} \1_{\{\delta_D(w)\le \delta_D(u)\}}  p(s,x,u)p(t-s,y,w)\\
	&\qquad \times   \Phi_1\bigg(\frac{\delta_D(w)}{r}\bigg) \Phi_2 \bigg( \frac{\delta_D(y) \vee t^{1/\alpha}}{r}\bigg) \ell\bigg( \frac{\delta_D(w)}{\delta_D(y)\vee t^{1/\alpha}}\bigg) dw\, du\, ds \\
	&\le c_2 \Phi_2 \bigg( \frac{\delta_D(y) \vee t^{1/\alpha}}{r}\bigg) \\
	&\quad \times \int_0^{t/2}\int_{B_D(x,r)} p(s,x,u) \Phi_1\bigg(\frac{\delta_D(u)}{r}\bigg)  \ell\bigg( \frac{\delta_D(u)}{\delta_D(y) \vee t^{1/\alpha}}\bigg)  du\, \int_{B_D(y,r)} p(t-s,y,w) dw\, ds \\
	& \le c_3 \bigg( 1\wedge \frac{\delta_D(y)}{t^{1/\alpha}}\bigg)^q  \Phi_2 \bigg( \frac{\delta_D(y) \vee t^{1/\alpha}}{r}\bigg)\\
	&\quad \times  \int_0^{t/2} \bigg(1 \wedge \frac{\delta_D(x)}{s^{1/\alpha}}\bigg)^q \Phi_1\bigg(\frac{\delta_D(x) \vee s^{1/\alpha}}{r}\bigg) \ell\bigg( \frac{\delta_D(x) \vee s^{1/\alpha}}{\delta_D(y) \vee t^{1/\alpha}}\bigg)   ds\\
	& \le c_4t\bigg(1\wedge  \frac{\delta_D(x)}{t^{1/\alpha}}\bigg)^q \bigg(1\wedge  \frac{\delta_D(y)}{t^{1/\alpha}}\bigg)^q \Phi_1\bigg(\frac{\delta_D(x)\vee  t^{1/\alpha}}{r}\bigg)  \Phi_2 \bigg( \frac{\delta_D(y) \vee t^{1/\alpha}}{r}\bigg)\ell\bigg( \frac{\delta_D(x) \vee t^{1/\alpha}}{\delta_D(y) \vee t^{1/\alpha}}\bigg)   \\	& \le c_5t\bigg( 1\wedge  \frac{\delta_D(x)}{t^{1/\alpha}}\bigg)^q  \bigg(1\wedge  \frac{\delta_D(y)}{t^{1/\alpha}}\bigg)^q  A_{\Phi_1,\Phi_2,\ell}(t,x,y).
\end{align*}
For $s\in (0,t/2)$ and $u\in B_D(x,r)$ with $\delta_D(u)\ge 2(\delta_D(y)\vee t^{1/\alpha})$, since $\delta_D(x)\le \delta_D(y)\vee  t^{1/\alpha}$, we have  $|x-u| \ge \delta_D(u) - \delta_D(x) \ge \delta_D(u)/2 \ge \delta_D(y)\vee t^{1/\alpha}$.  Thus, by 
Lemmas \ref{l:UHK-case1-induction} and \ref{l:Phi-ell-monotonicity}(ii),  
we get
\begin{align}\label{e:case1-I2}
	\begin{split} 
	&p(s,x,u)  \Phi_2 \bigg( \frac{\delta_D(u)}{r}\bigg) \ell\bigg( \frac{\delta_D(y)\vee t^{1/\alpha}}{\delta_D(u)}\bigg)\\
	 &\le \frac{c_5t}{|x-u|^{d+\alpha}} \bigg(1\wedge \frac{\delta_D(x)}{s^{1/\alpha}}\bigg)^q \Phi_0 \bigg(\frac{\delta_D(x) \vee s^{1/\alpha}}{|x-u|}\bigg)\Phi_2 \bigg( \frac{|x-u|}{r}\bigg) \ell\bigg( \frac{\delta_D(y) \vee t^{1/\alpha}}{|x-u|}\bigg).
	\end{split}
\end{align}
Using Lemma \ref{l:UHK-case1-main1}(i) (with $t-s\asymp t$ for $s\in (0,t/2)$) in the first line below, \eqref{e:case1-I2} and the Fubini's theorem in the second, Lemma \ref{cal:00} in the third, we obtain
\begin{align*}
	I_2&\le c_6  \bigg(1\wedge  \frac{\delta_D(y)}{t^{1/\alpha}}\bigg)^q\Phi_1\bigg(\frac{\delta_D(y) \vee t^{1/\alpha}}{r}\bigg)\\
	&\quad \times  \int_0^{t/2}\int_{B_D(x,r):|x-u|\ge  \delta_D(y) \vee t^{1/\alpha}} p(s,x,u)   \Phi_2 \bigg( \frac{\delta_D(u)}{r}\bigg) \ell\bigg( \frac{\delta_D(y) \vee t^{1/\alpha}}{\delta_D(u)}\bigg) \, du\, ds\\
	&\le c_7  \bigg(1\wedge  \frac{\delta_D(y)}{t^{1/\alpha}}\bigg)^q\Phi_1\bigg(\frac{\delta_D(y) \vee t^{1/\alpha}}{r}\bigg)  \int_{B_D(x,r):|x-u|\ge \delta_D(y) \vee t^{1/\alpha}} \frac{t}{|x-u|^{d+\alpha}}  \\
	&\quad \times  \Phi_2 \bigg( \frac{|x-u|}{r}\bigg) \ell\bigg( \frac{\delta_D(y) \vee t^{1/\alpha}}{|x-u|}\bigg)\int_0^{t/2}  \bigg(1\wedge \frac{\delta_D(x)}{s^{1/\alpha}}\bigg)^q \Phi_0 \bigg(\frac{\delta_D(x) \vee s^{1/\alpha}}{|x-u|}\bigg) ds du\\
		&\le c_8t^2  \bigg( 1\wedge  \frac{\delta_D(x)}{t^{1/\alpha}}\bigg)^q  \bigg(1\wedge  \frac{\delta_D(y)}{t^{1/\alpha}}\bigg)^q\Phi_1\bigg(\frac{\delta_D(y) \vee t^{1/\alpha}}{r}\bigg)  \\
	&\quad \times \int_{B_D(x,r):|x-u|\ge  \delta_D(y) \vee t^{1/\alpha}} \Phi_0 \bigg(\frac{ \delta_D(x) \vee t^{1/\alpha}}{|x-u|}\bigg) \Phi_2 \bigg( \frac{|x-u|}{r}\bigg) \ell\bigg( \frac{\delta_D(y) \vee t^{1/\alpha}}{|x-u|}\bigg)     \frac{du}{|x-u|^{d+\alpha}}\\
		&\le c_9t^2  \bigg( 1\wedge  \frac{\delta_D(x)}{t^{1/\alpha}}\bigg)^q \bigg(1\wedge  \frac{\delta_D(y)}{t^{1/\alpha}}\bigg)^q\Phi_1\bigg(\frac{\delta_D(y) \vee t^{1/\alpha}}{r}\bigg) \\
		&\quad \times   \int_{( \delta_D(y) \vee t^{1/\alpha})\wedge r}^r \Phi_0 \bigg(\frac{\delta_D(x)\vee  t^{1/\alpha}}{l}\bigg)\Phi_2 \bigg( \frac{l}{r}\bigg) \ell\bigg( \frac{\delta_D(y)\vee t^{1/\alpha}}{l}\bigg)    \frac{dl}{l^{\alpha+1}}.
\end{align*}

 \noindent  \textbf{Case 2:} $\delta_D(x)  >\delta_D(y) \vee  t^{1/\alpha}$.  We decompose $\sI_{r,t}^{(2)}(x,y; \Phi_1,\Phi_2,\ell)$ as  $	\sI_{r,t}^{(2)}(x,y; \Phi_1,\Phi_2,\ell) = I_1'+I_2'$, 
 where
 \begin{align*}
 	&	I_1':=\int_0^{t/2}\int_{B_D(x,r)}\int_{B_D(y,r)}\1_{\{\delta_D(w)\le \delta_D(u)< 2\delta_D(x)\}}  \, \cdots \,  dw\, du\, ds,\\
 	&	I_2':=\int_0^{t/2}\int_{B_D(x,r)}\int_{B_D(y,r)}\1_{\{\delta_D(w) \le \delta_D(u), \,2\delta_D(x) \le \delta_D(u)\}}  \, \cdots \,  dw\, du\, ds.
 \end{align*} 
 For $I_1'$, using 
 Lemma \ref{l:Phi-ell-monotonicity}(ii)
  in the first inequality below, Lemma \ref{l:UHK-case1-main1}(i) (with $t-s\asymp t$ for $s\in (0,t/2)$, \eqref{e:Phi1-scaling} and \eqref{e:ell-scaling}) in the second, and \eqref{e:newlemma1-case1-comparability}  in the fourth, we get that 
 \begin{align*}
 	&I_1'\le c_{10}\int_0^{t/2}\int_{B_D(x,r)}p(s,x,u) \Phi_2 \bigg( \frac{\delta_D(x)}{r}\bigg) \int_{B_D(y,r)}  p(t-s,y,w)  \Phi_1\bigg(\frac{\delta_D(w)}{r}\bigg) \ell\bigg( \frac{\delta_D(w)}{\delta_D(x)}\bigg) dw du ds \\
 	&\le  c_{11}\bigg( 1 \wedge \frac{\delta_D(y)}{t^{1/\alpha}} \bigg)^q \Phi_1\bigg(\frac{\delta_D(y)\vee t^{1/\alpha}}{r}\bigg) \Phi_2 \bigg( \frac{\delta_D(x)}{r}\bigg)  \ell\bigg( \frac{\delta_D(y) \vee t^{1/\alpha}}{\delta_D(x)}\bigg) \int_0^{t/2}\int_{B_D(x,r)}p(s,x,u)  du ds \\
 	& \le\frac{c_{11}t}{2}\bigg( 1 \wedge \frac{\delta_D(y)}{t^{1/\alpha}} \bigg)^q \Phi_1\bigg(\frac{\delta_D(y)\vee t^{1/\alpha}}{r}\bigg) \Phi_2 \bigg( \frac{\delta_D(x)}{r}\bigg)  \ell\bigg( \frac{\delta_D(y) \vee t^{1/\alpha}}{\delta_D(x)}\bigg) \\	
 	& \le c_{12}t \bigg(1\wedge  \frac{\delta_D(y)}{t^{1/\alpha}}\bigg)^q  A_{\Phi_1,\Phi_2,\ell}(t,x,y).
 \end{align*}
 For $s\in (0,t/2)$ and $u\in B_D(x,r)$ with $\delta_D(u)\ge 2\delta_D(x)$,  we have  $|x-u| \ge  \delta_D(u)/2 \ge \delta_D(x)$.  
 Thus, by 
 Lemmas \ref{l:UHK-case1-induction} and \ref{l:Phi-ell-monotonicity}(ii),  
  we get
 \begin{align}\label{e:case2-I2}
 	\begin{split} 
 		&p(s,x,u)  \Phi_2 \bigg( \frac{\delta_D(u)}{r}\bigg) \ell\bigg( \frac{\delta_D(y)\vee t^{1/\alpha}}{\delta_D(u)}\bigg)\\
 		&\le \frac{c_{13}t}{|x-u|^{d+\alpha}}  \Phi_0 \bigg(\frac{\delta_D(x) }{|x-u|}\bigg)\Phi_2 \bigg( \frac{|x-u|}{r}\bigg) \ell\bigg( \frac{\delta_D(y) \vee t^{1/\alpha}}{|x-u|}\bigg).
 	\end{split}
 \end{align}
 Using Lemma \ref{l:UHK-case1-main1}(i) (with $t-s\asymp t$ for $s\in (0,t/2)$) in the first line below, and \eqref{e:case2-I2} in the second, we obtain
 \begin{align*}
 	I_2'&\le c_{14}  \bigg(1\wedge  \frac{\delta_D(y)}{t^{1/\alpha}}\bigg)^q \Phi_1\bigg(\frac{\delta_D(y) \vee t^{1/\alpha}}{r}\bigg)\\
 	&\quad \times  \int_0^{t/2}\int_{B_D(x,r):|x-u|\ge  \delta_D(x)} p(s,x,u)   \Phi_2 \bigg( \frac{\delta_D(u)}{r}\bigg) \ell\bigg( \frac{\delta_D(y) \vee t^{1/\alpha}}{\delta_D(u)}\bigg) \, du\, ds\\
 	&\le c_{15}  \bigg(1\wedge  \frac{\delta_D(y)}{t^{1/\alpha}}\bigg)^q \Phi_1\bigg(\frac{\delta_D(y) \vee t^{1/\alpha}}{r}\bigg)  \int_0^{t/2}  ds \\
 	&\quad \times\int_{B_D(x,r):|x-u|\ge \delta_D(x)} \frac{t}{|x-u|^{d+\alpha}}   \Phi_0 \bigg(\frac{\delta_D(x)}{|x-u|}\bigg)  \Phi_2 \bigg( \frac{|x-u|}{r}\bigg) \ell\bigg( \frac{\delta_D(y) \vee t^{1/\alpha}}{|x-u|}\bigg) du\\
 	&\le c_{16}t^2   \bigg(1\wedge  \frac{\delta_D(y)}{t^{1/\alpha}}\bigg)^q\Phi_1\bigg(\frac{\delta_D(y) \vee t^{1/\alpha}}{r}\bigg)  \int_{ \delta_D(x) \wedge r}^r \Phi_0 \bigg(\frac{\delta_D(x)}{l}\bigg)\Phi_2 \bigg( \frac{l}{r}\bigg) \ell\bigg( \frac{\delta_D(y)\vee t^{1/\alpha}}{l}\bigg)    \frac{dl}{l^{\alpha+1}}.
 \end{align*}
 
 For $\sI_{r,t}^{(1)}(x,y; \Phi_1,\Phi_2,\ell)$, we decompose it as   $	\sI_{r,t}^{(1)}(x,y; \Phi_1,\Phi_2,\ell) = I_3'+I_4'$, 
 where
 \begin{align*}
 	&	I_3':=\int_0^{t/2}\int_{B_D(x,r)}\int_{B_D(y,r)}\1_{\{\delta_D(u)\le \delta_D(w)< 2\delta_D(x)\}}  \, \cdots \,  dw\, du\, ds,\\
 	&	I_4':=\int_0^{t/2}\int_{B_D(x,r)}\int_{B_D(y,r)}\1_{\{\delta_D(u) \le \delta_D(w), \,2\delta_D(x) \le \delta_D(w)\}}  \, \cdots \,  dw\, du\, ds.
 \end{align*} 
 Using Lemma \ref{l:Phi-ell-monotonicity}(ii) with the scaling properties of $\Phi_2$ and $\ell$ in the first inequality below, Lemma \ref{l:Phi-ell-monotonicity}(i) in the second, Lemma \ref{l:UHK-case1-main1}(i) (with  $t-s\asymp t$ for $s\in (0,t/2)$)   in the third,  and \eqref{e:newlemma1-case1-comparability}  in the last,  we obtain
 \begin{align*}
 &I_3'\le c_{16} \int_0^{t/2}\int_{B_D(x,r)}\int_{B_D(y,r)} \1_{\{\delta_D(u)\le \delta_D(w)\}}  p(s,x,u)p(t-s,y,w)\\
 &\qquad \times   \Phi_1\bigg(\frac{\delta_D(u)}{r}\bigg) \Phi_2 \bigg( \frac{\delta_D(x)}{r}\bigg) \ell\bigg( \frac{\delta_D(u)}{\delta_D(x)}\bigg) dw\, du\, ds \\
 &\le c_{17} \Phi_2 \bigg( \frac{\delta_D(x)}{r}\bigg)  \int_0^{t/2}\int_{B_D(x,r)} p(s,x,u)   du\, \int_{B_D(y,r)} p(t-s,y,w) \Phi_1\bigg(\frac{\delta_D(w)}{r}\bigg)  \ell\bigg( \frac{\delta_D(w)}{\delta_D(x)}\bigg) dw\, ds \\
 & \le c_{18} \bigg( 1\wedge \frac{\delta_D(y)}{t^{1/\alpha}}\bigg)^q \Phi_1\bigg(\frac{\delta_D(y) \vee t^{1/\alpha}}{r}\bigg) \Phi_2 \bigg( \frac{\delta_D(x) }{r}\bigg)\ell\bigg( \frac{\delta_D(y) \vee t^{1/\alpha}}{\delta_D(x)}\bigg)  \int_0^{t/2}      ds\\	& \le c_{19}t \bigg(1\wedge  \frac{\delta_D(y)}{t^{1/\alpha}}\bigg)^q A_{\Phi_1,\Phi_2,\ell}(t,x,y).
 \end{align*}
For $s\in (0,t/2)$ and $w\in B_D(y,r)$ with $\delta_D(w)\ge 2\delta_D(x)$, since $|y-w|\ge \delta_D(w)/2$, we see from 
Lemmas \ref{l:UHK-case1-induction} and \ref{l:Phi-ell-monotonicity}(ii),  
that 
 \begin{align}\label{e:case2-I4}
 	\begin{split}
 	&p(t-s,y,w) \Phi_2 \bigg( \frac{\delta_D(w)}{r}\bigg) \ell\bigg( \frac{\delta_D(x)}{\delta_D(w)}\bigg)\\
 	& \le \frac{c_{20}t}{|y-w|^{d+\alpha}}  \bigg(1\wedge  \frac{\delta_D(y)}{t^{1/\alpha}}\bigg)^q \Phi_0\bigg( \frac{\delta_D(y)\vee t^{1/\alpha}}{|y-w|} \bigg) \Phi_2 \bigg( \frac{|y-w|}{r}\bigg)  \ell\bigg( \frac{\delta_D(x)}{|y-w|}\bigg) .
 	\end{split} 
 \end{align}
 Using Lemma \ref{l:UHK-case1-main1}(i) in the first inequality, \eqref{e:case2-I4} in the second
 \begin{align*}
 	I_4'&\le c_{21}\Phi_1\bigg(\frac{\delta_D(x)}{r}\bigg)  \int_0^{t/2}\int_{B_D(y,r):|y-w|\ge 
	 \delta_D(w)/2} 
	p(t-s,y,w) \Phi_2 \bigg( \frac{\delta_D(w)}{r}\bigg) \ell\bigg( \frac{\delta_D(x)}{\delta_D(w)}\bigg) dw\,  ds \\
 	&\le c_{22} t\bigg(1\wedge  \frac{\delta_D(y)}{t^{1/\alpha}}\bigg)^q\Phi_1\bigg(\frac{\delta_D(x)}{r}\bigg) \int_0^{t/2}ds\\
 	&\quad \times\int_{B_D(y,r):|y-w|\ge \delta_D(x)}  \Phi_0 \bigg( \frac{\delta_D(y)\vee t^{1/\alpha}}{|y-w|}\bigg) \Phi_2 \bigg( \frac{|y-w|}{r}\bigg) \ell\bigg( \frac{\delta_D(x)}{|y-w|}\bigg) \frac{dw}{|y-w|^{d+\alpha}}
 \\
 	&\le  c_{23} t^2\bigg(1\wedge  \frac{\delta_D(y)}{t^{1/\alpha}}\bigg)^q\Phi_1\bigg(\frac{\delta_D(x)}{r}\bigg) \int_{\delta_D(x)\wedge r}^r  \Phi_0\bigg( \frac{\delta_D(y)\vee t^{1/\alpha}}{l} \bigg) \Phi_2 \bigg( \frac{l}{r}\bigg)  \ell\bigg( \frac{\delta_D(x)}{l}\bigg) \frac{dl}{l^{\alpha+1}}.
 \end{align*}
 
 The proof is complete.  \qed

 Recall that the constant $\eps_1$ is defined in \eqref{e:def-eps1}.
 
 \begin{thm}\label{t:UHK}
	If $\beta_1^* <\alpha+\beta_1$, then for
	any $T>0$, there exists  $C=C(T)>0$ such that for for all $t\in(0,T]$ and $x,y \in D$,
 	\begin{align}\label{e:UHK-1}
 		\begin{split} 
 		&p(t,x,y)\le C \left(1 \wedge \frac{\delta_D(x)}{t^{1/\alpha}}\right)^q \left(1 \wedge \frac{\delta_D(y)}{t^{1/\alpha}}\right)^q \left( t^{-d/\alpha} \wedge \frac{t}{|x-y|^{d+\alpha}}\right)\bigg[ A_{\Phi_1, \Phi_2,  \ell} (t,x,y) \\
 		&    + \big( t \wedge |x-y|^\alpha\big)  \int_{(\delta_D(x) \vee \delta_D(y)\vee t^{1/\alpha}) \wedge (\eps_1|x-y|/2)}^{\eps_1|x-y|}A_{\Phi_1,\Phi_2,\ell}(t,x,x+u\mathbf{n}_x)\,A_{\Phi_1,\Phi_2,\ell}(t,x+u\mathbf{n}_x,y)\frac{du}{u^{\alpha+1}}\\
 			& + \big( t \wedge |x-y|^\alpha\big) \int_{(\delta_D(x) \vee \delta_D(y)\vee t^{1/\alpha}) \wedge (\eps_1|x-y|/2)}^{\eps_1|x-y|}A_{\Phi_1,\Phi_2,\ell}(t,y,y+u\mathbf{n}_y)\,A_{\Phi_1,\Phi_2,\ell}(t,y+u\mathbf{n}_y,x)\frac{du}{u^{\alpha+1}} \bigg] .
 		\end{split} 
 	\end{align}
 \end{thm}
 \pf By symmetry, we  assume without loss of generality that $\delta_D(x) \le \delta_D(y)$. Let $r:=|x-y|/4$.
 If $\delta_D(x) \vee t^{1/\alpha} \ge  r$, then   
 \eqref{e:UHK-1} 
 follows from Proposition \ref{p:UHK-rough} and \eqref{e:interior-lower-bound-A}. 
 
 Suppose $\delta_D(x) \vee t^{1/\alpha} <r$.   Let $V_1:=B_{\overline D}(x,r)$, $V_3:= B_{\overline D}(y,r)$ and
 $V_2:=\overline{D} \setminus (V_1 \cup V_3)$.  Note that $\delta_D(y)\le \delta_D(x)+4r< 5r$, and $2r< |u-w|< 6r$ for all $u\in V_1$ and $w\in V_3$.
 By Lemmas \ref{l:firstpart} and \ref{l:two-jumps-lower-bound}(i),  we have 	\begin{align}\label{e:UHK-firstpart}
 	\begin{split} 	&\P_x(\tau_{V_1}<t< \zeta) \sup_{s \le t, \, z \in V_2} p(s,z,y)\le \P_x(\tau_{V_1}<t< \zeta) \sup_{s \le t, \, z \in D,  |z-y|\ge r} p(s,z,y)  \\	&\le \frac{c_2t^2}{r^{d+2\alpha}} \bigg(1 \wedge \frac{\delta_D(x)}{t^{1/\alpha}}\bigg)^q \bigg( 1 \wedge \frac{\delta_D(y)}{t^{1/\alpha}}\bigg)^q  A_{\Phi_0,\Phi_0,1}(t,x,y)\\
 		&\le  \frac{c_3t^2}{r^{d+\alpha}} \bigg(1 \wedge \frac{\delta_D(x)}{t^{1/\alpha}}\bigg)^q \bigg( 1 \wedge \frac{\delta_D(y)}{t^{1/\alpha}}\bigg)^q \\
 		&\quad \times \int_{(\delta_D(x) \vee \delta_D(y)\vee t^{1/\alpha}) \wedge (\eps_1|x-y|/2)}^{\eps_1|x-y|}A_{\Phi_1,\Phi_2,\ell}(t,x,x+u\mathbf{n}_x)A_{\Phi_1,\Phi_2,\ell}(t,x+u\mathbf{n}_x,y)\frac{du}{u^{\alpha+1}}.
 		\end{split} 
 \end{align}
 
 Define
 \begin{align*}
 		\sI&:= \int_0^t \int_{V_3}\int_{V_1} p^{V_1}(s,x,u) \sB(u,w) p(t-s,y,w)\,du\,dw\,ds.
 \end{align*}
 By using  \hyperlink{A1}{{\bf (A1)}} and \eqref{e:Phi1-scaling}--\eqref{e:ell-scaling} with the fact that $2r<|u-w|<6r$ for all $u \in V_1$ and $w\in V_3$,  we see that 
 \begin{align*}
 	\sI&\le c_4 \int_0^t \int_{V_1}\int_{V_3}\1_{\{\delta_D(u)\le \delta_D(w)\}} p(s,x,u)p(t-s,y,w) \\
 	&\qquad \times \Phi_1\bigg(\frac{\delta_D(u)}{r}\bigg)\Phi_2\bigg(\frac{ \delta_D(w)}{r}\bigg) 
 	\ell\bigg(\frac{\delta_D(u)}{\delta_D(w)}\bigg)dw\,du\,ds \\
 	&+ c_4 \int_0^t \int_{V_1}\int_{V_3}\1_{\{\delta_D(w)\le \delta_D(u)\}} p(s,x,u)p(t-s,y,w) \\
 	&\qquad \times\Phi_1\bigg(\frac{\delta_D(w)}{r}\bigg)\Phi_2\bigg(\frac{ \delta_D(u)}{r}\bigg)
 	\ell\bigg(\frac{\delta_D(w)}{\delta_D(u)}\bigg)dw\,du\,ds\\
 	& = c_4 \bigg(\int_0^{t/2}+ \int_{t/2}^{t} \bigg)
 	\int_{V_1} \int_{V_3} \1_{\{\delta_D(u)\le \delta_D(w)\}} \dots dw\,du\,ds\\
 	&+ c_4
 	\bigg(\int_0^{t/2}+ 
 	\int_{t/2}^{t} \bigg) \int_{V_1} \int_{V_3} \1_{\{\delta_D(w)\le \delta_D(u)\}} 
 	\dots dw\,du\,ds.
 \end{align*}
 Thus, by the change of variables $\wt s = t-s$ in the integrals $\int_{t/2}^t$, we get
 \begin{align*}
 	\sI	 \le c_4\sum_{i=1}^2	
	\mathcal{I}_{r,t}^{(i)}(x,y)
	+c_4\sum_{i=1}^2	
	 \mathcal{I}_{r,t}^{(i)}(y,x),
 \end{align*}
 where the functions $\mathcal{I}_{r,t}^{(i)}(x,y)$,  $1\le i \le 2$, are defined in \eqref{e:sI-1}--\eqref{e:sI-2}.  By  
 Lemmas \ref{l:newlemma1-case1} and \ref{l:two-jumps-equivalent},
 we deduce that 
 $|x-y|^{-d-\alpha}\sI$ 
 is bounded above by the right-hand side of \eqref{e:UHK-1}. 
 The result now follows from
  Lemma \ref{l:general-upper-2} and the fact that dist$(V_1,V_3)\ge |x-y|/2$.
  \qed


\section{Proofs of Theorems \ref{t:main}, \ref{t:main11}, and \ref{t:main2}}\label{s:proofs}

\noindent

\noindent \textbf{Proofs of Theorems \ref{t:main}(i) and \ref{t:main2}(i).}   
(a) If $t^{1/\alpha}\ge \eps_1 |x-y|/2$, then the results follow from  Proposition \ref{p:UHK-rough}  and  Lemma \ref{l:HKE-lower-1}.

\noindent (b)
 Suppose  $t^{1/\alpha}< \eps_1 |x-y|/2$. Then $t^{-d/\alpha}\wedge (t|x-y|^{-d-\alpha})\asymp t|x-y|^{-d-\alpha}$ and $t\wedge |x-y|^{\alpha}\asymp t$. Hence,  the lower bound in \eqref{e:HKE-off} follows from Theorem \ref{t:LHK}, and the upper bound follows from Theorem \ref{t:UHK} and  Lemma \ref{l:equivalences}. 

The second comparison in Theorem \ref{t:main2}(i) follows from Corollary \ref{c:life-2}. 
 \qed

\medskip
\noindent
\textbf{Proof of Theorem \ref{t:main11}.} (i)-(ii)  The lower 
bounds follow from Theorem \ref{t:LHK} combined with Lemma \ref{l:two-jumps-lower-bound}, and the upper 
bounds follow from Theorem \ref{t:UHK} combined with Lemma \ref{l:two-jumps-comparability}.  

\noindent  (iii) When $t^{1/\alpha} \ge \eps_1|x-y|/2$, using \eqref{e:upper-bound-A}, \eqref{e:interior-lower-bound-A} and \eqref{e:ell-scaling}, we see that the right-hand side of \eqref{e:main11-case3} is comparable with $t^{-d/\alpha}$. Thus, the result follows from Theorem \ref{t:main}(i). If $t^{1/\alpha} < \eps_1|x-y|/2$, then a direct computation, together with Theorem \ref{t:main}(i) and the comparability of the first and third terms in Lemma \ref{l:equivalences}, yields the desired result. \qed

\medskip
\noindent
\textbf{Proofs of Theorems \ref{t:main}(ii)  and  \ref{t:main2}(ii).}
By Proposition \ref{p:upper-heatkernel}, for anyy  $t>0$, there exists a constant $c_1(t)>0$ depending on $t$ such that for all $x,y\in D$,
\begin{align}\label{e:main2-0}
 p^{\kappa}(t,x,y) \le \bar{p}(t,x,y)\le c_1(t).
\end{align}
Thus, since $D$ is bounded,  both 
$(\overline P_t)_{t\ge 0}$ and 
$(P^{\kappa}_t)_{t\ge 0}$ are compact semigroups  on $L^2(D)$. Note that  the largest eigenvalue  of the infinitesimal generator of $(\overline P_t)_{t\ge 0}$ is 0, 
and the corresponding eigenfunction can be take to be $\overline \phi_1(x)=1$.
Since $(P^{\kappa}_t)_{t\ge 0}$ is non-conservative and positive preserving, 
 by  the Krein-Rutman theorem, 
 the largest eigenvalue $-\lambda_1$ of its infinitesimal generator is negative, 
 and the corresponding eigenfunction $\phi_1$ can be taken to be positive a.e. on $D$.

By Theorem \ref{t:main2}(i), for any  $t>0$, there exists a constant $c_2(t)\ge 1$ such that for all $x,y\in D$,
\begin{align}\label{e:main2-1}
	& c_2(t)^{-1}\left(1\wedge \frac{\delta_D(x)}{t^{1/\alpha}}\right)^q \left(1\wedge \frac{\delta_D(y)}{t^{1/\alpha}}\right)^q \bar{p}(t,x,y)\\ 
	&\le p^{\kappa}(t,x,y) \le c_2(t)\left(1\wedge \frac{\delta_D(x)}{t^{1/\alpha}}\right)^q \left(1\wedge \frac{\delta_D(y)}{t^{1/\alpha}}\right)^q \bar{p}(t,x,y).\nn
\end{align}
Recall from \cite[Remark 4.9]{CKSV24} that $\overline{P}_1$ is strongly Feller. Hence, $\overline{P}_1\big(\delta_D(\cdot)^p \phi_1\big)$ is positive and continuous on $\overline{D}$, and therefore $\overline{P}_1\big(\delta_D(\cdot)^q \phi_1\big)(x) \asymp 1$  for $x\in D$. Further, since $D$ is bounded, we have $1\wedge \delta_D(x) \asymp \delta_D(x)$ for $x\in D$. By taking $t=1$ in \eqref{e:main2-1}, we get that for $x\in D$,
\begin{align}\label{e:main2-2}
	\phi_1(x)&=e^{-\lambda_1}P_1^{\kappa}\phi(x)\asymp \int_D e^{-\lambda_1}\delta_D(x)^q \delta_D(y)^q\, \bar{p}(1,x,y)\ dy\\
	& = e^{-\lambda_1}  \delta_D(x)^q\, \overline{P}_1\big(\delta_D(\cdot)^q \phi_1\big)(x)\asymp \delta_D(x)^q. \nn
\end{align} 
Combining   \eqref{e:main2-0}, \eqref{e:main2-1} and \eqref{e:main2-2}, we get that for any $t>0$, 
$$
p^{\kappa}(t,x,y)\le c_1(t)c_2(t)  t^{-2q/\alpha} \delta_D(x)^q \delta_D(y)^q 
\le c_3(t) \phi_1(x)\phi_1(y).
$$
Hence, the semigroup $(P^{\kappa}_t)_{t\ge 0}$ is intrinsically ultracontractive (cf. \cite[Theorem 3.2]{DS84}). 

Notice that $(\overline P_t)_{t\ge 0}$ is also  intrinsically ultracontractive by \eqref{e:main2-0}. Now, it follows from \cite[Theorem 4.2.5]{Dav89} and \eqref{e:main2-2} 
that for any $T>0$, there exist comparison constants depending on $T$ such that for all $t\ge T$ and $x,y \in  D$,
\begin{align*}
	\overline p(t,x,y) \asymp \overline \phi_1(x) \overline \phi_1(y)= 
	1 \quad \text{and} \quad 
	p^{\kappa}(t,x,y) \asymp  e^{-\lambda_1 t} \phi_1(x) \phi_1(y) \asymp e^{-\lambda_1 t}\delta_D(x)^q \delta_D(y)^q.
\end{align*}
By the continuity of $\overline p(t,x,y)$, this completes the proofs. \qed


\section{Auxiliary results}\label{s:aux} 

Recall that the classes $\sM(\gamma,\gamma^*)$ and $\sM^\uparrow(\gamma,\gamma^*)$ are defined at the beginning of Section \ref{s:preliminary}.

Throughout this section, we let $D\subset \R^d$, $d\ge 2$, be a bounded Lipschitz open set,  let   $\Phi\in \sM^\uparrow(\beta,\beta^*)$ and  $\Psi\in \sM^\uparrow(\gamma,\gamma^*)$  for some $\beta^*\ge \beta\ge 0$ and $\gamma^*\ge \gamma\ge 0$, and let $\ell\in \sM(0,0)$. Observe that   for any $R\ge 1$, $\eta\in (0,1]$ and  $\eps>0$,   there exists $C= C(R,\eta, \eps)\ge 1$ such that 
\begin{align}
	C^{-1}\bigg( \frac{r}{s}\bigg)^{(\beta-\eps)_+}
	&\le \frac{\Phi(r)}{\Phi(s)}\le C
	\bigg( \frac{r}{s}\bigg)^{\beta^*+\eps} \quad \text{for all $0<\eta s\le r\le R$},\label{e:Phi-aux-scaling}\\
		C^{-1}\bigg( \frac{r}{s}\bigg)^{(\gamma-\eps)_+}
	&\le \frac{\Psi(r)}{\Psi(s)}\le C
	\bigg( \frac{r}{s}\bigg)^{\gamma^*+\eps} \quad \text{for all $0<\eta s\le r\le R$},\label{e:Psi-aux-scaling}
	\\	
	C^{-1} \bigg( \frac{r}{s}\bigg)^{ -\eps}
	&\le \frac{\ell(r)}{\ell(s)} 
	\le C \bigg( \frac{r}{s}\bigg)^{  \eps } \qquad \;\;\;\text{for all $0<\eta s\le r$}.\label{e:ell-aux-scaling}
\end{align}

\begin{lemma}\label{l:Phi-ell-monotonicity}
(i) Assume further that $\ell$ is almost increasing on $(0,1]$	if $\beta=0$. Then for any $a\ge 1$, there exists $C=C(a)>0$ such that
\begin{align*}
	\Phi(s)\ell(s/k) \le C\Phi(r)\ell(r/k) \quad \text{for all $0<s\le ar\le 1$ and $k>0$.}
\end{align*}

\noindent (ii)	Assume further that $\ell$ is almost decreasing on $(0,1]$	if $\gamma=0$.  Then for any $a\ge 1$, there exists $C>0$ such that
\begin{align*}
	\Phi(s)\ell(k/s) \le C\Phi(r)\ell(k/r) \quad \text{for all $0<s\le ar\le 1$ and $k>0$.}
\end{align*}
\end{lemma}
\pf Since the proofs are similar, we only prove (i).  Since $\Phi(r)\ell(r/k) \asymp \Phi(ar)\ell(ar/k)$ for $0<r\le 1/a$ by \eqref{e:Phi-aux-scaling} and \eqref{e:ell-aux-scaling}, it suffices to show that
\begin{align}\label{e:Phi-ell-monotonicity-claim}
		\Phi(s)\ell(k/s) \le C\Phi(r)\ell(k/r) \quad \text{for all $0<s\le r\le 1$ and $k>0$.}
\end{align}
 If $\beta=0$, then \eqref{e:Phi-ell-monotonicity-claim} follows from the almost increasing properties of $\Phi$ and $\ell$. Assume $\beta>0$ and let $\eps:=\beta/2$. Using \eqref{e:Phi-aux-scaling} and \eqref{e:ell-aux-scaling}, we obtain for all $0<s\le r\le 1$ and $k>0$,
\begin{align*}
	\frac{\Phi(s)\ell(s/k)}{\Phi(r)\ell(r/k)} \le c_1 \bigg( \frac{s}{r}\bigg)^{\beta-2\eps} = c_1.
\end{align*}
\qed

\begin{lemma}\label{cal:00}
 Let  $q\ge 0$, $r>0$, $0<t\le r^\alpha$, $k>0$  and $m>0$. 
 Suppose that either $k\ge t^{1/\alpha}$ or $q<\alpha+\beta$. Then there exists $C>0$, independent of $r,t,k$ and $m$,
 such that
 \begin{align}\label{e:cal:00-estimate-I}
 	&\int_0^t \left(1\wedge \frac{k}{s^{1/\alpha}}\right)^{q}\Phi\bigg(\frac{k\vee s^{1/\alpha}}{r}\bigg)\ell\bigg(\frac{m}{k\vee s^{1/\alpha}}\bigg) \, ds\le C t \left(1\wedge \frac{k}{t^{1/\alpha}}\right)^{q}\Phi\bigg(\frac{k\vee t^{1/\alpha}}{r}\bigg) \ell\bigg(\frac{m}{k\vee t^{1/\alpha}}\bigg)
 \end{align}
 and
\begin{align}\label{e:cal:00-estimate-J}
&\int_0^t \left(1\wedge \frac{k}{s^{1/\alpha}}\right)^{q}\Phi\bigg(\frac{k\vee s^{1/\alpha}}{r}\bigg) \ell\bigg(\frac{k\vee s^{1/\alpha}}{m}\bigg) \, ds \le C t \left(1\wedge \frac{k}{t^{1/\alpha}}\right)^{q}\Phi\bigg(\frac{k\vee t^{1/\alpha}}{r}\bigg)\ell\bigg(\frac{k\vee t^{1/\alpha}}{m}\bigg).
\end{align}
\end{lemma}
\pf  Since the proofs are similar, we only prove \eqref{e:cal:00-estimate-I}. The inequality clearly holds if $k\ge t^{1/\alpha}$. 
Assume  $k< t^{1/\alpha}$ and $q<\alpha+\beta$. Then the left-hand side of  \eqref{e:cal:00-estimate-I} is equal to
\begin{align*}
&\Phi(k/r) \ell(m/k)\int_0^{k^{\alpha}}    ds + k^q \int_{k^{\alpha}}^t s^{-q/\alpha}\Phi(s^{1/\alpha}/r)\ell(m/s^{1/\alpha}) ds=:I_1+I_2.
\end{align*}
Let $\eps>0$ be such that $q<\alpha+\beta - 2\eps$. . 
By  \eqref{e:ell-aux-scaling} and \eqref{e:Phi-aux-scaling},  we have 
\begin{align}\label{e:cal:00-1}
\frac{	\Phi(t^{1/\alpha}/r) \ell (m/t^{1/\alpha})}{ \Phi(s^{1/\alpha}/r) \ell (m/s^{1/\alpha})} \ge c_1 \bigg(\frac{t}{s}\bigg)^{(\beta-2\eps)/\alpha} \quad \text{for all $0<s\le t$.} 
\end{align}
Using \eqref{e:cal:00-1} with $s=k^\alpha$, since  $q<\alpha+\beta-2\eps$, we obtain that
\begin{align*}
	\frac{k^q t^{(\alpha-q)/\alpha} \Phi(t^{1/\alpha}/r) \ell(m/t^{1/\alpha})}{I_1}& =\frac{ \Phi(t^{1/\alpha}/r) \ell(m/t^{1/\alpha})}{\Phi(k/r) \ell(m/k)} \bigg(\frac{t^{1/\alpha}}{k}\bigg)^{\alpha-q} \ge c_1\bigg(\frac{t^{1/\alpha}}{k}\bigg)^{\beta-2\eps + \alpha -q} \ge c_1.
\end{align*}
This gives the correct bound for $I_1$. 
By  using \eqref{e:cal:00-1}, we also get that
\begin{align*}
I_2 &\le c_1^{-1} 	\Phi(t^{1/\alpha}/r) \ell (m/t^{1/\alpha}) \int_{k^{\alpha}}^t \left(\frac{k}{s^{1/\alpha}}\right)^{q} \bigg( \frac{s}{t}\bigg)^{(\beta-2\eps)/\alpha} ds \le c_2 t  \left(\frac{k}{t^{1/\alpha}}\right)^{q}	\Phi(t^{1/\alpha}/r) \ell (m/t^{1/\alpha}) .
\end{align*}
In the integration we used the fact that $(\beta-2\eps-q)/\alpha >-1$.  This complete the proof.
\qed

\begin{lemma}\label{cal:new1} 
For any $\eps \in (0,1)$, there exists $C=C(\eps)>0$  such that for all $x\in D$ and   $r>0$,
	\begin{align}\label{e:cal:new1} 
		\int_{z\in D: |x-z|<r} \delta_D(z)^{-\eps} dz \le C \left(\delta_D(x)\vee r\right)^{-\eps} r^d.
	\end{align}
\end{lemma}
\pf Since $D$ is a bounded Lipschitz open set,   $\partial D$  has Hausdorff dimension $d-1$. Thus,  by \cite[Lemma 2.1]{Le08}, we  obtain $	\int_{z\in D: |x-z|<r} \delta_D(z)^{-\eps} dz \le c_1 r^{d-\eps}.$ This proves
\eqref{e:cal:new1} for the case $\delta_D(x)<2r$.  If $\delta_D(x)\ge 2r$, then $\delta_D(z)\ge \delta_D(x)-r \ge  \delta_D(x)/2$ for all $|x-z|<r$. Hence,
$$
	\int_{z\in D: |x-z|<r} \delta_D(z)^{-\eps} dz \le (\delta_D(x)/2)^{-\eps}	\int_{B(x,r)} dz = c_2 \delta_D(x)^{-\eps}r^{d}.
$$
\qed

\begin{lemma}\label{cal:new2}
For any $\eps\in (0,1)$ and $\delta>0$, there exists $C=C(\eps,\delta)>0$  such that for all $x\in D$ and $r>0$,
\begin{align*}
		\int_{z\in D:  |x-z| \ge r} \delta_D(z)^{-\eps} |x-z|^{-d-\delta}\, dz \le C \left(\delta_D(x)\vee r\right)^{-\eps} r^{-\delta}.
\end{align*}
\end{lemma}
\pf   Applying Lemma \ref{cal:new1}, we obtain
\begin{align*}
		&\int_{z\in D:  |x-z| \ge r} \delta_D(z)^{-\eps} |x-z|^{-d-\delta}\, dz \le  \sum_{n=1}^\infty  (2^{n-1}r)^{-d-\delta}	\int_{z\in D:  |x-z| \in [2^{n-1}r, 2^nr)} \delta_D(z)^{-\eps} \, dz \\
		&\le c_1 \sum_{n=1}^{\infty} (2^{n}r)^{-\delta} (\delta_D(x) \vee (2^nr))^{-\eps} \le  c_1r^{-\delta} (\delta_D(x) \vee r)^{-\eps} \sum_{n=1}^\infty 2^{-n\delta} =  c_2r^{-\delta} (\delta_D(x) \vee r)^{-\eps}.
\end{align*}
\qed

For $x\in D$ and $t>0$, let $\delta_D(x,t):=\delta_D(x)\vee t^{1/\alpha}$.

\begin{lemma}\label{l:analog-of-10.10}
  Let $r\in (0, \diam(D))$,  let $x\in D$ satisfy $\delta_D(x)<5r$, and  $B:=B_D(x,r)$.

  \noindent (i) Assume that $\gamma^*<\alpha+\beta$ and  
that $\ell$ is almost increasing on $(0,1]$	if $\gamma=0$.  Then  there exists $C>0$ independent of $r$ and $x$ such that for all $0<t\le r^\alpha$ and $k>0$, 
\begin{align*}
	\begin{split} 
&\int_{B} \left(t^{-d/\alpha}\wedge \frac{t}{|x-z|^{d+\alpha}}\right) \Phi\bigg(\frac{\delta_D(x,t)}{|x-z|}\bigg) \Psi\left(\frac{\delta_D(z)}{r}\right) \ell\left(\frac{\delta_D(z)}{k}\right) dz\\
& \le C \Psi\left(\frac{\delta_D(x,t)}{r}\right)\ell\left(\frac{\delta_D(x,t)}{k}\right).
\end{split} 
\end{align*}
 
\noindent (ii) Assume that $\ell$ is almost decreasing on $(0,1]$	if $\gamma=0$. Then  there exists $C>0$ independent of $r$ and $x$ such that for all $0<t\le r^\alpha$ and $k>0$, 
\begin{align*}
	\begin{split} 
		&\int_{B} \left(t^{-d/\alpha}\wedge \frac{t}{|x-z|^{d+\alpha}}\right) \Phi\bigg(\frac{\delta_D(x,t)}{|x-z|}\bigg) \Psi\left(\frac{\delta_D(z)}{r}\right) \ell\left(\frac{k}{\delta_D(z)}\right) dz\\
		& \le C\left[ \Psi\left(\frac{\delta_D(x,t)}{r}\right)\ell\left(\frac{k}{\delta_D(x,t)}\right)
		+	t\int_{\delta_D(x,t)\wedge r}^{r}  
		\Phi\bigg(\frac{\delta_D(x,t)}{l} \bigg)\Psi\bigg(\frac{l}{r}\bigg)\ell\bigg(\frac{k}{l}\bigg)\frac{dl}{l^{\alpha+1}} \right] .
	\end{split} 
\end{align*}
\end{lemma}
\pf  We first give the proof of (ii). We first consider the integral over $|x-z|<t^{1/\alpha}$. In this case, we have that
$\delta_D(z)\le \delta_D(x)+ |x-z| \le 
2\delta_D(x,t)$. Using the boundedness of $\Phi$ in the first inequality, 
the almost increasing and scaling properties of $\Psi$ and \eqref{e:ell-aux-scaling} (with $R=5$, $\eta=1/2$ and  $\eps=1/2$)
in the second, and Lemma \ref{cal:new1}  in the third, we obtain
\begin{align*}
	\begin{split}
	&\int_{z\in B: |x-z|<t^{1/\alpha}} \left(t^{-d/\alpha}\wedge \frac{t}{|x-z|^{d+\alpha}}\right)\Phi\bigg(\frac{\delta_D(x,t)}{|x-z|}\bigg) \Psi\left(\frac{\delta_D(z)}{r}\right) \ell\left(\frac{k}{\delta_D(z)}\right) dz\\
	&\le c_1 t^{-d/\alpha}\int_{z\in B:|x-z|<t^{1/\alpha}}\Psi\left(\frac{\delta_D(z)}{r}\right)   \ell\left(\frac{\delta_D(z)}{k}\right) 
	\left(\frac{\delta_D(z)}{r}\right)^{1/2}  \left(\frac{\delta_D(z)}{r}\right)^{-1/2}\, dz \\
	&\le c_2 t^{-d/\alpha}\left(\frac{\delta_D(x,t)}{r}\!\right)^{1/2} \Psi\bigg(\frac{\delta_D(x,t)}{r}\bigg) \ell\left(\frac{k}{\delta_D(x,t)}\right)\int_{z\in D:|x-z|<t^{1/\alpha}} \left(\frac{\delta_D(z)}{r}\right)^{-1/2} dz\\
	&\le c_3\Psi\bigg(\frac{\delta_D(x,t)}{r}\bigg)\ell\left(\frac{k}{\delta_D(x,t)}\right).
\end{split}
\end{align*}
When $t^{1/\alpha}\le |x-z|<\delta_D(x)$, we have that  $\delta_D(z)\le \delta_D(x)+ |x-z|< 2\delta_D(x)= 2\delta_D(x,t)$.  Hence, by the similar argument as above and by using Lemma \ref{cal:new2} in the last inequality, we get
\begin{align*}
	\begin{split} 
	&\int_{z\in B: t^{1/\alpha}\le |x-z|<\delta_D(x)} \left(t^{-d/\alpha}\wedge \frac{t}{|x-z|^{d+\alpha}}\right)\Phi\bigg(\frac{\delta_D(x,t)}{|x-z|}\bigg) \Psi\left(\frac{\delta_D(z)}{r}\right) \ell\left(\frac{k}{\delta_D(z)}\right) dz\\
	&\le c_4 t \left(\frac{\delta_D(x,t)}{r}\right)^{1/2}\Psi\left(\frac{\delta_D(x,t)}{r}\right)\ell\left(\frac{k}{\delta_D(x,t)}\right) \int_{z\in D: t^{1/\alpha}\le |x-z|<\delta_D(x)} \left(\frac{\delta_D(z)}{r}\right)^{-1/2}\frac{dz}{|x-z|^{d+\alpha}}\\
	&\le c_5 \Psi\left(\frac{\delta_D(x,t)}{r}\right)\ell\left(\frac{k}{\delta_D(x,t)}\right).
	\end{split}
\end{align*}

 We now consider the integration over $|x-z|\ge \delta_D(x,t)$. Since $\ell$ is almost decreasing if $\gamma=0$, the function $\Psi(u/r)\ell(k/r)$ is almost increasing on $r$. Thus,  for all $z\in B$ with $|x-z|\ge \delta_D(x,t)$, using $\delta_D(z) \le \delta_D(x) + |x-z| \le 2|x-z|$, \eqref{e:Psi-aux-scaling} and \eqref{e:ell-aux-scaling}, we see that
\begin{align*} 
&\Phi\bigg(\frac{\delta_D(x,t)}{|x-z|}\bigg) \Psi\left(\frac{\delta_D(z)}{r}\right) \ell\left(\frac{k}{\delta_D(z)}\right)\le c_6\Phi\bigg(\frac{\delta_D(x,t)}{|x-z|}\bigg) \Psi\left(\frac{|x-z|}{r}\right) \ell\left(\frac{k}{|x-z|}\right).
\end{align*} 
It follows that 
\begin{align*}
	&\int_{z\in B: |x-z| \ge \delta_D(x,t)} \left(t^{-d/\alpha}\wedge \frac{t}{|x-z|^{d+\alpha}}\right)\Phi\bigg(\frac{\delta_D(x,t)}{|x-z|}\bigg) \Psi\left(\frac{\delta_D(z)}{r}\right) \ell\left(\frac{k}{\delta_D(z)}\right) dz\\
	&\le  c_6t \int_{z\in D:  \delta_D(x,t) \le |x-z| \le r}\Phi\bigg(\frac{\delta_D(x,t)}{|x-z|}\bigg) \Psi\left(\frac{|x-z|}{r}\right) \ell\left(\frac{k}{|x-z|}\right) \frac{dz}{|x-z|^{d+\alpha}}\\
	&=c_7 t \int_{\delta_D(x,t) \wedge r}^{r}  \Phi\bigg(\frac{\delta_D(x,t)}{l} \bigg)\Psi\bigg(\frac{l}{r}\bigg)\ell\bigg(\frac{k}{l}\bigg)\frac{dl}{l^{\alpha+1}}.
\end{align*} 
This complete the proof for (ii).

For (i), by following the above arguments, we can deduce that 
\begin{align}\label{e:analog-of-10.10-1}
	\begin{split} 
		&\int_{B} \left(t^{-d/\alpha}\wedge \frac{t}{|x-z|^{d+\alpha}}\right) \Phi\bigg(\frac{\delta_D(x,t)}{|x-z|}\bigg) \Psi\left(\frac{\delta_D(z)}{r}\right) \ell\left(\frac{\delta_D(z)}{k}\right) dz\\
		& \le c_8\left[ \Psi\left(\frac{\delta_D(x,t)}{r}\right)\ell\left(\frac{\delta_D(x,t)}{k}\right)+	t\int_{\delta_D(x,t)\wedge r}^{r}  \Phi\bigg(\frac{\delta_D(x,t)}{l} \bigg)\Psi\bigg(\frac{l}{r}\bigg)\ell\bigg(\frac{l}{k}\bigg)\frac{dl}{l^{\alpha+1}} \right] .
	\end{split} 
\end{align}
Suppose $\gamma^*<\alpha+\beta$ and let $\eps:=(\alpha+\beta-\gamma^*)/4>0$. Using \eqref{e:Phi-aux-scaling}--\eqref{e:ell-aux-scaling}, we see that for all $\delta_D(x,t)\le l \le r$,
\begin{align*}
	\Phi\bigg(\frac{\delta_D(x,t)}{l} \bigg)\Psi\bigg(\frac{l}{r}\bigg)\ell\bigg(\frac{l}{k}\bigg) & \le c_9	\Phi(1)\Psi\bigg(\frac{\delta_D(x,t)}{r}\bigg)\ell\bigg(\frac{\delta_D(x,t)}{k}\bigg)  \bigg(\frac{l}{\delta_D(x,t)}\bigg)^{-\beta+\gamma^*+3\eps}\\
	&= c_9\Psi\bigg(\frac{\delta_D(x,t)}{r}\bigg)\ell\bigg(\frac{\delta_D(x,t)}{k}\bigg)  \bigg(\frac{l}{\delta_D(x,t)}\bigg)^{\alpha-\eps}.
\end{align*}
It follows that
\begin{align*}
&t\int_{\delta_D(x,t)\wedge r}^{r}  \Phi\bigg(\frac{\delta_D(x,t)}{l} \bigg)\Psi\bigg(\frac{l}{r}\bigg)\ell\bigg(\frac{l}{k}\bigg)\frac{dl}{l^{\alpha+1}}\\
 &\le  \frac{c_9t}{\delta_D(x,t)^{\alpha-\eps}}\Psi\bigg(\frac{\delta_D(x,t)}{r}\bigg)\ell\bigg(\frac{\delta_D(x,t)}{k}\bigg) \int_{\delta_D(x,t)}^r \frac{dl}{l^{1+\eps}}\\
	&\le \frac{c_9t}{\eps\delta_D(x,t)^{\alpha}}\Psi\bigg(\frac{\delta_D(x,t)}{r}\bigg)\ell\bigg(\frac{\delta_D(x,t)}{k}\bigg)\le \frac{c_9}{\eps}\Psi\bigg(\frac{\delta_D(x,t)}{r}\bigg)\ell\bigg(\frac{\delta_D(x,t)}{k}\bigg).
\end{align*}
Combining this with \eqref{e:analog-of-10.10-1}, we obtain the desired result.
\qed

\textbf{Acknowledgement:} Part of this work was done while Zoran Vondra\v{c}ek was visiting Seoul National University. The support by the National Research Foundation of Korea (NRF) grant funded by the Korea government (MSIP) (No. 2022H1D3A2A01080536) is gratefully acknowledged.

\vspace{.1in}


\end{document}